\font\twelverm=cmr12
\font\twelveit=cmti12
\font\twelvett=cmtt12
\font\twelvesl=cmsl12 
\font\twelvebx=cmbx12 
\font\erm=cmmi12 
\font\esy=cmsy10 scaled 1200 
\font\eex=cmex10 scaled 1200
\font\symbol=msam10
\font\field=msbm10
\textfont\bffam=\twelvebx
\textfont\itfam=\twelveit
\textfont\slfam=\twelvesl
\textfont\ttfam=\twelvett
\textfont1\erm
\textfont2\esy
\textfont3\eex
\font\chapter=cmbx10 scaled 2200
\font\section=cmbx10 scaled 1600
\font\title=cmr10 scaled 1600
\def\bigjump{\vskip 1.27cm}
\def\medjump{\vskip 0.8cm}
\def\X{${\bf X}=(X,{\cal B},\mu,T)$\ } 
\def\Y{${\bf Y}=(Y,{\cal C},\nu,S)$\ } 
\def\proclaim #1. #2\par{\medbreak\noindent{\twelvebx#1.\enspace}{\twelvesl#2}\par\medbreak}
\def\joinXY{$J^{+}(\bf X, \bf Y)$}

\def\Bp{$B^{+}(\bf {p})$}

\def\closeproof{ \symbol \char'04 }
\def\whole{{\bf N}^{\ast}}
\baselineskip=20pt
\nopagenumbers
\def\leaderfill{\leaders\hbox to 1em{\hss. \hss}\hfill}
\twelverm
$ $
\vskip 6cm
\centerline{\title On the Isomorphism Problem of {\section p}-Endomorphisms}
\vskip 4cm 
\centerline{By}
\vskip 2cm 
\centerline{Peter Jong}
\vskip 4cm
\centerline {A thesis submitted in conformity with the requirements}
\centerline {for the degree of Doctor of Philosophy,} 
\centerline {Graduate Department of Mathematics, University of Toronto}
\medjump
\centerline {\copyright Peter Jong 2003} 
\vfill \eject
\footline={\hss\tenrm\folio\hss}
\pageno=-2
\centerline {\section Abstract}
\centerline {On the Isomorphism Problem of {\twelvebx p}-Endomorphisms}
\centerline {Peter Jong, Ph.D.}
\centerline {Department of Mathematics, University of Toronto, 2003}
\medskip
Let \X be a measure-preserving system on a Lebesgue probability space.  Given a fixed probability vector ${\bf p}=(p_1,\ldots, p_s)$, we say that \X is a $\bf p$-{\twelvebx endomorphism} if $T$ is $s$-to-1 a.e. and the conditional probabilities of the preimages are precisely the components of ${\bf p}$. Two measure-preserving systems \X and \Y are {\twelvebx isomorphic} if there exist a measure-preserving bijective map $\varphi : X \to Y$ such that $\varphi T = S \varphi$ a.e.   

This thesis considers the isomorphism problem of $\bf p$-endomorphisms, generalizing the work of Hoffman and Rudolph [H,R] which treats the case when $\bf p$ is a uniform probability vector, i.e. ${\bf p}=(\hbox{1}/p,\ldots,\hbox{1}/p)$.  In particular, we generalize the tvwB criterion introduced in Hoffman and Rudolph to prove two results.  

The first result is Theorem 2.4.1, which generalizes the main theorem in [H,R] to $\bf p$-endomorphisms.  We paraphrase this as follows:

\proclaim Theorem 2.4.1. Let \X be a $\bf p$-endomorphism.  Then \X is one-sided Bernoulli if and only if \X is tvwB. \par

We give two proofs of this result.  The first follows Ornstein's classical proof of his famous theorem that two shifts of equal entropy are isomorphic, and a second proof which follows the joinings proof as given in [H,R].  As a corollary of the joinings proof, we show that there are uncountably many automorphisms of the one-sided Bernoulli shift $B^+ ({\bf p})$ unless the components of $\bf p$ are pairwise distinct.  We also give examples of tvwB $\bf p$-endomorphisms such as mixing one-sided Markov shifts and a generalization of the $[T,Id]$ transformation.

The second main result is Theorem 5.1.1, which in view of Theorem 2.4.1, reduces to the statement that for any two tvwB finite group extensions of one-sided Bernoulli shifts, there is an isomorphism between them in a stronger sense than that asserted in Theorem 2.4.1.  Specifically, we have the following theorem in Chapter 5 which we paraphrase as follows:  

\proclaim Theorem 5.1.1$'$. Let $G$ be a finite group. For any two tvwB $G$-extensions of the one-sided shift $B^+ ({\bf p})$, there is an isomorphism which preserves the Bernoulli factor algebra and maps fibres over points in the factor to other such fibres by group rotations. \par
\vfill \eject 

\centerline {\section Acknowledgements}
\twelverm
\medjump
First and foremost, I would like to thank my thesis advisor, Professor Andr\'es del Junco, for his insights, advice and encouragement throughout the program. Without his assistance and reassurance (especially when I have self-doubts), this thesis would not have been possible.  I also wish to thank Professor Dan Rudolph for some helpful conversations during conferences in Memphis and Toronto, and for his groundbreaking work with Christopher Hoffman on which this thesis is based.  I also appreciate the numerous helpful suggestions from Professors George Elliott, Michael Yampolsky, Chandler Davis and Misha Lyubich, who also served on my final examination committee.    

In addition, I would like to thank NSERC and the Department of Mathematics in the University of Toronto for their generous financial support over the past few years. I owe special thanks to Ida Bulat for overseeing all the administrative aspects of the program.

Finally, I wish to thank my junior high school mathematics teacher, Charles Ledger, my high school mathematics teachers, Bill Bisset and Robert Velazquez, as well as my older brother, Philip, for instilling an interest in mathematics in me during my adolescent years.  Mathematics is such a beautiful (albeit challenging) subject that I feel fortunate to have the opportunity to study it in my lifetime.    
\vfill \eject
$ $
\vskip 1.50cm
\centerline{\chapter TABLE OF CONTENTS}
\bigskip
\line{\section 1.  Introduction \leaderfill 1}
\line{\itemitem{1.1} Background\leaderfill 1}
\line{\itemitem{1.2} Organization and Contributions of Thesis \leaderfill	3}
\line{\itemitem{1.3} The Tree Very Weak Bernoulli Condition \leaderfill	4}
\medskip
\line{\section 2.  TvwB p-Endomorphisms on Lebesgue spaces \leaderfill	 8}
\line{\itemitem {2.1} Tree-Adapted Factors \leaderfill	9}
\line{\itemitem {2.2} The Tree Ergodic Theorem \leaderfill	11}
\line{\itemitem {2.3} The Tree Rokhlin Lemma \leaderfill	12}
\line{\itemitem {2.4} Proof of the Isomorphism Theorem \leaderfill	15}
\line{\itemitem {2.5} An Elementary Proof of a Special Case of Theorem 2.4.1	 \leaderfill 30}
\medskip
\line{\section 3.  A Joinings Proof of the Isomorphism Theorem \leaderfill 33} 
\line{\itemitem{3.1}One-Sided Joinings \leaderfill	34}
\line{\itemitem {3.2} The Copying Lemma \leaderfill	38}
\line{\itemitem{3.3} The \=t    Distance \leaderfill	43}
\line{\itemitem{3.4} A Joinings Proof of Theorem 2.4.1 \leaderfill	48}
\medskip
\line{\section 4.  Examples of TvwB p-Endomorphisms \leaderfill	51}    
\line{\itemitem{4.1}  One-sided Markov Shifts \leaderfill	51}
\line{\itemitem{4.2} A Generalization of $[T,Id]$ \leaderfill	56}
\medskip
\line{\section 5.  Finite Group Extensions of One-sided Bernoulli Shifts \leaderfill 59}
\line{\itemitem{5.1} An Isomorphism Theorem on TvwB $G$-Extensions \leaderfill 59}
\line{\itemitem{5.2} Some Applications of Theorem 5.1.1 \leaderfill 67}      
\medskip
\line {\section 6.  Some Open Problems	\leaderfill 70}
\medskip
\line{\section References \leaderfill 72}
\vfill \eject
\footline={\hss\tenrm\folio\hss}
\pageno=1
$ $
\vskip 5.90cm
\noindent {\chapter Chapter 1:} 
\medjump
\noindent {\chapter Introduction}
\bigjump
\twelverm
\noindent {\section 1.1. Background}
\medjump

Let $X$ be a compact metric space and ${\cal B}$ be its Borel $\sigma$-algebra. Suppose $(X,{\cal B},\mu )$ is a {\twelvebx nonatomic Lebesgue probability space} (a probability measure space isomorphic to 
the unit interval with Lebesgue measure). An {\twelvebx endomorphism} of $X$ is a measure-preserving map 
$T:X \to X$, i.e. for all $B \in {\cal B}$, $\mu (T^{-1} B)=\mu (B)$. 
For us, a {\twelvebx measure-preserving system (m.p.s.)} is a quadruple $(X,{\cal B},\mu ,T)$ where $(X,{\cal B},\mu )$ is a nonatomic Lebesgue probability space defined on a compact metric space $X$, and $T$ is an endomorphism. 

This thesis is mainly concerned with the isomorphism problem of measure-preserving systems. Two measure-preserving systems \X and \Y are ({\twelvebx measure-theoretically}) {\twelvebx isomorphic} if after deleting null sets 
$X_{0}$ and $Y_{0}$ from $X$ and $Y$ respectively, there is a bijection $\phi 
:X\backslash X_{0} \to Y\backslash Y_{0}$ such that $\phi $ is 
measure-preserving and $\phi T=S\phi $ on $X\backslash X_{0}$. We say 
that $\phi $ is an {\twelvebx isomorphism} from ${\bf X}$ to ${\bf Y}$. In case ${\bf X}={\bf Y}$, we 
say that $\phi $ is an {\twelvebx automorphism} of $\bf X$. More generally, a {\twelvebx factor map} $\phi :{\bf X} \to 
{\bf Y}$ is a measure-preserving map $\phi :X \to Y$ such that $\phi 
T=S\phi$ a.e. In this case, we say that $\bf  Y$ is a {\twelvebx factor} of $\bf X$, and $\bf X$ is an {\twelvebx extension} of $\bf Y$.

The classical isomorphism problem in ergodic theory is the classification of {\twelvebx two-sided Bernoulli shifts}. To define a two-sided Bernoulli shift, fix a probability vector ${\bf 
p}=(p_{1},\ldots ,p_{s})$. Consider the finite set $I=\{ \hbox{1},\ldots ,s \}$ 
and a measure ${\bf m}$ on $I$ defined by ${\bf m}(j)=p_{j}$. Construct the 
product space $I^{\bf Z}$. Let ${\cal B}$ denote the Borel sigma-algebra and $\sigma $ denote the left-shift 
transformation defined by $\sigma (x)_j=x_{j+1}$, where $x_t$ is the $t$-th 
coordinate of $x$. The two-sided Bernoulli shift $B({\bf p})$ is 
the m.p.s. $(I^{\bf Z},{\cal B},{ \bf m}^{ \bf Z},\sigma )$.

The breakthrough in the classification problem of two-sided Bernoulli shifts 
came when Kolmogorov introduced the concept of {\twelvebx entropy} into ergodic theory. It is not difficult to show that the entropy of a m.p.s. is invariant under 
isomorphism. Moreover, the entropy of a Bernoulli shift is easy to compute. Indeed, for 
a probability vector ${\bf p}=(p_{1},\ldots ,p_{s})$, define the {\twelvebx entropy} of $\bf p$, denoted $h({\bf p})$, to be $\sum_{j=1}^s  - p_{j}\hbox{log}_{2}p_{j}$. It can 
be shown that the entropy of the two-sided Bernoulli shift $B({\bf p})$ is 
just $h({\bf p})$. It thus follows that the two-sided Bernoulli shift on two 
symbols with equal weights, which has entropy $\hbox{log}_{2}\hbox{2}$, is not isomorphic to 
the two-sided Bernoulli shift on three symbols with equal weights, which has 
entropy $\hbox{log}_{2}\hbox{3}$. 

The solution of the classification problem was finally achieved around 1970 
when Ornstein [Orn] showed that entropy is in fact a {\twelvebx complete invariant} for 2-sided 
Bernoulli shifts; that is, two 2-sided Bernoulli shifts are isomorphic if 
and only if they have the same entropy. More generally, it can be shown that 
entropy is also a complete invariant for two-sided shift spaces with {\twelvebx  finitely determined measures}. Since Ornstein's original proof, other criteria such as the {\twelvebx weak Bernoulli} and the {\twelvebx very weak Bernoulli} conditions 
have been developed which also turn out to be sufficient for a 
m.p.s. to be isomorphic to a two-sided Bernoulli shift. The reader is referred 
to Shields [Shi] for an excellent account of the proof of Ornstein's theorem 
as well as a discussion of finitely determined measures. Examples of weak 
Bernoulli and very weak Bernoulli systems such as ergodic toral automorphisms 
and two-sided mixing Markov shifts are discussed in Petersen [Pet]. 

What we will deal with in this thesis is the isomorphism problem in the case 
when the endomorphism $T:X \to X$ is not invertible. Throughout this thesis, ${\bf N}^{\ast}$ will denote the set of nonnegative integers. One example of a 
non-invertible endomorphism is a {\twelvebx one-sided Bernoulli shift}, which is derived from a two-sided Bernoulli shift by restricting the shift space to $\{\hbox{1},\ldots 
,s\}^{{\bf N}^{\ast}}$.  Let $B^{ + }({\bf p})$ denote the m.p.s. obtained from $B({\bf p})$
in this way. In this connection, we remark that 
Ornstein's theorem does not hold for one-sided Bernoulli shifts. Indeed, it 
is easy to construct two probability vectors with different numbers of 
components which have the same entropy. This leads to two 1-sided shifts 
having the same entropy; however, they are not isomorphic because they have 
distinct numbers of inverse images.
\vfill \eject
In a recent paper, Hoffman and Rudolph made a fundamental contribution to the 
isomorphism problem of measure-preserving systems with non-invertible maps. In [H,R], Hoffman and 
Rudolph considered a class of endomorphisms called the {\twelvebx uniformly p-to-1 endomorphisms}, and introduced a 
condition called {\twelvebx tree very weak Bernoulli }({\twelvebx tvwB}) on this class. They then proved that the tvwB condition is 
necessary and sufficient for a uniformly $p$-to-1 endomorphism to be isomorphic 
to the one-sided Bernoulli shift $B^{ + }({\bf p})$, where ${\bf p}$ is the 
uniform vector $(\hbox{1}/p,\ldots ,\hbox{1}/p)$. As the authors remarked, entropy 
turns out to have no role in this theory.  Perhaps as a result of this, many of the arguments are much simpler than those in Ornstein's proof. This thesis will consider 
generalizations of the tvwB criterion which will allow us to prove 
isomorphism theorems for more general measure-preserving systems with non-invertible 
transformations.
\medjump
\noindent{\section 1.2. Organization and Contributions of Thesis} 
\medjump
Without delving into definitions and details, which will be presented after 
this section and in subsequent chapters, we outline below the contents and 
main results in each chapter of this thesis.

Chapter 1 is this introduction. In \S 1.3, we will define the class of 
objects of interest in this thesis, the ${\bf p}$-endomorphisms, and extend 
the tvwB criterion to a general probability vector ${\bf p}$ (i.e. ${\bf 
p}$ not necessarily uniform).

Chapters 2 and 3 give two proofs of one of the main results in this thesis, 
Theorem 2.4.1, which generalizes the main result in [H,R] to ${\bf 
p}$-endomorphisms for a general probability vector ${\bf p}$:
\proclaim Theorem 2.4.1. Let ${\bf p}=(p_{1},\ldots ,p_{s})$ be a probability vector. If \X is a ${\bf p}$-endomorphism, then $\bf X$ is tvwB if and only if ${\bf X} \cong B^{ + }({\bf p})$. \par

In Chapter 2, we give a proof which follows the classical argument used by Ornstein in his isomorphism theorem. The 
proof is based on a joint paper by del Junco and me [J,J] which treats the 
case when ${\bf p}$ is a uniform probability vector. In particular, it does 
not require the machinery of joinings introduced in [H,R]. The main 
additional ingredient which enables us to extend the proof in [J,J] is 
proposition 2.1.3, which essentially says that under certain factor maps 
(tree-adapted factor maps) between ${\bf p}$-endomorphisms, conditional 
probabilities of inverse images are preserved.  In \S 2.5, we give an entirely different and quite elementary proof of Theorem 2.4.1 in the special case when 
the components of ${\bf p}$ are pairwise distinct. 

In Chapter 3, we will mirror the proof of the main result in [H,R] and extend their notion of {\twelvebx
one-sided joinings} to give a proof of Theorem 2.4.1. The new ingredient 
needed is an additional condition in the definition of one-sided joinings 
which is trivial in the case of a uniform probability vector. We will also 
show that the joinings proof implies that there are uncountably many 
automorphisms of $B^{ + }({\bf p})$, unless the components of ${\bf p}$ are 
pairwise distinct (in which case the identity is the only one). 

Chapter 4 illustrates some examples of tvwB ${\bf p}$-endomorphisms such 
as certain classes of one-sided Markov shifts and certain extensions of 
one-sided Bernoulli shifts.  It follows that these are all one-sided Bernoulli by Theorem 
2.4.1. 

Chapter 5 proves an isomorphism theorem for tvwB finite group extensions of one-sided 
Bernoulli shifts.  The main result is Theorem 5.1.1, which in view of Theorem 2.4.1 reduces to the following result which we paraphrase as follows:
\proclaim Theorem 5.1.1$'$. For any two tvwB finite group extensions of the one-sided shift $B^+ ({\bf p})$, there exists an isomorphism which preserves the Bernoulli factor algebra and maps fibres over points in the factor to other such fibres by group rotations. \par

\medjump
\noindent {\section 1.3. The Tree Very Weak Bernoulli Condition}
\medjump
Suppose that $(X,{\cal B},\mu )$ is a nonatomic Lebesgue space and $T:X \to X$ 
is measure-preserving. Fix a probability vector ${\bf p}$ with finitely many components. Let $\left|{\bf p}\right|$ denote the number of components of ${\bf p}$. A 
m.p.s. \X is a {\twelvebx $\bf p$-endomorphism} if $T$ is $\left|{\bf p}\right|$ to 
1 a.e., the conditional probabilities of the $\left|{\bf p}\right|$ inverse images of $x$ are 
the components of ${\bf p}$ for a.a. $x$, and the entropy of $\bf X$ is $h({\bf p})$.  Throughout 
this thesis, we will let $End({\bf  p})$ denote the collection of ${\bf 
p}$-endomorphisms.

The standard element in $End({\bf p})$ is the one-sided Bernoulli 
shift $B^{ + }({\bf p})$. Let us recall its definition here. For 
the finite set $I=\{ \hbox{1},\ldots ,\left| {\bf p} \right| \}$, define a measure ${\bf m}$ on $I$ 
by ${\bf m}(j)=p_{j}$. Consider the product space $I^{{\bf N} ^{\ast} }$ with the Borel sigma-algebra $\cal B$, product measure ${\bf 
m}^{{\bf N} ^{\ast} }$ and shift transformation $\sigma $. We then define $B^{ 
+ }({\bf p})$ to be the m.p.s. $(I^{{\bf N}^{\ast }},{\cal B},{\bf 
m}^{{\bf N} ^{\ast} } ,\sigma )$.  (The entropy condition in the definition of $\bf p$-endomorphism is thus a natural one to make since $B^+ ({\bf p})$ has entropy $h({\bf p})$.) 

In [H,R], Hoffman and Rudolph introduced a condition called tree very weak 
Bernoulli (tvwB) which they proved to be necessary and sufficient for 
a {\twelvebx uniformly p-to-one endomorphism}, i.e. a ${\bf p}$-endomorphism with ${\bf p}=({\hbox{1} / p},\ldots 
,{\hbox{1}/ p})$, to be isomorphic to the one-sided Bernoulli shift 
on $p$ symbols with equal weights. We will now extend the various definitions in 
[H,R] in order to handle the case that the components of the probability 
vector are not all equal. 

For a probability vector ${\bf p}$, we define the {\twelvebx ${\bf p}$-tree} to be the set, 
denoted ${\cal T}$, consisting of all finite sequences (including the empty 
sequence) of integers in $\{\hbox{1},\ldots , \left| {\bf p} \right| \}$. We define a {\twelvebx node} to be an 
element of the ${\bf p}$-tree. The {\twelvebx length} of a node $v$, denoted $\left| v \right|$, is the number of 
integers in the sequence $v$ (thus, the empty sequence has length zero). We 
will refer to the empty sequence as the {\twelvebx root node} and denote it by $\emptyset$. Given any 
two nodes $u$ and $v$, we define the node $uv$ by concatenating the sequence $u$ to 
the left of the sequence $v$. We will use ${\cal T}'$ to denote the set ${\cal T}\backslash\emptyset$.  Define the map $\sigma : {\cal T}'  \to {\cal T}$ by setting $\sigma (v)$ to be the sequence obtained by 
deleting the leftmost symbol in the sequence $v$. Note that if $\left| v \right| =\hbox{1}$, then $\sigma (v)=\hbox{\char'10}$. Moreover, $\sigma $ is a $\left|{\bf p}\right|$ to 1 
surjection. It is helpful to picture a ${\bf p}$-tree as a tree in the graph-theoretic sense with vertices 
corresponding to the nodes and  with an edge between $u$ and $\sigma (u)$ for each $u \ne \emptyset$. 
The picture we get is the usual $\left|{\bf p}\right|$-ary tree.  In view of this pictorial representation, we refer to the set $\{uv \mid \hbox{0}\le \left| u \right| \le N \}$ as the {\twelvebx subtree of height $N$ rooted at $v$} and the set $\{uv \mid \left| u \right| \ge \hbox{0} \}$ as the {\twelvebx subtree rooted at $v$}.

To each node $v$ in ${\cal T}$, we assign a weight, denoted $w_{v}$, as 
follows. For the root node ${\emptyset}$, we set $w_{\emptyset}=\hbox{1}$. For any other node $v=(a_{1},\ldots ,a_{j})$, set 
$w_v = \prod\nolimits_{i = 1}^j {p_{a_i } } $. We define a {\twelvebx tree automorphism } to be a 
bijection $A:{\cal T} \to {\cal T}$ such that $A \circ \sigma 
(v)=\sigma  \circ A(v)$ and 
$w_{v}=w_{Av}$ for $v \in {\cal T}'$ (i.e. $A$ preserves the tree structure and weights).  Note that this implies $A(\emptyset)=\emptyset$. It is obvious that the set ${\cal A}$ of tree automorphisms 
forms a group under composition. For $N \in {\bf N}$, let ${\cal T}_N \subseteq {\cal T}$ denote the set of nodes of length $\le N$ and let ${\cal T}_N' $ denote the set ${\cal T}_N \backslash \emptyset$.  Let ${\cal A}_{N}$ be the subgroup of bijections of ${\cal T}_{N}$ which is ${\cal A}|_{{\cal T}_N}$. 

Given a fixed compact metric space $(R,d)$, let us say that a {\twelvebx $R$-tree name} is a function $g:{\cal T'} \to R$, and for each $N \in {\bf N}$, a {\twelvebx $R$,$N$-tree name} is a function $g':{\cal T}_N' \to R$. As in [H,R], we define a distance function, $\bar 
{t}_N $, on the space of $R$,$N$-tree names as follows: if $h:{\cal 
T}_{N}'  \to R$ and $h':{\cal 
T}_{N}' \to R$, let 
\vfill \eject
$$
{\bar t}_N {(h,h')} = {\hbox{1} \over N} \inf \limits_{A \in {\cal A}_N} \sum\limits_{0 < \left| v \right|  \le N } {d( 
h(v),h'(Av))w_v}.
$$

Note that ${\bar t}_N $ is not a metric but it does satisfy the triangle 
inequality. We say that two $R$,$N$-tree names $h$ and $h'$ are the same {\twelvebx up to tree automorphism} if 
there exists some tree automorphism $A \in {\cal A}_{N}$ such that 
$h(v)=h'(Av)$ for all $v$.  Note that $\bar t_N{(h,h')}=\hbox{0}$ if and only if $h$ and $h'$ are the same up to tree automorphism. We shall let $R^{N\tau }$ denote the set of $R$,$N$-tree names 
and $R^{N\nabla }$ denote the equivalence classes of $R$,$N$-tree names modulo 
tree automorphism. 

Suppose that \X $\in End({\bf p})$. We wish to 
construct $R$-tree names for each point in $X$. To do this,
fix a measurable partition $K_{\bf X}:X \to \{ \hbox{1},\ldots ,\left|{\bf p}\right|\}$ 
such that for a.a. $x$, the $\left|{\bf p}\right|$ inverse images of $x$ have distinct $K_{\bf X}$ 
values and the conditional probability of the inverse image $x'$ of $x$ 
with $K_{\bf X}(x')=j$ is $p_{j}$.  We will refer to $K_{\bf X}$ as a {\twelvebx tree partition} of $\bf X$. Note that unless the components of ${\bf p}$ are pairwise distinct, 
there will in general be many distinct tree partitions.

We now use $K_{\bf X}$ to define a collection of {\twelvebx partial inverses} $T_{v}$ for each node $v \in 
{\cal T}'$. For each node $v$ of 
length one (i.e., $v \in \{\hbox{1},\ldots ,\left|{\bf p}\right|\}$), define the map $T_{v}:X 
\to X$ by setting $T_{v}x$ to be the inverse image $x'$ of $x$ with 
$K_{\bf X}(x')=v$. We may then extend the definition of $T_{v}$ to an 
arbitrary node $v \in {\cal T}'$ as follows: if $v=(a_{1},\ldots ,a_{n})$, 
set $T_v (x) = T_{a_1 } (T_{a_2 }\ldots(T_{a_n } x)\ldots)$. Note that $T_{v}$ is 
injective and $v \mapsto T_{v}x$ maps ${\cal T'}$ to $\{T^{-j}x \mid {j > \hbox{0}} \}$ for a.a. $x$. If $y \in T^{ - 
j}x$ for some $j \in {\bf N}$, then there exists a unique $v$ of length $j$ such 
that $y=T_{v}x$, and $w_{v}$ is the conditional probability of the preimage $y$ of $x$ under the map $T^{j}$. 

Suppose we have some function $g:X \to R$.  Using the tree partition $K_{\bf X}$, we may now associate to 
a.a. $x$ in $X$ the $R$-tree name $\tau _x^g :{\cal T'} \to R$ by 
setting $\tau _x^g (v)=g(T_{v}x)$. We shall refer to the $R$-tree name $\tau _x^g 
$ as the $g$-{\twelvebx tree name }of $x$ and the restriction $\tau _x^g \vert _{{\cal T}_{N}'} $ as the $g$-$N$-{\twelvebx tree name} of $x$. For any function $f:X \to R$ and $n \in {\bf N}$, define $f^{n\tau }:X \to R^{n\tau }$ by sending $x$ to its $f$-$n$-tree name and $f^{n\nabla }:X \to R^{n\nabla }$ by sending $x$ to the equivalence class in $R^{n\nabla }$ 
containing its $f$-$n$-tree name.

Note that for any $g:X \to R$, the $R$-tree name $\tau _x^g 
$ depends on the choice 
of the tree partition $K_{\bf X}$. However, the tvwB condition, which we 
now define for $\bf p$-endomorphisms, is not affected by the choice (since for any $g:X \to R$, different choices yield the same $g$-$N$-tree name of $x$ up to tree automorphism).  Note that this definition is essentially the same as the definition in [H,R].  The only difference is that our definition of $\bar t _n$ is slightly different as the group of tree automorphisms is more restrictive in the case of a general probability vector $\bf p$.
\proclaim Definition 1.3.1. Let \X $\in End({\bf p})$ and $g:X 
\to R$ for a compact metric space $(R,d)$. We say that  $({\bf X},g)$ is {\twelvebx tree very weak Bernoulli} ({\twelvebx tvwB}) if for each $\varepsilon 
>\hbox{0}$, there exists some $N$ such that, whenever $n \ge N$, we have some set $G \subseteq X$ of measure at least 1$-\varepsilon $ with $\bar {t}_n (\tau _x^g ,\tau 
_y^g ) < \varepsilon $ for all $x$ and $y$ in $G$. \par

We shall say that the a m.p.s. ${\bf X} \in End({\bf p})$ is {\twelvebx tvwB} 
if for all measurable functions $g:X\to R$ for a compact metric space $R$, $({\bf X},g)$ is tvwB. 
It is immediate that the tvwB property is preserved 
under isomorphism.
\vfill \eject
$ $
\vskip 5.90cm
\noindent {\chapter Chapter 2:} 
\medjump
\noindent {\chapter TvwB p-Endomorphisms on Lebesgue}
\medjump
\noindent {\chapter Spaces }
\bigjump

\twelverm
The goal in this chapter is to prove the following theorem stated in \S 2.4, which 
generalizes the main result in [H,R] to an arbitrary probability vector.

\proclaim Theorem 2.4.1.  Let ${\bf p}=(p_{1},\ldots ,p_{s})$ be a probability vector. If \X $\in End({\bf p})$, then ${\bf X}$  is tvwB if and only if ${\bf X} \cong  B^{ + }({\bf p})$.\par

This chapter is organized into several sections. The first three sections 
establish some basic tools that will be necessary in the proof of Theorem 
2.4.1. In \S 2.1, we define the concept of a tree-adapted factor map and 
establish a crucial property shared by these factor maps. \S 2.2 proves the 
tree ergodic theorem, which may be considered as a ``backward'' version of 
Birkhoff's ergodic theorem. \S 2.3 establishes the tree Rokhlin lemma and 
its strong form, which are analogues of the standard Rokhlin lemma for 
invertible transformations. \S 2.4 contains the proof of the main theorem. 
The proof in \S 2.4 does not require the machinery of one-sided joinings
introduced in [H,R] and is patterned on the proof of Theorem 2.4.1 given in a joint 
paper by del Junco and me [J,J] in the case that $\bf p$ is a uniform probability vector. \S 2.5 gives a completely different, and quite simple, proof of Theorem 
2.4.1 in the special case that $\bf p$ has pairwise distinct components. 

Throughout this chapter, $\bf p$ will denote the probability vector 
$(p_{1},\ldots ,p_{s})$.
\medjump
\noindent{\section 2.1. Tree-Adapted Factors}
\medjump

\proclaim Definition 2.1.1. Let \X $ \in End({\bf 
p})$, then a m.p.s. \Y is a {\twelvebx tree-adapted factor} of $\bf X$ if there is a factor map $\phi : {\bf X} 
\to {\bf Y}$  such that for $a.a.x$ in $X$, $\phi $ gives a bijection of the inverse 
images of $x$ and those of $\phi (x)$. We will refer to the factor map $\phi : {\bf X} 
\to {\bf Y}$   as a {\twelvebx tree-adapted factor map}.\par
Note that this definition is more restrictive than in [H,R], which requires only that $\phi$ maps $T^{-1}x$ one-to-one {\twelvebx into} $S^{-1}\phi (x)$.  Nonetheless, this definition will suffice since the factor maps which we construct in the proof of the isomorphism theorem (in particular, in Proposition 2.4.9) is tree-adapted in our sense. 

The following proposition shows that a tree-adapted factor of a $\bf p$-endomorphism is itself a $\bf p$-endomorphism.  It is the extension of Lemma 2.3 in [H,R] to $\bf p$-endomorphisms.  For a ${\bf p}$-endomorphism \X, define the 
function $p_{\bf X}:X \to $(0,1) by setting $p_{\bf X}(x)$ to be the conditional probability of the preimage $x$ of $Tx$. We shall refer to $p_{\bf X}$ as 
the {\twelvebx $p$-function on $\bf X$}. 

\proclaim Proposition 2.1.2. A tree-adapted factor of a $\bf p$-endomorphism is a $\bf p$-endomorphism. \par
\noindent {\twelvebx Proof: }Let \X $ \in End({\bf p})$.  Suppose \Y is a tree-adapted factor of $\bf X$ and $\phi : X \to Y$ is a tree-adapted factor map.  By definition, $S$ is also $\left| {\bf p} \right|$-to-one a.e.  Fix a point $y \in Y$, and suppose $y_1,\ldots,y_s$ are the inverse images of $y$.  For each $x \in \phi^{-1}y$, the conditional measure on $T^{-1}x$ given $x$ pushes forward via $\phi$ to a measure on $S^{-1}y$.  The conditional measure on $S^{-1}y$ given $y$ is just an average of these image measures on $S^{-1}y$ over all $x \in \phi^{-1}y$.  As $\phi$ is tree-adapted and ${\bf X} \in End({\bf p})$, these image measures assign the $s$ inverse images of $y$ with measures equal to the components of $\bf p$.  Thus, $(p_{\bf Y} (y_1)),\ldots,p_{\bf Y}(y_s))$ is an average of probability vectors $(p_{\sigma (1)},\ldots ,p_{\sigma (s)})$ for a permutation $\sigma : \{ \hbox{1},\ldots,s \} \to \{ \hbox{1},\ldots,s\}$.  The fact that the function $h(t)=-t \hbox{log}_{2}t$ is strictly concave implies that the entropy of the probability vector $(p_{\bf Y} (y_1)),\ldots,p_{\bf Y} (y_s))$ is at least $h(\bf p)$ with equality holding if and only if this probability vector is some permutation of $\bf p$.  The fact that the entropy of $\bf X$ is $h(\bf p)$ implies that equality must hold for a.a. $y$ and so the entropy of $\bf Y$ is $h({\bf p})$.  But this is precisely saying that \Y is a $\bf p$-endomorphism.  {\symbol \char'04}

The following proposition gives an important property shared by tree-adapted factor 
maps; namely, they preserve the $p$-function.   
\proclaim Proposition 2.1.3. Suppose \X $\in End({\bf p})$. If $\phi :{\bf X} \to {\bf Y}$ is a tree-adapted factor map, then for $a.a.x$ in $X$, 
$p_{\bf X}(x)=p_{\bf Y}(\phi (x))$.\par
\noindent {\twelvebx Proof: } By proposition 2.1.2, \Y $\in End({\bf p})$.  Let $y \in Y$  and let $\mu_{y}$ be the conditional measure on 
$\phi ^{-1}(y)$ given $y$. With no loss of generality, suppose $p_{1} \ge 
\ldots  \ge p_{s}$.  Let $y_1,\ldots ,y_s$ be the inverse images of 
$y$ such that $p_{\bf Y}(y_j)=p_j$. Since $\phi $ is tree-adapted, for each 
$x \in \phi ^{ - 1}(y)$, let $x_{y_j}$ be the unique inverse image of $x$ 
such that $\phi (x_{y_j})=y_j$. Now, $p_{1}=p_{\bf Y}(y_1)$ is an 
average of $p_{\bf X}(x_{y_1})$ over $x \in \phi ^{-1}(y)$. As 
$p_{\bf X}(x_{y_1}) \le p_{1}$, it follows that $p_{\bf Y}(y_1)=p_{1}= 
p_{\bf X}(x_{y_1})$ for $\mu _{y}$-a.a. $x$. Next, note that since $x_{y_1} \ne x_{y_2}$ by tree-adaptedness, $p_{\bf X}(x_{y_2}) \le p_2$ for $\mu_{y}$-a.a. $x$.  Hence, $p_{\bf Y}(y_2)=p_{2}=p_{\bf X}(x_{y_2})$ for $\mu 
_{y}$-a.a. $x$. Inductively, we see that $p_{\bf Y}(y_j)=p_{\bf X}(x_{y_j})$ 
for $\mu_{y}$-$a.a.x$ for each $\hbox{1} \le j \le  s$. As $y$ is arbitrary, there is a set $G$ of full measure in $X$ such that 
for each $x \in G$, if $x'  \in T^{-1}x$, then $p_{\bf X}(x' 
)=p_{\bf Y}(\phi (x'))$. Then $T^{ - 1}G$ is also a set of full 
measure satisfying the statement of the proposition. {\symbol \char'04}
\medbreak
If \X $\in End({\bf p})$ and $x \in X$,  
define the {\twelvebx $p_{\bf X}$-name} of $x$ to be the sequence $(p_{\bf X}(x),p_{\bf X}(Tx),\ldots )$. Proposition 2.1.3 implies that if $\phi :{\bf X} \to {\bf Y}$ is a tree-adapted factor map, then the $p_{\bf X}$-name 
of $x$ equals the $p_{\bf Y}$-name of $\phi (x)$ for a.a. $x$. 

The fact that the $p$-function is preserved under tree-adapted factor maps leads to the following proposition which will be used in the proof of Theorem 2.4.1.  

\proclaim Proposition 2.1.4. Suppose \X $\in End({\bf p})$ and \Y is a tree-adapted factor of $\bf X$.  Consider a tree-adapted factor map $\phi : {\bf X} \to {\bf Y}$.  For any measurable function $f:Y \to R$ and any $m \in {\bf N}$, we have for $a.a.x$,
$$
(f\circ \phi) ^{m\nabla }(x)=f^{m\nabla }(\phi (x)).
$$ \par
\noindent{\twelvebx Proof:} For $x \in X$, define the map $\pi _{x}:{\cal T}_{m}' 
\to {\cal T}_{m}'$ by setting $\pi _{x}(v)=u$ if $\phi (T_{v}x) 
=S_u \phi(x)$ and $\left| v \right| = \left| u \right|$. Note that this map is well-defined and is a bijection as 
$\phi$ is tree-adapted. 
Moreover, since $\phi$ is a factor map, it follows that $\pi 
_{x} \circ \sigma =\sigma \circ \pi _{x}$. By the remarks following proposition 
2.1.3, $p_{\bf X}$-name of $z = p_{\bf Y}$-name of $\phi (z)$ for a.a. $z$ and thus $\pi _{x}$ preserves the 
weights of nodes so that $\pi _{x}$ defines a tree automorphism. It follows 
that for every $v \in {\cal T}_{m}'$,
$$
\tau _x^{f \circ \phi} (v)  = f \circ \phi (T_v x) = f(S_{\pi _x (v)} \phi (x)) 
= \tau _{\phi (x)}^f (\pi _x (v)).
$$ 
{\symbol \char'04}

\medjump
\noindent {\section 2.2. The Tree Ergodic Theorem}
\medjump
The goal in this section is to demonstrate a tree version of the ergodic theorem (proposition 2.2.3).  To begin, notice that the following is a trivial consequence of the definition of tvwB of a m.p.s. ${\bf X} \in End({\bf p})$. 

\proclaim Proposition 2.2.1.  Suppose ${\bf X} \in End(\bf {p})$ is tvwB, then $\bf X$ is ergodic. \par

\noindent {\twelvebx Proof:} Consider the characteristic function $g=\chi _{G}$ for an invariant set $G$. Then $({\bf X},\chi _{G})$ can only be 
tvwB if $G$ has measure zero or one. {\symbol \char'04}
\medbreak
For a m.p.s. \X and a real-valued function $g$ 
on $X$ and $x \in X$, set
$$
A_N (g)(x) = {\hbox{1} \over N}\sum\limits_{i = 0}^{N - 1} {g(T^{i} x)}.
$$
The following proposition is proved using Birkhoff's ergodic theorem. 
\proclaim Proposition 2.2.2. Suppose \X $\in End(\bf{p})$ is ergodic. Let $B \subseteq X$. Given $\varepsilon>\hbox{0}$, there exists $N \in 
\bf N$ such that for all $n \ge N$, we have a set $G$ of measure $>\hbox{1}-\varepsilon $ such that for each $x \in G$, 
$$
\left| {\hbox{1} \over n}\sum\limits_{0 < \left| v \right| \le n} {w_{v}  \chi _{B} (T_v x) -\mu (B)} \right| < \varepsilon.
$$\par
\noindent{\twelvebx Proof:} Note that
$$
{\hbox{1} \over n}\sum\limits_{0 < \left| v \right|  \le n} {w_v 
 \chi _B  (T_v x)} = \sum\limits_{\left| v\right|  = n} {w_v A_n (\chi _B ) (T_v x)}. \eqno(\hbox{1})
$$
Since $\bf X$ is ergodic, we have for $a.a.x$, $A_{n}(\chi _{B})(x) 
\to \mu (B)$ as $n \to \infty $. In particular, for $\delta >\hbox{0}$, there 
exists some $N$ such that for all $n \ge N$ and for all $y$ in a set $G$ of 
measure $ \ge \hbox{1}-\delta $,
$$
\left| A_n (\chi _B )(y) - \mu (B) \right| < \delta. 
\eqno(\hbox{2})
$$
Hence, if $\mu _{x}$ is the conditional measure on $T^{- n}x$ given $x$, then 
$$
\mu _x (G) \ge \hbox{1} - \sqrt \delta 
\eqno(\hbox{3})
$$
\noindent
for all $x$ in a set $G'$ of measure at least $\hbox{1}-\sqrt \delta$. For 
each $x \in G'$, by (2) and (3), we have
$$
\eqalign
 { | \sum\limits_{\left| v\right|  = n} { w_v A_n (\chi _B 
)(T_v  x)} - \mu (B) | &\le | \sum\limits_{T_v x \in G  \atop 
\left| v\right|  = n} {w_v A_n (\chi _B )(T_v x )}- \sum\limits_{T_v x \in G  \atop \left| v\right|  = n}{w_v \mu (B)} | \cr  
& \quad  + | \sum\limits_{T_v x \notin G  \atop 
\left| v\right|  = n} {w_v A_n (\chi _B )(T_v x )}- \sum\limits_{T_v x \notin G  \atop \left| v\right|  = n}{w_v \mu (B)} | \cr
 &< \delta + \hbox{2}\sqrt \delta. \cr}
$$
Since $\delta $ can be chosen to be arbitrarily small, the result now follows from (1). {\symbol \char'04} 
\medjump
Let \X $\in End({\bf p})$ and consider a partition $P:X \to C$ 
for some finite set $C$. Given $M \in {\bf{N}}$ and $x \in X$, define a 
measure $\theta _{x,M,P}$ on $C$ by
$$
\theta _{x,M,P} (c) = {\hbox{1} \over {M}}\sum\limits_{\{v \in {{\cal T}_M '} \mid 
P(T_v x) = c\}} {w_v }
$$
for each $c \in C$. We say that $x \in X$ is $\varepsilon$,$M$-{\twelvebx generic for $P$} if
$$
\sum_{c \in C}\left| {\theta _{x,M,P}(c) - \mu(P^{-1}(c))} \right| < \varepsilon.
$$

The following is the Tree Ergodic Theorem and follows easily from proposition 2.2.2.
\proclaim Proposition 2.2.3 (Tree Ergodic Theorem). Let \X be an ergodic $\bf p$-endomorphism and $P: X \to C$ be a finite partition.  For every $\varepsilon > \hbox{0}$, we have for all large $M$, a set $G$ of measure $\hbox{1}-\varepsilon$ such that for each $x \in G$, $x$ is $\varepsilon$,$M$-generic for $P$. \par
\medjump
\noindent {\section 2.3. The Tree Rokhlin Lemma}
\medjump
Let \X be an ergodic m.p.s. A $\varepsilon$-{\twelvebx tree Rokhlin tower of height} $N+\hbox{1}$ in $\bf X$ is a collection of pairwise disjoint measurable sets $\{ B_{0},B_{1},\ldots 
,B_{N}\}$ in $X$ such that $B_{j}=T^{-j}B_{0}$ and whose union $\cup _{i = 
0}^N B_i $ has measure $>\hbox{1}-\varepsilon $.  If $\bf M$ is such a tree Rokhlin tower, we will let $\cup \bf M$ denote the union of the sets in $\bf M$.  We will refer to the set $B_{0}$ as 
the {\twelvebx base of the tower} and the set $B_{j}=T^{-j}B_{0}$ as the {\twelvebx j-th level of the tower}. If $B' \subseteq B$, we will refer to the union $ \cup _{i = 1}^N T^{ - i}B' $ as the {\twelvebx column of the tower over} {$B'$}. We now establish the analogues for 
$\bf p$-endomorphisms of the standard Rokhlin lemma and the Strong 
Rokhlin lemma. Their proofs follow along the same lines as propositions 5.2 
and 5.3 in [H,R], although their results were only stated for uniformly $p$-to-1 endomorphisms. For completeness, we include the proofs below.
\proclaim Proposition 2.3.1 (Tree Rokhlin Lemma). Let \X be ergodic, then for each $\varepsilon > \hbox{0}$ and $N \in {\bf N}$, there exists a $\varepsilon$-tree Rokhlin tower of height N+1 in $\bf X$.\par
\noindent{\twelvebx Proof:} For a set $D \subseteq X$ of positive measure, define the set
$$
B = \{x \in X \mid \hbox{min} (i \ge \hbox{0} \mid T^{i} x \in D) \equiv \hbox{0}\;\hbox{mod} 
(N + \hbox{1})\}\backslash D
$$
Note that $B$ is well-defined since by ergodicity, for a.a. $x$ in $X$, $T^{i}x 
\in D$ for some $i \in {\bf N}$. It suffices to show that $B \cap T^{ - 
j}B=\emptyset $, for each $\hbox{1} \le j \le N$. Suppose on the contrary, $x 
\in B \cap T^{ - j}B$ for some $\hbox{1} \le j \le N$. Then there exists $y 
\in B$ such that $y=T^{j}x$. By the definition of $B$, the set $\{x,Tx,\ldots 
,T^{j - 1}x,T^{j}x=y\}$ is disjoint from $D$, hence we have
$$
\hbox{min} (t \ge \hbox{0} \mid T^{t}x \in D )=\hbox{min} (t \ge \hbox{0} \mid T^{t}y \in 
D)+j.
$$
Since $x \in B$ and $y \in B$, taking mod($N+\hbox{1}$) of both sides forces $j=\hbox{0}$, 
which is a contradiction. Hence $B \cap T^{ - j}B=\emptyset $. Let $ 
{\bf M}= \{ T^{ - i}B \mid \hbox{0}\le i \le N \}$. 

Now, as $\bf X$ is ergodic, $\mu (B) \ge \hbox{1}/(N+ \hbox{1}) - \mu (D)$ and so
$$\mu (\cup_{i=0}^N T^{-i}B)=(N+\hbox{1}) \mu (B) \ge \hbox{1}-(N+\hbox{1})\mu(D).$$
Choosing $D$ such that $\mu (D)< \varepsilon /(N+\hbox{1})$ 
shows that the Rokhlin tower ${\bf M}$ has the desired property. {\symbol \char'04}
\medbreak
Before we prove the Strong Tree Rokhlin lemma, we introduce some notations 
which we will use in the rest of this thesis concerning distributions of 
partitions on a probability measure space. For a finite set $C$, define the 
{\twelvebx total variation norm} on the probability measures on $C$ by
$$
|\rho - \theta| = \sum\limits_{c \in C} {|\rho (c) - \theta (c)|}. 
$$
We write $\rho 
\mathop \sim \limits^\varepsilon \theta $ if $\left| \rho  - \theta \right|
<\varepsilon$. Suppose   $(\hbox{\char'12},{\cal B},\lambda )$ is a probability measure space. For a subset $G \subseteq \hbox{\char'12} $ of positive measure and a measurable finite partition $P:\hbox{\char'12}  \to 
C$, let $P | G$ denote the restriction of $P$ to $G$. Let $\lambda _{G}$ be a 
probability measure on $\cal B$ defined by $\lambda _{G}(B)=\lambda 
(B \cap G)/\lambda (G)$ for $B \subseteq \hbox{\char'12}$. 
Obviously, $\lambda =\lambda _{G}$ if $G=\hbox{\char'12}$. For a partition $Q:G 
\to C$, let $dist_{\lambda_{G}}(Q)$ be the probability measure on $C$ defined by 
$dist_{\lambda_{G}}(Q)(c)=\lambda _{G}(Q^{ - 1}(c))$ for each $c \in C$. We will refer 
to the sets $Q^{-1}(c) \subseteq G$, $c \in C$, as the {\twelvebx atoms} of $Q$. Where there is no ambiguity as to the measure on the domain of the partition, we will generally omit the subscript and just write $dist(Q)$. For a measure 
preserving transformation $T$ on $(\hbox{\char'12} ,{\cal B},\lambda )$ and a 
partition $P:\hbox{\char'12}  \to C$, we define the $P$-{\twelvebx name} of $x$ to be the infinite sequence 
$(P(x),P(Tx),\ldots )$ and the $P$-{\twelvebx n}-{\twelvebx name} of $x$ to be the finite sequence $(P(x),\ldots 
,P(T^{n-1}x))$.

If $P:\hbox{\char'12}  \to C$ and $Q:\hbox{\char'12}  \to C'$ are partitions 
into finite sets $C$ and $C'$, we say that $Q \le P$ if there exists a 
map $\pi :C \to C'$ such that $Q=\pi  \circ P$, i.e. knowing 
$P(x)$ determines $Q(x)$. 

We will need the following technical lemma that guarantees the measurability 
of the atoms of the various partitions that will be constructed in the proof 
of the Strong Tree Rokhlin lemma and in the proof of Theorem 2.4.1.

\proclaim Lemma 2.3.2. Suppose \X  $\in End(\bf {p})$. Consider a measurable finite partition $P:X 
\to C$. For any $C$,$n$-tree name $\tau $, the set $\{x \in X \mid \tau _x^P =\tau  \; \hbox{on} \; {\cal T}_n  ' \}$ is measurable.  \par

\noindent {\twelvebx Proof:} Note that $\tau $ induces a subset in $C^{n}$ whose elements are indexed by the nodes of length $n$. In fact, if $\left| v \right| =n$, we may associate 
the sequence $\bar {\tau }(v) =  (\tau (v),\ldots ,\tau (\sigma 
^{n - 1}v)) \in C^{n}$ to $v$. Let $\bar {C}(v)$ denote the set $\{ x 
\in X \mid \hbox {$P$-$n$-name of}\; x\; \hbox {is}\; \bar {\tau }(v)\}$. 

Let $G= \cup _{\left| v \right| = n} T_v X \cap \bar {C}(v)$, which is 
clearly measurable by the definition of the partial inverses $T_{v}$ in \S 
1.3. For $x \in X$, let $\mu _{x}$ be the conditional measure on $T^{ - n}x$ given $x$. Notice 
that the function $g:X \to $ [0,1] defined by $g(x)=\mu _{x}(G)$ is measurable. In particular, the set $\{x \in X \mid g(x)=\hbox{1} \}$ is 
measurable. However, this is precisely the set $\{x \in X \mid \tau _x^P 
=\tau \; \hbox{on}\; {\cal T}_{n}' \}$. {\symbol \char'4}
\medbreak
\proclaim Proposition 2.3.3 (Strong Tree Rokhlin lemma). Let ${\bf X} \in End({\bf p})$ be ergodic. Consider a measurable partition $P:X \to C$.  Then for each $\varepsilon >\hbox{0}$ and $N \in 
{\bf N}$, there exists a $\varepsilon $-tree Rokhlin tower of height $N$+1 in ${\bf X}$ whose base $B$ is independent of $P$, i.e. $dist(P|B)=dist(P)$.\par

\noindent{\twelvebx Proof:} Let $N'$ be an integer and $\varepsilon ' > \hbox{0}$, both to 
be specified later. Use the Tree Rokhlin lemma to build a $\varepsilon ' / \hbox{2}$-tree Rokhlin tower $\bf {M}'$ of height $N'+\hbox{1}$ with base $D$. Let $\bar 
P$ be the partition $P^{N'\tau }|D$. By Lemma 2.3.2, $\bar {P}$ 
defines a measurable partition of $D$. For each atom $\alpha $ in $\bar {P}$, 
divide $\alpha $ into $N+\hbox{1}$ measurable disjoint sets $\beta _0^\alpha ,\ldots 
,\beta _{N}^\alpha$ of equal measure. Consider the set
$$
B' = \bigcup_{\alpha \in \bar{P}} \bigcup_{0 \le i \le N} \bigcup_{{j \equiv i\bmod (N + 1)} \atop {0<j\le N'-N}} {T^{ - j}\beta 
_i^\alpha }.
$$
\noindent Notice that $B'$ is disjoint from $T^{ - j}B'$ for all $\hbox{0}<j\le N$. Let ${\bf M}= \{ T^{ - 
i}{B}' \mid \hbox{0}\le i \le N \}$. As $\bf M$ contains the $(N+\hbox{1})^{st}$ to $(N' 
- N)^{th}$ levels of ${\bf {M}'}$, we may choose $N'$ large enough and $\varepsilon '$ small enough such that the union of the sets in $\bf M$ has measure $>\hbox{1}-\varepsilon /\hbox{2}$.
 
From our definition of $B'$, we have
$$
dist(P|B') = dist(P|\bigcup_{i = 1}^{N' - 
N} T^{ - i}D).
$$
For any $\eta > \hbox{0}$, by decreasing $\varepsilon '$ and increasing $N'$ if necessary, we have 
$$
\left| dist(P|B') - dist(P) \right| \le \eta .
$$

If $\eta > \hbox{0}$ is small enough, we may remove at most $\varepsilon /\hbox{2}$ 
fraction of $B'$ to arrive at a set $B$ such that $dist(P|B)=dist(P)$ and 
$B$ forms the base of a $\varepsilon $-tree Rokhlin tower of height $N+{\hbox {1}}$. {\symbol \char'04}
\medjump
\noindent{\section 2.4. Proof of the Isomorphism Theorem }
\medjump
Unless otherwise specified, all m.p.s. in this section are 
$\bf p$-endomorphisms. The purpose of this section is to prove the following 
theorem. 
\proclaim Theorem 2.4.1. Let \X $\in End({\bf{p}})$, then $\bf X$ is tvwB if and only if ${\bf X}\cong B^{ + }({\bf{p}})$.\par

Our first goal is to prove that $B^{ + }({\bf{p}})$ is tvwB, which immediately proves half of Theorem 2.4.1.  To do this, we will prove an extension of Lemma 2.5 in [H,R] to $End(\bf p)$, which shows that any tree-adapted factor of a tvwB $\bf p$-endomorphism remains tvwB.  Before we prove this (Proposition 2.4.3), a technical lemma is in order. 

\proclaim Lemma 2.4.2. Given a m.p.s. \X, for each $n \in {\bf N}$, $\delta >\hbox{0}$ and 
$B \subseteq X$ such that $\mu (B)>\hbox{1}-\delta $, there exists a set $G$ of measure $>\hbox{1}-\sqrt \delta $ such that 
for each $x \in G$,
$$
{\hbox{1} \over n}\sum\limits_{0 < \left| v \right| \le n} {w_v \chi _B (T_v x)}  > \hbox{1} - \sqrt \delta  .
$$\par
\noindent {\twelvebx Proof: }For each $k \ge \hbox{0}$, let 
$$w_k(x,B)=\sum _{{\left| v \right| =k} \atop {T_v x \in B}} w_v .$$ 
Clearly, we have 
$$\mu (B) = \int {w_k (x,B)\,d\mu (x)}.$$ 
\noindent Hence, 
$$
\mu (B) = \int {{\hbox{1} \over n}\sum\limits_{i = 1}^n {w_i (x,B)\,d\mu (x)} } > \hbox{1} - 
\delta . 
$$
\noindent So,
$$
{1 \over n} \sum\limits_{i = 1}^n {w_i (x,B)} > \hbox{1} - \sqrt \delta 
$$
\noindent except on a set of measure at most $\sqrt \delta $. We finish the result by 
observing that 
$$
{\hbox{1} \over n} \sum\limits_{i = 1}^n {w_i (x,B)} = {\hbox{1} \over n} \sum\limits_{0 < 
\left| v \right|  \le n} {w_v } \chi _B (T_v x).
$$
{\symbol \char'04}
\medbreak
Suppose we have a measurable function $g:X \to R$ for a compact metric space 
$(R,d)$. We say that $g$ is {\twelvebx tree-adapted} if for $a.a.x$, $g$ assigns different values to the 
inverse images of $x$. If ${\cal D}$ is the Borel $\sigma $-algebra of $(R,d)$, 
we say that $g$ is {\twelvebx generating} if ${\cal B}=\mathop \vee \limits_{i = 0}^\infty T^{ - 
i}g^{ - 1}({\cal D})$. It is clear that a generating $g$ is necessarily 
tree-adapted.  We now prove the extension of Lemma 2.5 in [H,R].  The proof follows along similar lines, though we will also need the fact that tree adapted maps preserve the $p$-function (Proposition 2.1.3). 
\proclaim Proposition 2.4.3. Let \X $\in End({\bf p})$ and let \Y be a 
tree-adapted factor of $\bf X$. Suppose the function $h:X \to R$ is generating. If $({\bf X},h)$ is tvwB, then 
for any function $g:Y \to R' $ into a compact metric space $(R',d')$, $({\bf Y},g)$ is tvwB.\par

\noindent{\twelvebx Proof: }By normalizing $d$ and $d'$, we may assume that the 
metric spaces $R$ and $R'$ have unit diameter.  Let $\phi : X \to Y$ be a tree-adapted factor map.    Now, $h$ is generating and $g \circ \phi$ is $\cal B$-measurable.  Hence, for each $\varepsilon > \hbox{0}$, there exist some $s \in {\bf N}$, $\hbox{0}<\delta<\varepsilon$ and a set $G_\varepsilon$ of measure $> \hbox{1}-\varepsilon$ such that whenever $x$ and $x'$ are in $G_\varepsilon$ and $d(h(T^i x),h(T^i x'))<\delta$ for all $\hbox{0} \le i \le s$, then $d'(g(\phi(x)),g(\phi(x'))<\varepsilon$. 

Since $({\bf X},h)$ is tvwB, for a sufficiently large $N$, there exists a set $G\subseteq X$ of measure $> \hbox{1} -\delta^2$ such that whenever $z$ and $z'$ are in $G$, 
$$\bar t_N (\tau_z^h , \tau_{z'}^h)<\delta^2. \eqno(\hbox{1})$$
Moreover, by Lemma 2.4.2, we have a set $G'\subseteq X$ of measure $> \hbox{1} - \sqrt \varepsilon$ such that for each $x \in G'$ 
$$ {\hbox{1} \over N} \sum_{0<\left| v \right| \le n} w_v \chi _{G_{\varepsilon}} (T_v x) > \hbox{1}-\sqrt {\varepsilon}. \eqno(\hbox{2}) $$
Fix a pair of points $x$ and $x'$  in $G \cap G'$.  Using (1), we have some tree automorphism $A \in {\cal A}$ such that
$${\hbox{1} \over {N}} \sum_{0<\left| v \right| \le N} w_v d(h(T_v x), h(T_{Av} x')) < \delta^2.\eqno(\hbox{3})$$
Let $D_\delta^{s+}=\{ v \in {\cal T}' \mid d(h(T^i(T_v x)),h(T^i(T_{Av} x')))\ge\delta \hbox{ for some 0} \le  i \le s \}$.  Then by (3), 
$$ {\hbox{1} \over {N}} \sum_{ v \in D_\delta^{s+}} w_v < (s+\hbox{1}) \delta + s/N. \eqno(\hbox{4})$$
Let $V=\{ v \in {\cal T}' \mid v \notin D_\delta^{s+}, T_v x \in G_{\varepsilon}, T_{Av} x' \in G_{\varepsilon} \}$.  Then notice that whenever $v \in V$, by our choice of $\delta$, 
$$d'(g(\phi (T_v x)), g(\phi (T_{Av} x'))) < \varepsilon. \eqno(\hbox{5})$$ 
Note that the tree automorphism $A$ corresponds to a bijection of the trees of inverse images of $x$ and $x'$.  This in turn pushes down via $\phi$ to a bijection of the trees of inverse images of $\phi (x)$ and $\phi (x')$.  Since $\phi$ is tree-adapted and hence preserves the $p$-function (proposition 2.1.3), this last bijection in turn yields a tree automorphism $A'$.  By (2), 
(4) and (5), we have 
$$
\eqalign{
  {\hbox{1} \over N}\sum\limits_{0 < \left| v\right|  \le N} {w_v d'(g(S_v \phi (x)),g(S_{A'v} \phi (x')))} &= {\hbox{1} \over N}\sum\limits_{0 < \left| v \right| \le N} {w_v d'(g(\phi (T_v x)),g(\phi (T_{Av} x')))} \cr 
& \le {\hbox{1} \over N}\sum\limits_{{0 < \left| v \right| \le N} \atop 
  {T_v x \notin G_\varepsilon}}  {w_v }  + {\hbox{1} \over N}\sum\limits_{{0 < \left| v \right| \le N} \atop 
  {T_{Av} x' \notin G_\varepsilon}}  {w_v }  + {\hbox{1} \over N}\sum\limits_{{0 < \left| v \right| \le N} \atop {v \in D_\delta ^{s + }}}   {w_v }  \cr 
 & \quad + {\hbox{1} \over N}\sum\limits_{v \in V} {w_v d'(g(\phi (T_v x)),g(\phi (T_{Av} x')))}  \cr &< \sqrt \varepsilon   + \sqrt \varepsilon   + (s + \hbox{1})\delta + s/N + \varepsilon  < \hbox{5} \sqrt \varepsilon   \cr}
$$
for all sufficiently large $N$ and small $\delta$.

As $\mu (G' \cap G)>\hbox{1} - \hbox{2}\sqrt \varepsilon $, if $\varepsilon $ is 
sufficiently small, then there will be a large set $H \subseteq Y$ such that 
each $y \in H$ is the image of some point in $G'\cap G$. Thus, the result follows 
from the preceding calculation. {\symbol \char'04}
\medbreak

Proposition 2.4.3 implies that if $({\bf X},f)$ is tvwB for some 
generating $f$, then $({\bf X},g)$ is also tvwB for any compact-valued 
function $g$. We thus see that the m.p.s. $\bf X$ is tvwB if and only if 
$({\bf X},f)$ is tvwB for some generating $f$. Hence, $B^{ + 
}({\bf p}) $ is tvwB since its standard generator (zero 
coordinate partition) yields the same tree name for all points in $B^{ + 
}({\bf p}) $  by choosing the tree partition $K_{B^{+}({\bf p})} $ to be the 
zero coordinate partition.  We have thus proved the easier half of Theorem 2.4.1, which we record as the proposition below.

\proclaim Proposition 2.4.4. Suppose \X $\in 
End({\bf p})$. If ${\bf X} \cong B^{ + }({\bf p})$, then $\bf{X}$ is tvwB. {\symbol \char'04} \par
Our goal now is to show that if a $\bf{p}$-endomorphism is tvwB, 
then it is one-sided Bernoulli. It will be convenient for our presentation 
to work with functions defined on $\bf{p}$-endomorphisms with range in the metric space [0,1) (with the absolute value 
metric $\left| x-y \right|$).  As we will see below, this gives us a fairly natural way to partition the range space and to define tree name distributions induced by a function on a $\bf p$-endomorphism.  For a m.p.s. \Y and two 
functions $g, h:Y\to [\hbox{0},\hbox{1})$, let $\left| g-h \right| ={\Vert {g - h} \Vert} _{L^1} $. 

For any function $g:Y \to [\hbox{0},\hbox{1})$ and $N \in {\bf N}$, we will construct a 
``discretize'' version of $g$ which closely approximates it. Specifically, 
for each $N$, let $P_{N}$ be the set of dyadic intervals 
$\{ [\hbox{0},\hbox{1}/\hbox{2}^{N}),\ldots,[(\hbox{2}^{N} - \hbox{1})/\hbox{2}^{N},\hbox{1})\}$ of length 
$\hbox{2}^{-N}$. Let $D_{N}$ be the set of midpoints of the intervals in 
$P_{N}$, i.e. $D_{N}=\{(\hbox{2}t+\hbox{1})/\hbox{2}^{N + 1}, \hbox{0} \le t < \hbox{2}^{N}\}$. Define the function $g_{N}:Y \to D_{N}$ by setting
$$
g_N (x) = {{\hbox{2}t + \hbox{1}} \over \hbox{2}^{N + 1}}\quad {\hbox {if}}\; x \in [{t \over \hbox{2}^N} , {{t + \hbox{1}} \over \hbox{2}^N} ),\; \hbox{0} \le t < \hbox{2}^N 
$$
\noindent Clearly, $\left| g_{N} - g \right| < \hbox{1}/\hbox{2}^{N}$. Observe that 
$g_{N}$ assumes values in the finite set $D_N$ so that we may regard $g_{N}$ as a partition of $Y$. 

Given a function $g: Y \to [\hbox{0},\hbox{1})$ and positive integers $m$ and $n$, recall from \S 1.3 
that we have the partition $g_{m}^{n\tau }:Y \to (D_{m})^{n\tau 
}$ by mapping a point $y$ to its $g_{m}$-$n$-tree name and the partition $g_{m}^{n\nabla 
}:Y \to (D_{m})^{n\nabla }$ by mapping a point $y$ to the equivalence 
class in $(D_{m})^{n\nabla }$ containing $g_{m}^{n\tau }(y)$. Let us 
also define the partition $g_{m}^{n + }:Y \to (D_{m})^{n}$ by 
mapping a point $y$ to its $g_{m}$-$n$-name. Observe that 
$g_{n'}^{n\tau } \le g_{m'}^{m\tau }$, $g_{n'}^{n\nabla } \le g_{m'}^{m\nabla }$ and 
$g_{n'}^{n + } \le g_{m'}^{m + }$ if $n' \le m'$ and $n \le m$.  The following proposition shows that closeness in tree name distributions implies closeness in forward name distributions.

\proclaim Proposition 2.4.5. Let \Y $\in End({\bf p})$.  For $n$, $m \in {\bf N}$ and any functions $g, h: Y \to [\hbox{0},\hbox{1})$, if
$$dist(g_n^{m\nabla })\mathop \sim \limits^\varepsilon 
dist(h_n^{m\nabla }),$$ 
then 
$$dist(g_n^{m + })\mathop \sim 
\limits^\varepsilon dist(h_n^{m + }).$$ \par
\noindent {\twelvebx Proof: } For each element $\zeta  
\in (D_{n})^{m}$ and $\alpha  \in (D_{n})^{m\nabla }$, and for any 
representative $\beta  \in \alpha $, let $w(\zeta ,\alpha )$ be 
the total weights of nodes $v$ of length $m$ such that $(\beta (v), \beta 
(\sigma v),\ldots ,\beta (\sigma ^{m - 1}(v))=\zeta $. (This is clearly independent of the representative chosen.)  Then
$$
\eqalign
{\sum\limits_{\zeta \in (D_n )^m} {\left| dist(g_n^{m + })(\zeta ) - 
dist(h_n^{m + })(\zeta ) \right|} & \le \sum\limits_{\zeta \in (D_n )^m} 
{\sum\limits_{\alpha \in (D_n )^{m\nabla }} {w(\zeta ,\alpha )\left| 
dist(g_n^{m\nabla })(\alpha ) - dist(h_n^{m\nabla })(\alpha ) \right|}} \cr 
 &= \sum\limits_{\alpha \in (D_n )^{m\nabla }} {\sum\limits_{\zeta \in (D_n 
)^m} {w(\zeta ,\alpha )\left| dist(g_n^{m\nabla })(\alpha ) - dist(h_n^{m\nabla 
})(\alpha ) \right| }} \cr 
&= \sum\limits_{\alpha \in (D_n )^{m\nabla }} {\left| dist(g_n^{m\nabla 
})(\alpha ) - dist(h_n^{m\nabla })(\alpha ) \right|} < \varepsilon . \cr }
$$
{\symbol \char'04}
\medbreak
Let \X and \Y be ${\bf p}$-endomorphisms.  Suppose $g:X \to [\hbox{0},\hbox{1})$ and $h:Y \to [\hbox{0},\hbox{1})$, we wish to define a $\bar {t}$ 
distance between the processes $({\bf X},g)$ and $({\bf Y},h)$ analogous 
to the $\bar {d}$ distance in Ornstein's Theory. To do this, we first let 
$$
\bar {t}_n (({\bf X},g),({\bf Y},h)) = \int {\bar {t}_n } (\tau _x^g ,\tau _y^h )\,d\mu 
(x)d\nu (y).
$$
\noindent Then define
$$
\bar {t}(({\bf X},g),({\bf Y},h)) = \hbox{liminf}\;\bar {t}_n (({\bf X},g),({\bf Y},h)).
$$

We remark that this is {\twelvebx not} the definition in [H,R] in which 
they define $\bar {t}_n (({\bf X},g),({\bf Y},h))$ as an infimum of 
the integral of $\bar {t}_n (\tau _x^g ,\tau _y^h )$ over the set of 
one-sided couplings of ${\bf X}$ and ${\bf Y}$ (see Chapter 3). Here, we only define 
it with respect to the product measure. Note that closeness in $\bar {t}$ 
essentially means that for some large $n$, $\bar {t}_n (\tau _x^g ,\tau _y^h )$ 
is small for any $x$ in a large set in $X$ and for any $y$ in a large set in $Y$. 

The following proposition shows that if $\bf{X}$ is tvwB, then $\bar 
{t}(({\bf X},g),({\bf Y},h))$ is small provided that $dist(g_n^{n\nabla 
})$ is sufficiently close to $dist(h_n^{n\nabla })$ for some $n$. The argument 
follows along the same lines as Lemma 4.4 in [H,R], though the notations 
found there need to be modified for our situation, since we are considering $\bf p$-endomorphisms so that the weights of the nodes within a given level are not uniform.  

\proclaim Proposition 2.4.6. Suppose ${\bf X}$ is tvwB and $g: X\to [\hbox{0},\hbox{1})$. For all $\varepsilon >\hbox{0}$, there exist $\delta >\hbox{0}$ and $N \in 
{\bf N}$ with the following property: for any ${\bf Y} \in End({\bf p})$ and any function $h:Y\to [\hbox{0},\hbox{1})$, if $dist(g_N^{N\nabla })\mathop 
\sim \limits^\delta dist(h_N^{N\nabla })$, then $\bar 
{t}(({\bf X},g),({\bf Y},h))<\varepsilon $.\par
\noindent{\twelvebx Proof:} Let $\eta >\hbox{0}$, to be specified later. Since $\bf X$ is tvwB, we may choose $N$ large enough such that there exists a set $G$ of measure $>\hbox{1}-\eta$ such that $\bar {t}_N (\tau _{z}^g ,\tau _{z'}^g ) < 
\eta $ whenever $z$ and $z'$ are in $G$.  Clearly, we may assume that 
$\hbox{1}/\hbox{2}^{N}<\eta $. Fix $x' \in G$.  If $dist(g_N^{N\nabla 
})\mathop \sim \limits^\eta dist(h_N^{N\nabla })$, then we have a set 
$H$ in $Y$ of measure at least $\hbox{1}-\hbox{2}\eta$ such that if $y \in 
H$, then 
$$\bar {t}_N (\tau _{x'}^g ,\tau _y^h ) < \eta + (\hbox{1} / \hbox{2}^{N}) 
< \hbox{2}\eta .$$
Let $\tau' = \tau_{x'}^g$.  Create a tree name $\tau : {\cal T}' \to \hbox{[0,1)}$ by tiling with copies of $\tau'$ as follows:  If $\hbox{0}<\left| v \right| \le N$ and 
$\left| v' \right| =kN $ for some integer $k \ge \hbox{0}$, we define $\tau $ by setting $\tau 
(vv')=\tau' (v)$. We would like to show that $\tau$ is close to $\tau_x^g$ and $\tau_y^h$ for large sets in $X$ and in $Y$ in ${\bar t}_{kN}$ for all $k \in {\bf N}$.  Roughly speaking, the idea is to think of each tree of height $kN$ as consisting of subtrees of height $N$ and to realize that the tvwB condition guarantees that the tree names of these subtrees are close to $\tau '$ in $\bar t_N$ on average.  

We will now define inductively, $N$-levels at a time, a tree automorphism $A$ 
which makes $\bar {t}_{kN} (\tau _x^g , \tau)$ small on a large set in $X$ for 
all $k$. For $x \in X$, let $A_{x}$ be a tree automorphism which realizes the 
minimum in the definition of $\bar {t}_N (\tau _x^g ,\tau )$. For $\hbox{0}<\left| v \right| \le N$, set $A(v)=A_{x}(v)$. Inductively, for each $v \in {\cal T}$ such that 
$jN<\left| v \right| \le (j+\hbox{1})N$, $j \ge \hbox{1}$, write $v=v'u$ for unique nodes 
$v'$ and $u$ such that $\left| u \right| =jN$ and $\hbox{0}< \left| v' \right| \le N$, and define 
$A(v)=A_{T_u x} (v')A(u)$.

For $j \ge \hbox{0}$ and $x \in X$, let $S_{j}(x)=\{v \in {\cal T} \mid \left| v \right| =jN \; \hbox{and} \; 
T_{v} {x} \in G\}$.  Let $w(S_{j}(x))$ denote the sum $\sum\nolimits_{v \in 
S_{j} (x)} {w_v } $.  Then by Lemma 2.4.2, for all $k \ge \hbox{0}$, there exists a set $G_k \subseteq X$ of measure $\ge \hbox{1} - \sqrt {\eta}$ such that 
$${\hbox{1} \over k} \sum_{j=0}^{k-1} w(S_{j}(x)) \ge \hbox{1}-\sqrt{\eta}$$
for each $x \in G_k$.   By our construction of the tree automorphism $A$, $\tau_y^g$  and $\tau '$ are matched to within $\eta$ in $\bar t_N$ whenever $y=T_v x$ for $v \in S_{j}(x)$.  For each $x \in G_k$, by calculating $\bar t_{kN}(\tau_x^g, \tau)$ as an average of $\bar t_N(\tau_y^g,\tau')$ over all $y \in T^{-jN}x, \hbox{0} \le j \le k-\hbox{1}$, we have 
$$
\eqalign{
\bar {t}_{kN} (\tau _x^g ,\tau ) &={\hbox{1} \over kN}\sum\limits_{j=0}^{k-1} \sum\limits_{\left|v \right|=jN} w_v N \bar t_N (\tau_{T_v x}^g, \tau')  \cr &={\hbox{1} \over kN}\sum\limits_{j=0}^{k-1}\Bigl( \sum\limits_{v \in S_{j}(x)} w_v N \bar t_N (\tau_{T_v x}^g, \tau') + \sum\limits_{v \notin S_{j}(x)} w_v N \bar t_N (\tau_{T_v x}^g, \tau') \Bigr) \cr
&\le {\hbox{1} \over kN}\sum\limits_{j = 0}^{k - 
1} {N\eta w(S_{j} (x))}  + N(\hbox{1} - w(S_{j} (x))) \cr 
&\le \hbox{1} - {1 \over k}\sum\limits_{j = 0}^{k - 1} {w(S_{j} (x)) (\hbox{1} - \eta )} 
\cr
&\le \hbox{1} - (\hbox{1} - \sqrt{\eta})(\hbox{1} - \eta ). \cr }
$$
\noindent Hence for all but a set of measure $\hbox{2} \sqrt {\eta}$ in $X$, we have $\bar 
{t}_{kN} (\tau _x^g ,\tau )<\hbox{2} \sqrt {\eta} .$ 

By the same argument, for all but a set of measure $\hbox{4}\sqrt \eta$ in $Y$, we 
have $\bar {t}_{kN} (\tau _y^h ,\tau )<\hbox{4}\sqrt \eta $. Hence, by choosing 
$\eta $ sufficiently small, we have 
$$
\eqalign{
 \int {\bar {t}_{kN} } (\tau _x^g ,\tau _y^h )\,d\mu (x)d\nu (y) &\le \int 
{\bar {t}_{kN} } (\tau _x^g ,\tau )\,d\mu (x)d\nu (y) + \int {\bar {t}_{kN} } 
(\tau ,\tau _y^h )d\mu (x)d\nu (y) \cr 
 &\le \hbox{4}\sqrt \eta + \hbox{8}\sqrt \eta < \varepsilon \cr 
 }
$$
\noindent Since $k$ is arbitrary, the result now follows. {\symbol \char'04}
\medbreak
Let \X and consider a partition $P:X \to C$ 
for some finite set $C$. Given $M \in {\bf{N}}$ and $x \in X$, recall from \S 2.2 that we defined a
measure $\theta _{x,M,P}$ on $C$ by
$$
\theta _{x,M,P} (c) = { \hbox{1} \over M }\sum\limits_{\{v \in {{\cal T}'_M} \mid 
P(T_v x) = c\}} {w_v } 
$$
for each $c \in C$. Also recall that $x \in X$ is $\varepsilon$,$M$-generic for $P$ if
$$
\left| {\theta _{x,M,P} - dist(P)} \right| < \varepsilon .
$$
Note that if $\bf X$ is ergodic, the Tree Ergodic Theorem implies that for all sufficiently large $M$, there is a large set  $G \subseteq X$ such that all points $x \in G$ are $\varepsilon$,$M$-generic for $P$.  

Combined with proposition 2.4.6, the next proposition (the perturbation lemma), shows that 
if $\bf X$ is tvwB and $\bar {t}(({\bf X},f),({\bf Y},g))<\varepsilon $, then we only 
need to perturb $g$ slightly (depending on $\varepsilon$) to 
get a function $h$ on $Y$ such that the processes $({\bf X},f)$ and $({\bf Y},h)$ are as close in $\bar {t}$ as we like. This will be the key 
observation in the proof of Sinai's Theorem (Proposition 2.4.9). 

For the proof of the perturbation lemma (Proposition 2.4.7) and the copying lemma (Proposition 
2.4.12), we will be constructing a function on a chosen tree Rokhlin tower 
${\bf M}$ using one or more chosen tree names in a way analogous to 
``painting'' columns of a Rokhlin tower with a name as in the proof of 
Ornstein's theorem. Since this construction is central to the proofs, let us now define it explicitly. Suppose \Y $\in End({\bf p})$, $N \in 
{\bf N}$ and $A \in {{\cal A}_{N}}$. Given a tree name $\tau :{\cal 
T}_{N}'  \to R$ and a point $y \in Y$, we may define 
a function $h: \{S^{-t}y \mid {\hbox{1} \le t \le N} \} \to R$ such that 
$h(S_{v}y)=\tau (Av)$ for all $v \in {\cal T}_{N} '$. 
We will refer to $h$ as a {\twelvebx laying of $\tau$ on} $\{S^{-t}y \mid \hbox{1} \le t \le N \}$ {\twelvebx via $A$}. Doing this for all $y$ in the base $B$ of a tree Rokhlin tower ${\bf M}$ of height $N+ \hbox{1}$ for fixed $A \in {\cal A}_N$ and tree name $\tau$ defines a function on ${\cup \bf M} \backslash B$.

\proclaim Proposition 2.4.7 (Perturbation Lemma). Let \X and \Y be ergodic. Suppose $\delta > \hbox{0}$ and $f:X\to [\hbox{0},\hbox{1})$.  For any $\varepsilon > \hbox{0}$, $N \in {\bf N}$ and for any function $g:Y\to [\hbox{0},\hbox{1})$ such that $\bar 
{t}(({\bf X},f),({\bf Y},g))<\delta ^{4}$ , there exists some function $h:Y\to [\hbox{0},\hbox{1})$ such that $\left| g-h \right| <\hbox{8}\delta$  and $dist(h_N^{N\nabla})\mathop \sim \limits^\varepsilon dist(f_N^{N\nabla })$.\par 
\noindent{\twelvebx Proof:} There is no loss of generality in assuming that $N$ is large enough such that $\hbox{1}/\hbox{2}^{N} \le \delta ^{4}$. Let $\eta<\delta $ be specified later. 
Choose $M \in {\bf N}$ and a set $G \subseteq X$ of measure $>\hbox{1}-\eta$ such that the following conditions hold: 
\itemitem{a)} $N/M<\eta /\hbox{2}$
\itemitem{b)} $x$ is $\eta$,$(M-N)$-generic for $f_N^{N\nabla }$ for all $x \in G$
\itemitem{c)} $\int {\bar {t}_M } (\tau _x^f ,\tau _y^g )\,d\mu (x)d\nu (y) < \delta ^4.$

\noindent The Tree Ergodic theorem and the hypothesis show that such an integer $M$ and set $G$ exist. 

From c), we have
$$
\int {\bar {t}_M } (\tau _x^f ,\tau _y^g )\,d\nu (y) < \delta ^2
$$
\noindent except on a set $H \subseteq X$ of measure $<\delta ^{2}$. Choose some point $x'  \in G\backslash 
H$. It follows that  on a set $G' \subseteq Y$ of measure $\ge 
\hbox{1}-\delta$, we have $
{\bar {t}_M } (\tau _{x'}^f ,\tau _y^g ) \le \delta 
$ for $y \in G'$. By the Strong Tree Rokhlin Lemma, we may build a $\eta 
/\hbox{2}$-tree Rokhlin tower ${\bf M}$ of height $M+\hbox{1}$ in $\bf{Y}$ with base $B$ 
such that $\nu_B (G') \ge \hbox{1}-\delta $.

We now construct the required function $h$. For each $y \in Y$, let $A_{y}\in {\cal A}_M$ be a tree 
automorphism which realizes $\bar {t}_M (\tau _y^g ,\tau _{x'}^f )$. 
For each atom $\alpha$ in the partition $g_M^{M\tau }|B$, we pick a 
representative point $y(\alpha ) \in \alpha$. We define $h$ on ${\cup{\bf M}}\backslash B$ by laying $\tau _{x'}^f $ on 
$\{S^{-t}y \mid {\hbox{1} \le t \le N} \}$ via $A_{y(\alpha )}$ for each point 
$y \in \alpha $ and then for each atom $\alpha $.  Note that for 
a fixed node $v$ and a fixed atom $\alpha $ in $g_M^{M\tau }|B$, $h$ is 
constant on $S_{v}\alpha $, thus $h$ defines a measurable function on ${\cup \bf 
M}\backslash B$. We then extend $h$ measurably to the rest of the space 
in any way we desire.

For each atom $\alpha  \in g_M^{M\tau }|B$, let $C_{\alpha }$ be 
the column of the tower over $\alpha $. Let us say that $C_{\alpha }$ is a {\twelvebx good 
column} if $\alpha $ contains a point in $G'$. Now, if $y \in 
G'\cap B$ and $y \in \alpha $, then since $\hbox{1}/\hbox{2}^{M}< \delta$, we have
$$
\bar {t}_M (\tau _{y(\alpha )}^{g_M} ,\tau _{x'}^f ) \le\bar {t}_M (\tau _{y(\alpha )}^g ,\tau _{x'}^f )+\delta \le \bar {t}_M 
(\tau _y^g ,\tau _{x'}^f ) +\hbox{2} \delta < \hbox{3}\delta . \eqno(\hbox{1})
$$
\noindent For each good column $C_{\alpha}$, a simple calculation using (1) and our definition of $h$ shows that
$$
\int\limits_{C_\alpha } {\left|g_M (z) - h(z)\right| \, d\nu (z)} 
< \hbox{3}\delta \nu (C_\alpha ).  
 $$
\noindent Since $\nu_B (G') \ge \hbox{1}-\delta $, conditionally on $B$, the bases of the good columns is a set of measure $\ge \hbox{1}-\delta$. Thus,  
$$
\nu ( \{ \cup C_\alpha \mid C_\alpha \; \hbox{is good} \} ) \ge (\hbox{1} - \delta )\nu (\cup_{i=1}^M S^{-i}B)
> \hbox{1} - \hbox{2} \delta .  
$$
\noindent Hence, 
$$
\eqalign{
\left| g - h\right| &\le \left|g - g_M \right| + \left|g_M - h\right| \cr
 &\le \hbox{1} / \hbox{2}^M + \hbox{5} \delta < \hbox{8}\delta . \cr} 
 $$

For the second part of the conclusion, fix an atom $\alpha  \in h_{N}^{M\tau }|B$.  Then there 
exists some tree automorphism $A$ such that 
$h_{N}(S_{v}y)=f_{N}(T_{Av}x')$ for all $y \in \alpha $ and 
for all $v \in {\cal T}_{M}'$. Hence, for each $\xi  
\in (D_{N})^{N\nabla }$, if $f_N^{N\nabla }(T_{v}x')=\xi $ for 
$v \in {\cal T}_{M - N}'$, then $h_N^{N\nabla }(S_{A^{ - 1}v} y)=\xi $ 
for all $y \in \alpha $. Let $\bar {\xi }=\{v \in {\cal T}_{M-N}' \mid f_N^{N\nabla }(T_{v}x')=\xi \}$. \noindent Since tree 
automorphisms preserve weights of nodes, we have
$$
\eqalign{
 dist(h_N^{N\nabla}\mid \bigcup _{0 < j \le M - N} {S^{ - 
j}\alpha} )(\xi ) &= \sum\limits_{v \in \bar {\xi }} {\nu (S_{A^{ - 1}v} \alpha )\over {(M-N) \nu (\alpha )}} \cr 
&={\hbox{1} \over {M - N }}\sum\limits_{v \in \bar {\xi }} {w_v } \cr 
&= \theta _{x',M - N,f_N^{N\nabla }} (\xi ). \cr} 
$$

\noindent Thus, 
$$
dist(h_N^{N\nabla }|\bigcup_{0 < j \le M - N} {S^{ - 
j}\alpha} ) = 
\theta _{x',M - N,f_N^{N\nabla }} .
$$

\noindent Obviously, this is true for each $\alpha  \in h_{N}^{M\tau }| B$. 
Hence, (b) implies that 
$$dist(h_N^{N\nabla }|\bigcup _{0 < j \le M - N} {S^{ - j}B})
\mathop \sim \limits^\eta dist(f_N^{N\nabla }).$$

\noindent Using (a) and the fact that the measure of ${\cup \bf M}$ is large, we also have

$$
dist(h_N^{N\nabla })\mathop \sim \limits^\varepsilon dist(f_N^{N\nabla })
$$
\noindent for $\eta $ small enough. This completes the proof. {\symbol \char'04} 
\medbreak
By imitating the proof of the previous proposition, but without the need to make $h$ close to a predefined $g$, the following proposition is immediate.

\proclaim Proposition 2.4.8. Suppose \X and \Y are both ergodic $\bf p$-endomorphisms. For any function $f:X\to [\hbox{0},\hbox{1})$, and for all $\varepsilon $ and $N$, we have some function $h:Y\to [\hbox{0},\hbox{1})$ such that $dist(h_N^{N\nabla })\mathop \sim 
\limits^\varepsilon dist(f_N^{N\nabla })$.  {\symbol \char'04}

For \Y and a function $g:Y\to [\hbox{0},\hbox{1})$, we have a map 
$g^{{\bf N}^\ast}:Y \to \hbox{[0,1)}^{{\bf N}^\ast }$ defined by $g^{{\bf N}^\ast
}(y)=(g(y),g(Sy),\ldots )$. The measure $\nu $ pushes forward via $g^{{\bf N}^\ast}$ 
to the shift invariant measure $dist({\bf Y},g)=\nu  \circ (g^{{\bf N}^\ast
})^{ - 1}$ on $\hbox{[0,1)}^{{\bf N}^\ast}$.  We will refer to the 
map $g^{{\bf N}^\ast }$ as the {\twelvebx g-name map}. 

We will now prove the analogue of Sinai's theorem in Ornstein's theory for tvwB ${\bf p}$-endomorphisms.  In fact, we will need a slightly stronger version of it.  

\proclaim Proposition 2.4.9 (Strong Sinai's Theorem). Suppose \X is tvwB. Let $f:X\to [\hbox{0},\hbox{1})$ be a generating function. Given $\varepsilon >\hbox{0}$, there exist $\delta$ and $N$ such that if \Y is ergodic and if $g:Y\to [\hbox{0},\hbox{1})$ satisfies $dist(g_N^{N\nabla })\mathop \sim \limits^\delta 
dist(f_N^{N\nabla })$, then there exists a tree-adapted function $h:Y\to [\hbox{0},\hbox{1})$ such that $dist({\bf Y},h)=dist({\bf X},f)$ and $\left| h - 
g \right| <\varepsilon $.\par
\noindent{\twelvebx Proof:} By Proposition 2.2.1, $\bf{X}$ is ergodic. The strategy is 
to apply Propositions 2.4.6 and 2.4.7 repeatedly. Precisely, choose $N_{k} 
\nearrow \infty $ and $\varepsilon _{k} \searrow \hbox{0}$ such that for 
each $k$, $N_{k}$ and $\varepsilon _{k}$ correspond to $(\varepsilon /{\hbox{2}^{k 
+ 3}})^{4}$ in Proposition 2.4.6 applied to $({\bf X},f)$. We show that if 
$N=N_{1}$ and $\delta =\varepsilon _{1}$, the result holds. 

To see this, notice that if
$$
dist(g_{N_1}^{N_1 \nabla })\mathop \sim \limits^{\varepsilon _1 }  dist(f_{N_1}^{N_1 
\nabla }),
$$
\noindent
we have by Proposition 2.4.6,
$$
\bar {t}(({\bf X},f),({\bf Y},g)) < (\varepsilon / \hbox{16})^4.
$$
\noindent Hence, by Proposition 2.4.7, we have some function $g^{2}$ on $Y$ such that 
$\left| g-g^{2}\right|<\varepsilon / \hbox{2}$ and 
$$
dist((g^2)_{N_2}^{N_2 \nabla })\mathop \sim \limits^{\varepsilon _2 } 
dist(f_{N_2}^{N_2 \nabla }).
$$
\noindent Inductively, we obtain a sequence of functions $g^{j}$'s such that $\left| g^{j} - g^{j + 1}\right| <\varepsilon / \hbox{2}^{j }$ and 
$$
dist((g^j)_{N_j}^{N_j \nabla })\mathop \sim \limits^{\varepsilon _j } 
dist(f_{N_j}^{N_j \nabla }).
$$
Hence, the $g^{j}$ approach some function $h$ pointwise $a.e$. Note that $\left| g-h \right| <\varepsilon.$ 

To see that $dist({\bf Y},h)=dist({\bf X},f)$, note that for any $N \in {\bf N}$ and any  $\theta > \hbox{0} $, we have by proposition 2.4.5,
$$
dist((g^{j})_N^{N+})\mathop \sim \limits^\theta
dist(f_N^{N+})
$$
\noindent for all large $j$.  Thus, for any given $\gamma >\hbox{0}$, we have
$$
dist({\bf Y},g^{j})\mathop \sim \limits^\gamma 
dist({\bf X},f)
$$
in the $w$*-topology for all large $j$. Now, as $g^{j} \to h$ a.e., we also have 
for large $j$,
$$
dist({\bf Y},h)\mathop \sim \limits^\gamma 
dist({\bf Y},g^{j})\mathop \sim \limits^\gamma 
dist({\bf X},f).
$$
\noindent Since $\gamma$ is arbitrary, we have 
$dist({\bf Y},h)=dist({\bf X},f)$. 

It remains to show that $h$ is tree-adapted. By construction, for any $N \in {\bf N}$ and $\gamma > \hbox{0}$, $dist((g^j)_N^{1 \nabla}) \mathop \sim \limits^\gamma dist(f_N^{1 \nabla})$ for all large $j$.  As $g^j \to h$, it immediately follows that for any $M \in {\bf N}$, $dist(h_M^{1 \nabla}) = dist(f_M^{1 \nabla})$.  Since $f$ is generating and hence tree-adapted, so is $h$. 
{\symbol \char'04}
\medbreak
The following is an immediate consequence of propositions 2.4.8 and 2.4.9.

\proclaim Corollary 2.4.10 (Sinai's Theorem). Suppose \X is tvwB. Let $f:X\to [\hbox{0},\hbox{1})$ be a generating function. If \Y is ergodic, then there exists a tree-adapted function $h:Y\to [\hbox{0},\hbox{1})$ such that $dist({\bf Y},h)=dist({\bf X},f)$. {\symbol \char'04} \par

For a function $g$ 
on $Y$, consider the $g$-name map $g^{{\bf N}^\ast 
}:Y \to \hbox{[0,1)}^{\bf{N}^\ast }$ defined previously. Similarly, if $f$ is a function on $X$, consider the $f$-name map $f^{{\bf N}^\ast}:X \to \hbox{[0,1)}^{\bf{N}^\ast }$. If $dist({\bf X},f)=dist({\bf Y},g)$ for a generating $f$ on $X$, then we have a factor map $\pi : {\bf Y} \to {\bf X}$ defined by $\pi = (f^{{\bf N}^\ast}) 
^{ - 1} \circ g^{{\bf N}^\ast }$.  Note that $\pi$ is tree adapted if $g$ is a tree-adapted function on $Y$. Thus, Sinai's 
Theorem implies that if ${\bf X}$ is tvwB, then it is a tree-adapted 
factor of ${\bf Y}$. 

Before we establish the copying lemma, Proposition 2.4.12, we need the 
following result which shows that if $dist({\bf Y},g)=dist({\bf X},f)$ for a generating $f$ and a tree-adapted $g$ as in Sinai's Theorem, the tree name distributions induced by $g$ and $f$ are the same.

\proclaim Proposition 2.4.11. Suppose ${\bf X}$  and ${\bf Y}$  are in $End({\bf p})$. Let $f:X\to [\hbox{0},\hbox{1})$ be a generating function. Suppose $g:Y\to [\hbox{0},\hbox{1})$ is a tree-adapted function such that $dist({\bf Y},g)=dist({\bf X},f)$. Then for all $M 
\in {\bf N}$,
$$
dist(g_M^{M\nabla })=dist(f_M^{M\nabla }).
$$ \par
\noindent{\twelvebx Proof:} This follows directly from proposition 2.1.4 applied to the tree-adapted factor map $\pi : {\bf Y} \to {\bf X}$ defined by $\pi = (f^{{\bf N}^\ast}) 
^{ - 1} \circ g^{{\bf N}^\ast }$ and on noticing that $g_M=f_M \circ \pi$ {\symbol \char'04}
\medbreak
If $h, h':Y\to [\hbox{0},\hbox{1})$, let $h \vee h'$ be the joined 
function defined by $h \vee h' (y)=(h(y),h'(y))$. By analogy 
with our definitions of $h_M^{M\tau }$ and $h_M^{M\nabla }$, we can construct the 
partitions $(h_M \vee {h}_M')^{M\tau }:Y \to (D_{M}\times 
D_{M})^{M\tau }$ and $(h_M \vee h_M')^{M\nabla }:Y \to 
(D_{M}\times D_{M})^{M\nabla }$ by discretizing $h \vee h'$ as ${h_M} \vee {h_M'}$. 

\proclaim Proposition 2.4.12 (Copying Lemma). Let \Y $ \in End({\bf p})$  and suppose  \X $ \in 
End({\bf p})$  is ergodic. Let $f:X\to [\hbox{0},\hbox{1})$ and $g:Y\to [\hbox{0},\hbox{1})$ be generating. Suppose $f'$ is a tree-adapted function on $Y$  such that $dist({\bf X},f)=dist({\bf Y},f')$ . For all $\varepsilon $  and $N$, there exists a function $g'$ on $X$  such that
$$
dist(f_N \vee g'_N)^{N\nabla }\mathop \sim \limits^\varepsilon dist(f'_N \vee 
{g}_N)^{N\nabla }.
$$\par
\noindent{\twelvebx Proof:} By proposition 2.4.11, we have 
$$
dist({f'}_M^{M\nabla })=dist(f_M^{M\nabla }).
$$
\noindent Choose $M$ such that $N/M< \varepsilon /\hbox{2}$.  Construct a $\varepsilon /\hbox{2}$-tree Rokhlin tower ${\bf M}$ of height $M+\hbox{1}$ in 
${\bf X}$. Let $B$ be the base of ${\bf M}$. By the 
Strong Tree Rokhlin Lemma, we may assume that $f_M^{M\nabla }$ is independent 
of $B$. Hence, $dist(f_M^{M\nabla }|B)=dist({f'}_M^{M\nabla })$.  Let $\pi : {\bf Y} \to {\bf X}$ denote the tree adapted factor map defined by $\pi = (f^{{\bf N}^\ast})^{-1} \circ (f')^{{\bf N}^\ast}$.  

Fix an atom $q$ of $f_M^{M\nabla }$. Define a partition $P_{q}:q \cap B \to 
(D_{M}\times D_{M})^{M\nabla }$ such that
$$
dist(P_q ) = dist((f'_M \vee {g}_M)^{M\nabla }| \pi ^{ - 1} (q)).
$$
\noindent This gives a bijection $\varphi $ of the atoms of $P_{q}$ and 
those of $(f'_M \vee {g}_M)^{M\nabla } | \pi ^{ - 1} (q)$ such that 
corresponding atoms have the same conditional measures and are mapped to the 
same element in $(D_{M}\times D_{M})^{M\nabla }$.

For each atom  $\alpha \in P_{q}$, choose any representative point $y 
\in \varphi (\alpha )$. Note that as $f_M = f'_M \circ \pi$ and $\pi (y) \in q$, by proposition 2.1.4, $f_M^{M\nabla}(z)={f'}_M^{M\nabla }(y)$ for all $z \in \alpha $. For each $z \in \alpha 
$, choose a tree automorphism $A \in {\cal A}_{M}$ such that 
$f_{M}(T_{v}z)=f'_{M}(S_{Av}y)$ for each $v \in {\cal 
T}_{M}'$. We then define $g'$ by laying $\tau _y^{g} $ 
on $\{T^{-t}z \mid \hbox{1} \le t \le M \}$ via $A$. (We can ensure the 
measurability of $g'$ by choosing the same tree automorphism for all $z \in 
\alpha $ which are in the same atom of $f_M^{M\tau }$.) This defines the 
function $g'$ on the column over $\alpha $.

Clearly we have for every $z \in \alpha $ and $v \in {\cal 
T}_{M}'$,
$$
f_{M} \vee g'_{M}(T_{v}z)=f'_{M} \vee g_{M}(S_{Av}y),
$$
\noindent so that 
$$
(f_M \vee g'_M)^{M\nabla }(z)=(f'_M \vee g_M)^{M\nabla }(y)=P_{q}(z).
$$
\noindent Repeat the above procedure for each $\alpha  \in P_{q}$. We then have 
$(f_M \vee g'_M)^{M\nabla }=P_{q}$ on $q \cap B$. We thus have
$$
dist((f_M \vee g'_M)^{M\nabla }| q \cap B) = dist((f'_M \vee g_M)^{M\nabla 
}| \pi ^{ - 1} (q)) \eqno(\hbox{1})
$$
We then extend $g'$ to the rest of the tower by following the above procedure 
for each atom $q$ of $f_M^{M\nabla }$. We thus have a function $g'$ defined on ${\cup \bf M}\backslash B$. Extend $g'$ measurably to the rest of the space in any way 
we like.

We wish to prove that $dist(f_N \vee g'_N)^{N\nabla } \mathop \sim 
\limits^\varepsilon dist(f'_N \vee g_N)^{N\nabla }$. First, let us 
observe that by our construction and by (1), we have
$$
\eqalign{
 dist((f_M \vee g'_M)^{M\nabla }|B) &= \sum\limits_{q \in f_M^{M\nabla }} {\mu (q \cap B|B) dist((f_M \vee g'_M)^{M\nabla }|q \cap B)} \cr 
&= \sum\limits_{q \in f_M^{M\nabla }} {\mu (q) dist((f'_M \vee 
g_M)^{M\nabla }|\pi ^{ - 1} (q))}\cr 
 &= \sum\limits_{\pi ^{ - 1} (q) \in {f'}_M^{M\nabla }} {\nu (\pi ^{ - 1} (q)) 
dist((f'_M \vee g_M)^{M\nabla }|\pi ^{ - 1} (q))} \cr 
&= dist(f'_M \vee g_M)^{M\nabla }. \cr}
 \eqno(\hbox{2})
$$

We claim that for each $\hbox{0} \le j \le M-N$, 
$$ 
dist((f_N \vee g'_N)^{N\nabla }|T^{ - j}B) = dist(f'_N \vee g_N)^{N\nabla}.\eqno(\hbox{3})
$$
To see this, note that if $\beta $ and $\beta '$ are corresponding 
atoms of $(f_M \vee g'_M)^{M\nabla }$ and $(f'_M \vee g_M)^{M\nabla }$, 
then 
$$
dist((f_N \vee g'_N)^{N\nabla }|T^{ - j}\beta \cap T^{ - j}B) = dist((f'_N \vee 
g_N)^{N\nabla }|S^{ - j} \beta '). \eqno(\hbox{4})
$$
\noindent Moreover, by (2), we have
$$
\mu (T^{ - j}\beta |T^{ - j}B) = \mu (\beta | B ) = \nu 
(\beta ') = \nu (S^{ - j} \beta '). \eqno(\hbox{5})
$$
\noindent Thus (3) follows from (4) and (5) since the distributions in (3) are just 
weighted averages of the distributions equated in (4). As $N/M<\varepsilon 
/ \hbox{2}$ and the rest of the space is a set of measure $<\varepsilon / \hbox{2}$, we see 
that 
$$
dist(f_N \vee g'_N)^{N\nabla }\mathop \sim \limits^\varepsilon dist(f'_N \vee 
g_N)^{N\nabla }.
$$
\noindent This finishes the result.  {\symbol \char'04}
\medbreak
For \Y $ \in End({\bf p})$ and partitions $
P,R:Y \to C$ for a finite set $C$, we write $P\mathop \sim 
\limits^\varepsilon R$ if $\{ y \in Y \mid P(y) \ne R(y) \}$ has measure $
<\varepsilon $. For two partitions $P:Y \to C$ and $Q:Y \to C'$, we 
write $P\mathop \subset \limits^\varepsilon Q$ if there exists a partition 
$R:Y \to C$ with $R \le Q$ such that $P\mathop \sim \limits^\varepsilon R$. 

Given a function $g:Y \to \hbox{[0,1)}$, let $\sum (g)$ denote the 
sub-sigma-algebra of ${\cal C}$ generated by $g$, i.e. $\sum (g)=\bigvee 
\nolimits_{i = 0}^\infty S^{ - j}(g^{ - 1}({\cal D}))$ for the Borel 
sigma-algebra $\cal D$ of [0,1). If $g$ and $h$ are functions on $Y$, we write 
$h\mathop \subset \limits^{\varepsilon ,N} \sum (g)$ if there exists some 
integer $M$ such that $h_N \mathop \subset \limits^\varepsilon g_M^{M + }$. We 
write $h \subset \sum (g)$ if $h$ is $\sum (g)$-measurable. Clearly, $h 
\subset \sum (g)$ if and only if $h\mathop \subset \limits^{\varepsilon 
,N} \sum (g)$ for all $\varepsilon $ and $N$. 

For the remainder of this section, let \X and \Y be tvwB $\bf p$-endomorphisms. Moreover, we fix generating functions $f:X \to [\hbox{0},\hbox{1})$ and $g:Y\to [\hbox{0},\hbox{1})$. Before delving into details, we briefly describe the 
strategy to proving Theorem 2.4.1 from this point on. The idea is to 
construct a Cauchy sequence of functions $\{g^{j}\}$ on $Y$ such that 
$dist({\bf Y},g^{j})=dist({\bf X},f)$ for each $j$, and which converge 
to a function $\bar g$ on $Y$ such that $g  \subset \sum (\bar g)$. Thus 
$dist({\bf Y},\bar g)=dist({\bf X},f)$ and $\sum (\bar g)=\sum (g) = {\cal C}$. 
Now, the sigma-algebra ${\cal B}$ pulls backs to $\sum (g)$ via the factor 
map $\pi = (f ^{{\bf N}^\ast})^{-1} \circ {\bar g}^{{\bf N}^\ast} : {\bf Y} \to {\bf X}$ as constructed prior to Proposition 2.4.11. It thus follows that $\pi $ is an 
isomorphism as $\pi^{-1}$ gives an isomorphism of $\sigma$-algebras of Lebesgue spaces.  As we have previously remarked, $B^{ + }({\bf p})$ is tvwB. Theorem 2.4.1 thus follows. 

For technical reasons, we will choose the generators $f$ and $g$ such 
that $\{x \in X \mid f(x)=q\}$ is a null set for each rational $q$, and 
similarly for the sets $\{y \in Y \mid g(x)=q\}$. Observe that we 
can always do this as $X$ and $Y$ are nonatomic Lebesgue spaces. (Just choose any point isomorphism $f: X \to \hbox{[0,1)}$ and $g:Y \to \hbox{[0,1)}.$)  This assumption is needed to ensure that for the functions we construct via the copying lemma, a large enough set of points assume values sufficiently bounded away from the boundary points in the dyadic intervals of length $\hbox{1}/\hbox{2}^N$ in $P_N$. For any such function $h$, if $\left| \bar h - h \right|$ is sufficiently small, then $h_N(x)=\bar h_N(x)$ for all $x$ in a set of large measure.  

\proclaim Proposition 2.4.13. Let $f':Y \to [\hbox{0},\hbox{1})$ be a tree-adapted function such that $dist({\bf Y},f')= dist({\bf X},f).$
Then for any $\eta 
> \hbox{0}$, $\varepsilon > \hbox{0}$ and $M \in {\bf N}$, there exists a tree-adapted function $\bar f :Y\to [\hbox{0},\hbox{1})$  such that 
\itemitem {a)} $g \mathop \subset \limits^{\varepsilon ,M} \sum (\bar {f})$
\itemitem {b)} $\left| \bar f  - f' \right|<\eta $
\itemitem {c)} $dist({\bf Y},\bar f)=dist({\bf X},f)$.\par
\noindent{\twelvebx Proof:} Let $N \in {\bf N}$ be specified later. By Proposition 2.4.12, for any $\delta $ and $L$, we have some function $g'$ on $X$ such that 
$$
dist(f_L \vee g'_L)^{L\nabla }\mathop \sim \limits^\delta dist(f'_L \vee 
g_L)^{L\nabla }.\eqno(\hbox{1})
$$
\noindent As $f' \mathop \subset \sum (g)$, for each $\theta >\hbox{0}$, we can 
choose some $k$ such that
$$
f'_N \mathop \subset \limits^\theta {g_{k}^{k + }}.
$$
Now, by our choice of $g$, $\{y \in Y \mid g(y)=q \}$ is a 
null set for each rational $q$. Thus, for any $\beta >\hbox{0}$, there exist open intervals containing the rationals $\{ {t/ \hbox{2}^{k}} \mid t=\hbox{1},\ldots ,\hbox{2}^{k}-\hbox{1}\}$ such that $g(y)$ assumes a value in one of these intervals on a set of measure $<\beta$. Thus, (1) also implies that $g'(x)$ assumes a value in one of these intervals on a set of measure $< \hbox{2}\beta$ for all sufficiently small $\delta$ and $\hbox{1}/L$.  If $\beta > \hbox{0}$ is small enough, we may choose $\rho > \hbox{0}$ such that for any function $\bar g$ with $\left| g' - 
\bar g \right| <\rho $, $g'_{k} (x) = \bar {g}_{k} (x)$ on a 
set of sufficiently large measure to give 

$$
dist(f_{k} \vee g'_{k})^{k\nabla }\mathop \sim \limits^\theta dist(f_{k} \vee 
\bar {g}_{k})^{k\nabla }.\eqno(\hbox{2})
$$

\noindent From (1), $dist({g'}_L^{L\nabla })\mathop \sim \limits^\delta dist(g_L^{L\nabla})$. By the Strong Sinai's Theorem, if $\delta $ and $\hbox{1}/L$ are chosen small 
enough, then we may choose $\bar {g}$ with $\left| g'-\bar {g}\right|<\rho $, 
$dist({\bf X},\bar {g})= dist({\bf Y},g)$, and (2) holds.

Once more by proposition 2.4.12, for any $L'$ and $\delta '$, 
we have some function $\hat f$ on $Y$ such that
$$
dist(f_{L'} \vee \bar {g}_{L'})^{L'\nabla }\mathop \sim \limits^{\delta '} dist(\hat {f}_{L'} 
\vee g_{L'})^{L'\nabla }.\eqno(\hbox{3})
$$
\noindent Using the Strong Sinai's Theorem again, for any $\hbox{0}<\rho '<\theta$, we may choose $\delta '$ and $L'$ to give a tree-adapted function $\bar 
f$ on $Y$ such that $\left| \hat f  - \bar f \right|<\rho '<\theta $ 
with $dist({\bf Y},\bar f)=dist({\bf X},f)$. This gives (c) of the 
proposition. 

Choose $k' \in {\bf N}$ such that $\bar {g}_M \mathop 
\subset \limits^\varepsilon f_{k'}^{k'+}$. By (3) and decreasing 
$\delta '$ and $\hbox{1}/L' $ if needed, we have $g_M \mathop 
\subset \limits^\varepsilon \hat{f}_{k'}^{k' + }$. Once again, 
since $\{x \in X \mid f(x)=q\}$ is a null set for each rational $q$, we 
have $g_M \mathop \subset \limits^\varepsilon \bar {f}_{k'}^{k' + }$ by choosing $\rho' $, $\delta'$ and $\hbox{1}/L'$ 
small enough. This gives (a) of the proposition.

It remains to prove (b). Now, 
$$
\eqalign{
\left| f' - \bar {f} \right| &\le \left| f' - f'_N \right| + \left| f'_N - \hat {f}_N \right|  + 
\left| \hat {f}_N - \hat {f} \right| + \left| \hat {f} - \bar f \right|  \cr
& < \hbox{1} / \hbox{2}^N + \left| f'_N - \hat {f}_N \right| + \hbox{1} / \hbox{2}^N + \theta . \cr} 
$$
\noindent Using (1), (2) and (3) and ensuring $\hbox{min} (L,L') \ge k $, 
and $\delta $, $\delta '$ small enough we have
$$
dist(f'_{k} \vee g_k)^{k\nabla }\mathop \sim \limits^{2\theta } 
dist(\hat {f}_{k} \vee g_k)^{k\nabla }.
$$
\noindent Thus,
$$
dist(f'_{k} \vee g_k)^{k+}\mathop \sim \limits^{2\theta } 
dist(\hat {f}_k \vee g_k)^{k+}.
$$
\noindent Since $f'_N \mathop \subset \limits^\theta {g_k}^{k + }$, we may also conclude that $\left| {f'_N - \hat {f}_N } \right|<\hbox{3}\theta $. We may now conclude b) by choosing $\theta $ and $N$ in the 
beginning to satisfy $\hbox{4}\theta + \hbox{1}/{\hbox{2}^{N - 1}}<\eta $. {\symbol \char'04}
\medbreak
We are now ready to prove Theorem 2.4.1.
\medbreak
\noindent{\twelvebx Proof (Theorem 2.4.1):} The idea is to use Sinai's Theorem and 
Proposition 2.4.13 repeatedly to construct a Cauchy sequence of functions 
${\{g^{j}\}}$ on $Y$ converging to some function $\bar g$ pointwise a.e., and 
$dist({\bf Y},g^{j})=dist({\bf X},f)$ for each $j$, from which it follows that 
$dist({\bf Y},\bar g)=dist({\bf X},f)$. We now need to ensure that the 
functions $\{g^{j}\}$ are chosen in such a way that $g  \subset 
\sum (\bar g)$. From this, we may conclude that $\sum (\bar g)=\sum (g)$, 
thus concluding the proof of Theorem 2.4.1. 

To carry out the above plan, choose a sequence of reals $\varepsilon 
_{j} \searrow \hbox{0}$ and a sequence of integers $M_{j} \nearrow \infty 
$. Let $\eta _{j} \searrow \hbox{0}$ be a sequence of reals, to be specified 
later. Using Sinai's Theorem, we begin by choosing a tree-adapted function 
$g^{1}$ on $Y$ such that $dist({\bf Y},g^{1})=dist({\bf X},f)$. Using 
proposition 2.4.13, we construct a sequence of tree-adapted functions 
$\{g^{j}\}$ on $Y$, $j> \hbox{1}$, such that
\itemitem {i)} $dist({\bf Y},g^{j})=dist({\bf X},f)$, for $j \ge \hbox{1}$
\itemitem {ii)} $\left| g^{j} - g^{j + 1} \right|<\eta _{j }$, for $j \ge \hbox{1}$ 
\itemitem {iii)} $g\mathop \subset \limits^{\varepsilon _j ,M_j } \sum (g^j)$, for $j> \hbox{1}.$

\noindent By ii), we may obviously arrange the $\eta _{j}$'s such that the $g^{j}$'s 
converge to a function $\bar g$ pointwise a.e. Thus 
$dist({\bf Y},\bar g)=dist({\bf X},f)$. 
Since $g \mathop \subset \limits^{\varepsilon _j ,M_j } \sum 
(g^{j})$, for $j> \hbox{1}$, we have some integer $k_{j}$ such that
$$
g_{M_j } \mathop \subset \limits^{\varepsilon _j } {(g^j)_{k_j}^{k_j +}}.
$$
\noindent Now as $dist({\bf Y},g^{j})=dist({\bf X},f)$ for each $j$, the fact that 
$\{x \in X \mid f(x)=q\}$ is a null set for each rational $q$ implies 
$\{y \in Y \mid g^{j}(x)=q\}$ is also. Hence, we have some $\theta 
_{n}> \hbox{0}$ for each $n> \hbox{1}$ such that whenever
$$
\sum\limits_{j = n}^\infty  {\left| {g^j  - g^{j + 1} } \right|}  < \theta _n , \eqno(\hbox{1})
$$
\noindent we have
$$
g_{M_n } \mathop \subset \limits^{2\varepsilon _n } {\bar {g}_{k_n}^{k_n + }}.
$$

\noindent For any $m \in {\bf N}$ and $\varepsilon > \hbox{0}$, it follows that for all large n, 

\centerline{ $
{g \mathop \subset \limits^{\varepsilon ,m} {\sum  (\bar g)}.}
$}

\noindent We thus have $\sum (g)=\sum 
(\bar g)$, which gives Theorem 2.4.1.

Hence we are done if we can choose the $g^{j}$'s to satisfy the inequalities 
specified in (1) for all $n > \hbox{1}$. To see that this is possible, note that 
they are chosen in the order $g^{1} \to g^{2} \to g^{3} \ldots .$ In view of 
this, when $\theta _{n}$ is determined, only $g^{1},\ldots ,g^{n}$ have 
been chosen and $g^{n + 1}$ has not been chosen yet. Now, $\eta _{n}$ only need 
to be chosen when we choose $g^{n + 1}$. As a result, $\eta _{n}$ has not been declared at the time we specify $\theta _{n}$. 

Thus, we may choose $\eta _{n}< \hbox{min} ( {{\theta _1 }\over {2^n}},\ldots,{{\theta _n 
}\over {2}})$, it follows that
$$
\sum\limits_{j = n}^\infty {\left| {g^j - g^{j + 1}} \right|} <{ {\theta 
_n }\over {2}} +{ {\theta _n } \over {2^2}} + ... = \theta _n .
$$
\noindent Thus the inequalities in (1) can be satisfied. This completes the proof of 
Theorem 2.4.1. {\symbol \char'04}
\medjump
\noindent {\section 2.5. An Elementary Proof of a Special Case of Theorem}

{\quad ~ \section 2.4.1}
\medjump
In this section, we give a direct proof of Theorem 2.4.1 in the case that 
${\bf p}=(p_{1},\ldots ,p_{s})$ is a probability vector such that $p_{1}>\ldots 
>p_{s}$.  A simplified proof exists in this special case because the group $\cal A$ of tree automorphisms is trivial.   

Let \X $ \in End({\bf p})$. Let $I=\{ \hbox{1},\ldots ,s\}$. 
Consider the one-sided Bernoulli shift $B^{ + }({\bf p})$ represented as 
the shift space $(I^{{\bf N}^{\ast}},{\cal C},\nu,S)$ such that the $j$-th symbol 
has weight $p_{j}$. We then have a canonical factor map $\phi : {\bf X} \to B^{ 
+ }({\bf p})$ defined by setting $\phi (x)_k=j$ if $p_{\bf X}(T^{k}x)=p_{j}$, for 
$k \ge \hbox{0}$. 

Recall that a node of the ${\bf p}$-tree ${\cal T}$ is a finite sequence of 
integers in $\{ \hbox{1},\ldots ,s\}$. For each node $v$ and a point $z$ in $I^{{\bf N}^\ast 
}$, let $vz$ be the point in $I^{{\bf N}^\ast }$ obtained by concatenating $v$ to the 
left of $z$. For each $z$ in $I^{{\bf N}^\ast }$, let $z[\hbox{0},m]$ be the cylinder set 
$$
\{z'  \in I^{{\bf N}^\ast } \mid z'_j=z_j\; \hbox{for}\; \hbox{0} \le j \le 
m\}.
$$

Our goal is to prove the following special case of Theorem 2.4.1.
\proclaim Theorem 2.5.1. Let ${\bf p}$ be a probability vector with pairwise distinct components. Let \X $\in End({\bf p})$ be tvwB, then the canonical factor map $\phi :{\bf X} \to B^{ + }({\bf p})$ is an isomorphism.\par

To prove this, we need a preliminary lemma. Recall that for \X and \Y, if 
$\psi :{\bf X} \to {\bf Y}$ is a factor map, then we have fiber 
measures $\mu _{y}$ supported on $\psi ^{ - 1}(y)$ with the property 
that $\mu = \int {\mu _y \,d\nu (y)}$.

\proclaim Lemma 2.5.2.   Consider the one-sided Bernoulli shift $B^{ + }({\bf p})=(I^{{\bf N}^\ast },{\cal C},\nu,S)$. 
Let \X  $ \in End({\bf p})$. If $D \subseteq X$ and $v \in {\cal T}$, then for a.a. $z$ in $I^{{\bf N}^\ast }$, $\mu _{vz}(T_{v}D)=\mu _{z}(D)$ $($i.e. $T^{\left| v \right|}:(X,\mu 
_{vz}) \to (X,\mu _{z})$ is measure preserving$)$.\par

\noindent{\twelvebx Proof:} This follows easily from the fact that
$$\mu(T_v B | T_v C) = \mu (B | C)$$ 
for positive measurable sets $B \subseteq X$ and $C \subseteq X$. {\symbol \char'04}
\medbreak
\noindent{\twelvebx Proof (of Theorem 2.5.1):} Let $f:X \to \hbox{[0,1]}$ be generating. By 
assumption, $({\bf X},f)$ is tvwB. It suffices to prove that for 
a.a. $z$ in $I^{{\bf N}^\ast }$, there exists some $r \in \hbox{[0,1]}$ such that for $\mu 
_{z}$-a.a. $x$, $f(x)=r$ (i.e. $f$ is constant on fibres). Indeed, this implies that $f$ is $\phi 
^{-1}({\cal C})$-measurable. Since $f$ is generating, we thus have ${\cal 
B}=\phi ^{-1}({\cal C})$ and so $\phi $ is an isomorphism. 

We proceed by contradiction. Hence, we suppose that there exist $\eta >\hbox{0}$, a 
set $Z \subseteq I^{{\bf N}^\ast }$ of positive measure and two disjoint 
intervals $J$ and $J'$ in [0,1] separated by a distance of at least 
$\eta $ such that for each $z \in Z$, $\mu _{z}(x \in X \mid f(x) \in J)$ 
and $\mu _{z}(x \in X \mid f(x) \in J')$ are both at least $\eta $. 
We may clearly assume $\eta<\nu (Z)$. By the Tree Ergodic Theorem, there 
exists $N \in {\bf N}$ such that for all $n \ge N$, we have a set $K \subseteq I^{{\bf N}^\ast }$ of measure at least $\eta$ such that for each 
$z \in K$, 
$$
{\hbox{1} \over n}\sum\limits_{ \{0 < 
\left| v \right|  \le n \mid vz \in Z  \}} {w_v } \ge {{\nu (Z)}\over \hbox{2}} > 
{{\eta } \over \hbox{2}}. \eqno(\hbox{1})
$$
We now show that the integral
$$
\int {\bar {t}_n } (\tau _x^f ,\tau _y^f )\,d\mu _z (x)d\mu _z (y)
$$
\noindent is bounded away from zero for $z \in K$, $n \ge N$.
To see this, for $v \in {\cal T}$, by Lemma 2.5.2, 
$$
\eqalign{
\int {\left| \tau _x^f  (v)-\tau _y^f (v) \right| \,d\mu _z (x)d\mu _z (y)} &= \int 
{\left| f(T_v  x)-f(T_v y) \right| \,d\mu _z (x)d\mu _z (y)} \cr 
 &= \int {\left| f(x)-f(y) \right| \,d\mu _{vz} (x)d\mu _{vz} (y)}. \cr} 
$$
\noindent Hence, if $z \in K$ for $n \ge N$, then from (1) and the fact that $\cal A$ consists of only the identity automorphism,
$$
\eqalign{
 \int {\bar {t}_n (\tau _x^f , \tau _y^f )\,d\mu _z (x)d\mu _z (y)} &= 
{\hbox{1} \over n}\sum\limits_{0 < \left| v\right|  \le n} {w_v \int {\left| \tau _x^f  (v) - \tau _y^f (v) \right| \,d\mu _z (x)d\mu _z (y)}}\cr 
&\ge {\hbox{1} \over n}\sum\limits_{\{0 < \left| v \right|  \le n \mid vz \in Z \}} {w_v \int {\left| f(T_v x)-f(T_v y) \right| \,d\mu _z (x)d\mu _z (y)}} \cr 
& = {\hbox{1} \over n}\sum\limits_{\{0 < \left| v\right|  \le n \mid vz \in Z \} } {w_v  \int \left| f(x)-f(y) \right| \,d\mu _{vz} (x)d\mu _{vz} (y)} \cr 
&\ge {\hbox{1} \over n}\sum\limits_{\{0 < \left| v\right|  \le n \mid vz \in Z \} }{w_v  \eta ^{3}} \ge {\eta ^4 }/ \hbox{2}, \cr} 
$$
\noindent where the second last inequality follows from the fact for each $vz$ in $Z$, we have two disjoint sets $U$ and $U'$ with $\mu _{vz}$ measures at least 
$\eta $ and $\left| f(x)-f(y) \right| \ge \eta$ whenever $x \in U$ and $y \in 
U' $. Hence, 
$$
\int {\bar {t}_n } (\tau _x^f ,\tau _y^f )\,d\mu _z (x)d\mu _z (y) \ge {\eta ^4 
/  \hbox{2}}
$$
\noindent for all $z \in K$, $n \ge N$.

We will now use the fact that $({\bf X},f)$ is tree very weak Bernoulli to 
arrive at a contradiction. Let $\hbox{0}<\varepsilon < \hbox{1}$ be specified later. Choose 
$n$ corresponding to $\varepsilon $ in the definition of tvwB for 
$({\bf X},f)$. We may assume that $n \ge N$ (for the $N$ chosen in the last paragraph). By 
definition, we have a set $G$ of measure $>\hbox{1}-\varepsilon 
$ such that for all $x$ and $y$ in $G$, $\bar {t}_n (\tau _x^f 
,\tau _y^f ) \le \varepsilon $. Now, we have
$$
\mu (G^c) = \int { \mu_z ({G^c})\;d\nu (z)} \le \varepsilon .
$$
\noindent Hence, 
$$
\mu_z (G^c )  \le \sqrt \varepsilon 
$$
\noindent for all $z$ in a set $G' \subseteq I^{{\bf N}^\ast }$ of measure $\ge \hbox{1} 
- \sqrt \varepsilon$. Consequently, we see that for each $z$ in $G'$,
$$
\int {\bar {t}_n } (\tau _x^f ,\tau _y^f )\,d\mu _z (x)d\mu _z (y) \le 
\varepsilon + \hbox{2} \sqrt \varepsilon .
$$
Choosing $\varepsilon > \hbox{0}$ sufficiently small ensures that $\nu (K \cap 
G')> \hbox{0}$. Moreover, choose $\varepsilon $ such that $\varepsilon +  
 \hbox{2} \sqrt \varepsilon <{\eta ^{4}}/ \hbox{2}$. Now, if $z \in K \cap 
G'$, we have on the one hand,
$$
\int {\bar {t}_n } (\tau _x^f ,\tau _y^f )\,d\mu _z (x)d\mu _z (y) < {\eta ^4} / 
 \hbox{2}
$$
\noindent as $z \in G'$, while 
$$
\int {\bar {t}_n } (\tau _x^f ,\tau _y^f )\,d\mu _z (x)d\mu _z (y) \ge \eta ^4 
/  \hbox{2}
$$
\noindent as $z \in K$. This is a contradiction and thus completes the proof of the 
theorem. {\symbol \char'04}
\vfill \eject

$ $
\vskip 5.90cm
\noindent {\chapter Chapter 3: }
\medjump
\noindent{\chapter A Joinings Proof of the Isomorphism} 
\medjump
\noindent{\chapter Theorem}
\bigjump
\twelverm
The goal of this chapter is to present a joinings proof of the isomorphism 
theorem, Theorem 2.4.1, in the previous chapter. The proof is a modification 
of Hoffman and Rudolph's ([H,R]) proof of the isomorphism theorem in the case 
of the uniform probability vector. While the joinings proof is more 
technical in certain aspects than the proof presented in the previous chapter, it has the 
advantage of showing that there are uncountably many automorphisms of \Bp , unless the components of the probability vector ${\bf p}$ are 
pairwise distinct (see Proposition 3.4.2). Most of the definitions and 
theorems below are modelled after Hoffman and Rudolph. There are, however, 
two modifications that we will make to their proof which will allow us to 
extend their arguments to the general probability vector. First, we will 
need to extend the definition of one-sided joinings introduced in [H,R]. 
In particular, we need to impose an additional condition in the definition 
of a one-sided joining which is trivial when ${\bf p}$ is uniform. 
Second, our statement of the copying lemma will differ from that in Hoffman 
and Rudolph ([H,R]'s Lemma 5.4). In our presentation, we will need to copy tree 
distributions induced by the functions on the $\bf p$-endomorphisms under consideration. This approach 
will save us from tackling the technical issue of whether dist and tdist 
generate the same topology on tree-adapted functions on ${\bf 
p}$-endomorphisms for a general ${\bf p}$ ([H,R]'s Lemma 3.6). 

This chapter is organized into four sections. In \S 3.1, we define the 
notion of a one-sided joining of two ${\bf p}$-endomorphisms and discuss some 
topological properties of such joinings. \S 3.2 is devoted to the proof of 
the copying lemma, which is similar in form to Proposition 2.4.12 but it is 
more involved. \S 3.3 examines the $\bar 
{t}$ distance between two processes, as defined in [H,R]. \S 3.4 contains the proof of Theorem 
2.4.1, using the machinery of one-sided joinings developed in \S 3.1 to \S 
3.3. The idea is to show that for a tvwB ${\bf p}$-endomorphism 
\X and for $B^+({\bf p})=(Y,{\cal C},\nu,S),$ the set 
of one-sided joinings $\lambda$ such that ${\cal B}\mathop =
\limits^\lambda {\cal C}$ (the isomorphic joinings) is a dense $G_{\delta }$ in the space of one-sided 
joinings, which is a compact metric space in the $w^{\ast}$-topology. The 
Baire Category Theorem then implies that the set of isomorphic joinings is non-empty, provided that the space of one-sided joinings is non-empty (which will indeed be the case). However, each such joining gives an isomorphism and so ${\bf X} \cong 
B^+({\bf p})$. 
\medjump
\noindent {\section 3.1. One-sided Joinings }
\medjump
Let ${\bf p}=(p_{1},\ldots ,p_{s})$ be a fixed probability vector such 
that $p_{1} \ge \ldots \ge p_{s}$. Let $p_{0}=\hbox{1}$ and assume that 
$\hbox{0}=s_{0}<s_{1}<\ldots <s_{r}=s$ are chosen such that for all $\hbox{0} \le i 
\le r-\hbox{1}$, $p_{s_i + 1} =\ldots =p_{s_{i + 1} }$ and $p_{s_i } 
>p_{s_{i + 1} }$. Define the probability vector $\bar {\bf p}$ by 
summing the identical components of ${\bf p}$, i.e.
$$
\bar {\bf p} = (\sum\limits_{i = 1}^{s_1 } {p_i},\ldots,\sum\limits_{i = s_{r - 1} 
+ 1}^{s_r } {p_i } ) 
$$
Let $I(\bar {\bf p})=\{ \hbox{1},\ldots ,r\}^{{\bf N}^\ast }$ and assign $j$ with weight equal to the $j$-th component of $\bar {\bf p}$.  Construct the one-sided Bernoulli shift $B^{ + }(\bar {\bf p})=(I(\bar 
{\bf p}),{\bf m},\sigma )$ where ${\bf m}$ is the product measure and $\sigma$ is the 
shift defined in the usual way. Given \X $ \in 
End({\bf p})$, define a map $\psi _{\bf X}:X \to I(\bar {\bf p})$ by setting 
$\psi _{\bf X}(x)_t=j$ if $p_{\bf X}(T^t x)=p_{s_j}.$  Note that if ${\bf p}$ is a uniform 
probability vector, $\psi _{\bf X}$ is a constant map into a single point 
system. 
For two m.p.s. \X and \Y, a {\twelvebx coupling} of $\bf X$ and $\bf Y$ is a measure $\lambda $ on the product 
space $(X\times Y,{\cal B}\times {\cal C})$ such that $\lambda 
(B\times Y)=\mu (B)$ for all $B \subseteq X$ and $\lambda (X\times 
C)=\nu (C)$ for all $C \subseteq Y$ (i.e. $\lambda $ has marginals $\mu $ 
and $\nu$). A {\twelvebx joining} of ${\bf X}$ and ${\bf Y}$ is a coupling $\lambda$ which 
is also $T\times S$-invariant. For brevity, for a set $B \subseteq X$, we 
will also use $B$ to denote the subset $B\times Y$ in the product space 
$X\times Y$, with similar convention for a set $C \subseteq Y$. (It will 
always be clear from the context whether $B \subseteq X$ refers to $B$ or 
$B \times Y$.) 

Two examples of joinings are worth mentioning at this stage. First, the 
product measure $\mu \times \nu $ is always a joining of ${\bf X}$ and 
${\bf Y}$. Second,  the diagonal measure $\chi 
_{\Delta }$ on $X\times X$ defined by $\chi _{\Delta }(B \times C)=\mu 
(B \cap C)$ for any $B \subseteq X$ and $C \subseteq X$ is a 
self-joining of ${\bf X}$. 

\proclaim Definition 3.1.1. For \X and \Y in $End({\bf p})$, a {\twelvebx one-sided coupling (joining)} $\lambda$ of  ${\bf X}$ and ${\bf Y}$ is a coupling $($joining$)$ of ${\bf X}$ and ${\bf Y}$ such that
 \item{i)} for any pair of compact-valued generating functions $f:X \to R$ and $g:Y \to U$, if $f^i = f \circ T^{i}$ and $g^i = g \circ S^{i}$, then for each $j \ge \hbox{0}$,  we have 
$$\eqalignno{
Dist_{\lambda}(\{f^i\}_{0 \le i \le j} | \{f^i\}_{j < i} ,\{g^i\}_{j < i}) &= Dist_{\lambda}(\{f^i\}_{0 \le i \le j} | \{f^i\}_{j < i} )\cr
\noalign  {\noindent \hskip 0.7cm \hbox{and}}
Dist_{\lambda}(\{g^i\}_{0 \le i \le j} | \{f^i\}_{j < i} ,\{g^i\}_{j < i}) &= Dist_{\lambda}(\{g^i\}_{0 \le i \le j} | \{g^i\}_{j < i})\cr}$$

\noindent \item {ii)} $\lambda $ projects to the diagonal measure on $I(\bar {\bf p})\times I(\bar {\bf p})$ via the map $\psi _{\bf X}\times \psi _{\bf Y}:X\times 
Y \to  I(\bar {\bf p})\times I(\bar {\bf p})$.\par

\medbreak
Some remarks on definition 3.1.1 are in order. First, condition i) holds if 
the equations hold for {\twelvebx some} pair of generating functions.  Indeed, note that the equalities in i) hold if and only if 
$$\eqalignno{
Dist_{\lambda}(\{f^i\}_{ i \ge 0} | \{f^i\}_{j < i} ,\{g^i\}_{j < i}) &= Dist_{\lambda}(\{f^i\}_{i \ge 0} | \{f^i\}_{j < i} )\cr
\noalign  {\noindent \hbox{and}}
Dist_{\lambda}(\{g^i\}_{i \ge 0} | \{f^i\}_{j < i} ,\{g^i\}_{j < i}) &= Dist_{\lambda}(\{g^i\}_{i \ge 0} | \{g^i\}_{j < i}).\cr}$$  
Since the sigma-algebras generated by $\{f^i\}_{i \ge j}$ and $\{g^i\}_{i \ge j}$ for any $j \ge \hbox{0}$ do not depend on the choice of the generators, condition i) is just 
the Hoffman and Rudolph definition ([H,R]'s Definition 3.2). In the case that $\lambda$ is $T \times S$-invariant, condition i) essentially says that for $\lambda$-a.a $(x,y)$, if $x'$ is a preimage of $x$, then the conditional probability of $x'$ given $(x,y)$ is just the conditional probability of $x'$ given only $x$.  Condition ii) essentially says that for 
$\lambda \hbox{-}a.a$ $(x,y)$, the $p_{\bf X}\hbox{-name of}\; x = p_{\bf Y}\hbox{-name of}\; y$. Thus if ${\bf p}$ is a 
uniform probability vector, condition ii) always holds and Definition 3.1.1 
is just the Hoffman and Rudolph definition.  

For brevity, given ${\bf X}, {\bf Y} \in End({\bf p})$, let $C({\bf X},{\bf Y})$ and $C^{ + }({\bf X},{\bf Y})$ denote the set of 
couplings and one-sided couplings of ${\bf X}$ and ${\bf Y}$. Similarly, let 
$J({\bf X},{\bf Y})$ and $J^{ + }({\bf X},{\bf Y})$ denote the set of joinings 
and one-sided joinings of ${\bf X}$ and ${\bf Y}$. By viewing the set of 
probability measures on $(X\times Y,{\cal B}\times {\cal C})$ as a 
subset of bounded linear functionals on the continuous functions on $X\times 
Y$, we have a metrizable topology, the $w^{\ast}$-topology, on the set of couplings 
and joinings. It is standard that $C({\bf X},{\bf Y})$ and $J({\bf X},{\bf Y})$ 
are $w^{\ast}$-compact. 

We will now establish some basic facts of one-sided 
couplings and joinings which will be used in subsequent sections.  

\proclaim Proposition 3.1.2. $C^{ + }({\bf X},{\bf Y})$ and $J^{ + }({\bf 
X},{\bf Y})$ are $w^{\ast}$-closed convex subsets of $C({\bf X},{\bf Y})$ and $J({\bf X},{\bf Y})$ respectively.\par 

\noindent{\twelvebx Proof:} The statement that $J^{ + }({\bf X},{\bf Y})$ is a $w^{\ast}$-closed 
convex subset of $J({\bf X},{\bf Y})$ easily follows from the corresponding 
statement for $C^{ + }({\bf X},{\bf Y})$ and $C({\bf X},{\bf Y})$ and the fact 
that the set of joinings is a closed convex subset of all couplings. Hence, 
it is enough to prove that $C^{ + }({\bf X},{\bf Y})$ is a $w^{\ast}$-closed convex 
subset of $C({\bf X},{\bf Y})$. 

For convexity, notice that condition i) of Definition 3.1.1 says that for each $t \ge \hbox{0}$, each one-sided coupling couples ${\cal B}$ and $S^{-t}{\cal C}$ independently when conditioned on $T^{-t}{\cal B}$.  As each one-sided coupling projects to $\mu$ on $\cal B$, this implies that conditionally on $T^{-t}{\cal B}$, each one-sided coupling can be viewed as a product measure of the form $\mu \times \nu_i$ on ${\cal B} \times S^{-t}{\cal C}$.  It follows that any convex combinations of one-sided couplings is still a product measure of this form on ${\cal B} \times S^{-t}{\cal C}$ conditionally on $T^{-t}{\cal B}$.  Thus any convex combinations couples $\cal B$ and $S^{-t}{\cal C}$ independently over $T^{-t}{\cal B}$.  By symmetry, ${\cal C}$ and $T^{-t}{\cal B}$ are also coupled independently when conditioned on $S^{-t}{\cal C}$.  Thus, condition i) holds for convex combinations.  The fact that condition ii) of Definition 3.1.1 also holds for convex combinations of one-sided couplings is obvious since any convex combination of measures projecting to the diagonal measure on $I(\bar {\bf p})\times I(\bar {\bf p})$ also projects to the diagonal measure.  Hence, one-sided couplings are convex.  

We prove closure. Fix $t \in {\bf N}$. Let $\lambda _{i} \in C^{ + 
}({\bf X},{\bf Y})$ be a sequence of one-sided couplings such that $\lambda 
_{i} \to \lambda $ in $w^{\ast}$. Clearly, $\lambda  \in C({\bf X},{\bf 
Y})$. Suppose $P$ is a finite partition of $X$, and $Q$ is a finite partition of $Y$ 
which is $S^{ - t}{\cal C}$-measurable. Choose finite partitions 
$R_{j}$ of $X$ such that $ {R_j} \nearrow T^{-t}{\cal B}$. 
For condition i) of definition 3.1.1, it suffices to show that for each atom $\alpha \in P$, 
$$
E_\lambda (\alpha|{T^{ - t}{\cal B} \vee Q}) = E_\lambda (\alpha| 
{T^{ - t}{\cal B}}) \quad \hbox{a.e.}
$$
\noindent Indeed, this implies conditionally on $T^{-t}{\cal B}$, $\lambda $ couples $P$ and $Q$ independently. 
Since $Q$ and $P$ are arbitrary, condition i) of Definition 3.1.1 follows.

Note that for each $s \in {\bf N}$ and $\alpha \in P$, 
\vfill \eject
$$E_{\lambda _i } (\alpha|\mathop \vee 
\limits_1^s {R_j \vee Q}) \to E_\lambda (\alpha|\mathop \vee \limits_1^s {R_j 
\vee Q}) \quad \hbox{a.e.} $$
\noindent and
$$
E_{\lambda _i } (\alpha|\mathop \vee \limits_1^s R_j \vee Q) \ge E_{\lambda 
_i } (\alpha|T^{ - t}{\cal B} \vee Q) = E_\mu (\alpha|T^{ - t}{\cal B}) \quad \hbox{a.e.}
$$
\noindent as $\lambda _{i }$ is one-sided. Thus, we have $E_\lambda (\alpha|T^{ - 
t}{\cal B} \vee Q) \ge E_\mu (\alpha|T^{ - t}{\cal B})$ a.e. by the Martingale 
Convergence Theorem. Since $T^{-t}{\cal B}$ is a sub-$\sigma$-algebra 
of $T^{-t}{\cal B} \vee Q$, the reverse inequality also holds so that 
condition i) of Definition 3.1.1 is proved. For condition ii), note 
that for each cylinder set $c \subseteq I(\bar {\bf p})$, $(\psi _{\bf X}\times 
\psi _{\bf Y})^{-1}(c\times c)$ is a set of the form $B\times C$, 
where $B \subseteq X$ and $C \subseteq Y$. Since $\lambda _{i} \to 
\lambda $ in $w^{\ast}$,  we have $\lambda _{i}(B\times C) \to \lambda 
(B\times C)$.  Hence, as 
$$
\lambda _{i}((\psi _{\bf X}\times 
\psi _{\bf Y})^{-1}(c\times c))={\bf m}(c),
$$
the same holds for $\lambda $ so that $\lambda$ also projects to the diagonal measure on 
$I(\bar {\bf p})\times I(\bar {\bf p})$. {\closeproof}
\medbreak

The following proposition, which extends Lemma 3.7 in [H,R] to $\bf p$-endomorphisms, gives an example of a one-sided joining.  Given a factor map $\phi :{\bf X} \to {\bf Y}$, $\phi $ yields a 
probability measure $\lambda _{\phi }$ on $X\times Y$ defined by $\lambda 
_{\phi }(B\times C)=\mu (B \cap \phi ^{ - 1}C)$ for $B \in 
{\cal B}$ and $C \in {\cal C}$. It is easy to check that $\lambda _{\phi 
} \in J({\bf X},{\bf Y})$. We say that $\lambda _{\phi }$ is the 
{\twelvebx graphical joining} induced by $\phi $. 
\proclaim Proposition 3.1.3. Suppose \X and \Y are ergodic $\bf p$-endomor-phisms. Let $\phi :{\bf X} \to {\bf Y}$ be a tree-adapted factor map. Then $\lambda _{\phi } 
\in$ \joinXY. \par
\noindent {\twelvebx Proof:}  Note that for any $t \in {\bf N}$ and any node $u$ of length $t$, the 
conditional probability of $T_{u}x$ given $x$ and $\phi (x)$ equals that of $T_{u}x$ 
given $x$, as $\lambda_{\phi}$ is the graphical joining arising from the factor map $\phi : X \to Y$.  On the other hand, note that as 
$\phi $ is tree-adapted, for any $x \in \phi ^{-1}(y)$, there is a 
unique node $v$ of length $t$ such that $\phi (T_{v}x)=S_{u}y$. Clearly, the 
conditional probability of $S_{u}y$ given $x$ and $y$ equals that of $T_{v}x$ 
given $x$. By Proposition 2.1.3, since $\phi $ is a tree adapted factor map, 
$w_{u}=w_{v}$. Thus, the conditional probability of $S_{u}y$ given $x$ and $y$ 
equals that of $S_{u}y$ given $y$ (in fact, both are equal to $w_{u}$). Since 
$\lambda _{\phi }$ is stationary, it follows that condition i) of 
Definition 3.1.1 of one-sided joinings holds for $\lambda _{\phi 
}$. 

By Proposition 2.1.3, for any cylinder set $c$ in $I(\bar {\bf p})$, $\phi 
^{-1}(\psi _{\bf Y}^{ - 1} (c))=\psi _{\bf X}^{ - 1} (c)$.  Thus, 
$$\lambda _\phi (\psi _{\bf X}^{ - 1} (c)\times \psi _{\bf Y}^{ - 1} (c)) = 
\mu (\psi _{\bf X}^{ - 1} (c)) = {\bf m}(c),$$
\noindent This implies condition ii) of definition 3.1.1. Hence $\lambda 
_{\phi } \in$ \joinXY. {\closeproof} 
\medbreak
Note that Proposition 3.1.3 implies that if \X is a $\bf p$-endomorphisms, then $J^+ ({\bf X}, B^+({\bf p}))$ is non-empty.  Indeed, we can easily construct a tree-adapted factor map $\phi_1 : {\bf X} \to B^+({\bf p})$ by choosing a tree partition $K_{\bf X}$ of $X$ and mapping a point $x$ to its $K_{\bf X}$-name.  By Proposition 3.1.3, the graphical joining arising from this factor map gives a one-sided joining.  

\medjump
\noindent{\section 3.2. The Copying Lemma}
\medjump
In this section, we will establish the copying lemma for one-sided joinings 
which will be the key ingredient in the proof of Theorem 2.4.1.

Let ${\bf X}$ and ${\bf Y}$ be ${\bf p}$-endomorphisms and let $\lambda  \in 
J^{ + }({\bf X},{\bf Y})$. For $(x,y) \in X\times Y$, and a pair of 
nodes $(v,u)$ of the same length $j$, let $p_{x,y}(v,u)$ be 
the conditional mass of $(T_{v}x,S_{u}y)$ given $(x,y)$. By condition 
i) of Definition 3.1.1, we have 
$\sum _{\left| u \right| = j} 
{p_{x,y} (v,u)}=w_{v}$ and $\sum _{\left| v \right| = j} 
{p_{x,y} (v,u)}=w_{u}$ for $\lambda$-a.a. $(x,y)$. Moreover, if $w_v \ne w_{u}$, then $p_{x,y}(v,u)=\hbox{0}$ for 
$\lambda $-a.a. $(x,y)$ by condition ii) of Definition 3.1.1. For a 
tree automorphism $A$ and two nodes $v$ and $u$, let $A(v,u) 
=w_{v}$ if $u=Av$, and $A(v,u)=\hbox{0}$ otherwise. 

Before we state the next proposition, we recall the following well-known 
theorem which states that any doubly stochastic matrix (i.e. a square matrix such that the entries in each row and column sum to one) is expressible as an average of permutation matrices (i.e. square matrices such that each row and column consists of precisely a single entry of 1 with the rest 0's).  The proof is based on Hall's Marriage Lemma and can be found in many standard combinatorial texts (eg. [Ryd]).
\proclaim Proposition 3.2.1. Let $M$ be a doubly stochastic matrix $of$ order $n$.  Then $M$ can be expressed as a convex combination of  permutation matrices, i.e. there exist nonnegative reals $c_1, \ldots, c_t$ which sum to one, along with permutation matrices $P_1, \ldots, P_t$ such that 
$$
M=c_1 P_1 + \cdots + c_t P_t .
$$\par
\noindent {\twelvebx Remark:} If in Proposition 3.2.1, $M$ is a matrix such that the entries in each row and column sum to some fixed number $\alpha$, then the conclusion still holds with $P_1, \ldots, P_t$ replaced by ``permutation" matrices where each nonzero entry is $\alpha$.
\medbreak
\proclaim Proposition 3.2.2. Let $\lambda  \in J^{ + }({\bf X},{\bf Y})$ and $N \in {\bf N}$. Then for $\lambda $-a.a. $(x,y) \in X\times Y$, we have a probability measure $m_{x,y}$ on ${\cal A}_{N}$  such that for any $\hbox{1}\le j \le N$, if $\left| v \right| =\left| u \right| = j$, then 
$$
p_{x,y} (v,u) = \int {A(v,u)\,dm_{x,y} (A)} .
$$ \par
\noindent{\twelvebx Proof:} Fix $(x,y) \in X\times Y$. We prove the result by induction on $N$. First, 
suppose $N=\hbox{1}$. Let $V$ be the set of nodes of length 1. Define a measure $\rho $ 
on $V\times V$ by setting
$$
\rho (U) = \sum\limits_{(v,u) \in U} {p_{x,y} (v,u)} 
$$
for $U\subseteq V \times V$.
We wish to apply Proposition 3.2.1 to construct a measure on ${\cal A}_1$ from $\rho$.  To do this, recall that each node of length 1 is an integer in $\{ \hbox{1},\ldots,s \}$.  We can represent the measure $\rho$ as a $s \times s$ matrix $M$ whose columns are indexed by ${\hbox{1},\ldots,s}$ such that $M_{vu}=\rho (v,u)$.  By the one-sidedness of $\lambda$, note that $M$ is a block diagonal matrix such that for each block, the entries in each row and column have the same sum (in fact, the sum for the $j$-th block is just  $p_{s_j}$).  For each of the blocks $M_j$, we can apply the remark following Proposition 3.2.1 and express it as a convex combination of permutation matrices.  Doing this block by block, gives a decomposition of $M$ in the form
$$
M=\sum_{n=1}^{t} {a_n Q_n}, \eqno(\hbox{1})
$$
where the $a_n$'s sum to 1, $Q_n$'s are distinct block diagonal matrices such that for each $Q_n$, the $j$-th block is a permutation matrix with each nonzero entry being $p_{s_j}$, and the sum is taken over all such possible matrices.  We may then use this decomposition on $M$ to define a measure on ${\cal A}_1$ as follows.  For $A \in {\cal A}_1$, choose the unique matrix $Q_n$ such that $(Q_n)_{vu}=A(v,u)$ for all $v, u \in \{ \hbox{1},\ldots,s\}$, and set ${m}(A)=a_n$.  Then (1) immediately implies that 
$$
p_{x,y} (v,u) = \rho (v,u) = M_{vu}= \int {A(v,u)dm(A)} .
$$

Next, assume that the result holds for $N=t$.  We wish to build a measure on ${\cal A}_{t + 1}$ such that the asserted equality in the statement of the proposition continues to hold 
for nodes of length $ \le t+\hbox{1}$. Notice that each tree automorphism in ${\cal A}_{t + 1}$ is defined by a tree automorphism  in ${\cal A}_{t }$ combined with a collection of tree automorphisms in ${\cal A}_1$ indexed by the nodes of length $t$.  To define the required measure ${m}_{t+1}$ on ${\cal A}_{t + 1}$, we proceed as follows.  Fix $A' \in {\cal A}_{t + 1}$ and let $A$ denote the restriction of $A'$ to ${\cal T}_t$.  For each node $\left| v \right| =t$, consider the pair of points 
$(T_{v}x,S_{Av}y)$. Using the basis case, we have a measure ${m}_v$ on ${\cal A} _1$ such that for any $\left| u \right| =\left| u' \right| = \hbox{1}$,
$$
p_{T_v x,S_{Av} y} (u,u') = \int {\bar {A}(u,u')dm_{v} (\bar {A})}. 
$$
We then define
$$
m(A') = \prod\limits_{\left| v \right| = t} {m_{v} (B_{v} )m(A)},
$$
where $B_{v}$ is the tree automorphism in  ${\cal A}_{1}$ induced by $A'$ on the trees of height one rooted at $v$ and $Av$. A simple calculation using the basis case and the induction hypothesis when $N=t$ shows that for every pair of nodes $(v,u)$ of 
common length $ \le t+\hbox{1}$
$$
p_{x,y} {(v,u)} = \int\limits_{A' \in {{\cal A}_{t + 1}} } {A'(v,u)\,dm(A')}. 
$$
{\closeproof}
\medbreak
For a finite set $C$ and partitions $P:X \to C$ and $Q:Y \to C$, define 
the joined partition $P \otimes Q:X\times Y \to C\times C$ by $P 
\otimes Q(x,y)=(P(x),Q(y))$. If $\lambda  \in J({\bf 
X},{\bf Y})$, then $\lambda $ induces a stationary measure $\lambda _{P 
\otimes Q}$ on the shift space $(C\times C)^{{\bf N}^\ast }$ via the map $(x,y) 
\to  (P \otimes Q(T^{i}x,S^{i}y))_{i \ge 0}$. Clearly, we can extend 
$\lambda _{P \otimes Q}$ to a stationary measure on $(C\times C)^{\bf Z}$ 
and then restrict it to a measure on $(C\times C)^{ - \bf {N}}$. 

For any $N \in {\bf N}$ and $j \ge N$, and a pair of elements $(\alpha ,\beta )$ in 
$C^{N}$, let $(\alpha \times \beta)^{-j}$ denote the cylinder set
$$\{z \in (C\times C)^{ - \bf{N}} \mid z_t=(\alpha_{t+j+1},\beta_{t+j+ 1}), 
-j \le t \le  -j+N-\hbox{1}\}.$$ 
Given $A \in {\cal A}$ and tree names $h,h': {\cal T}' \to C$, we then define a measure $\lambda _{(h,h',A)} $ on the 
cylinder sets by setting 
$\lambda _{(h,h',A)} ((\alpha \times \beta)^{-j})$ to be the total weights of all nodes $\left| v \right| = j$ such that $$(h(v),h(\sigma (v)), \ldots, h(\sigma^{N-1}(v)))=\alpha$$ 
and 
$$(h'(Av),h'(\sigma (Av)), \ldots, h'(\sigma^{N-1}(Av)))=\beta.$$ 
The following proposition states that we can represent $\lambda _{P 
\otimes Q}$ as an average of measures of the form $\lambda _{(\tau _x^P 
,\tau _y^Q ,A)} $. This will be used to prove the copying lemma. 
\proclaim Proposition 3.2.3. With the notations above, for $\lambda \in J^{ + }({\bf X},{\bf Y})$ and partitions $P:X \to C$ and $Q:Y \to C$ for a finite set $C$, we have for each $N \in {\bf N}$, a family of probability measures $m_{x,y}$ on ${\cal 
A}_{N}$ such that 
$$
\lambda _{P \otimes Q} ((\alpha \times \beta)^{-N}) = \int {\int {\lambda _{(\tau _x^P 
,\tau _y^Q ,A)} ((\alpha \times \beta)^{-N})\,dm_{x,y} (A)d\lambda (x,y)}}
$$
\noindent for each pair of elements $(\alpha ,\beta )$  in $C^{N}$.\par
\noindent {\twelvebx Proof:} Abbreviate the cylinder set $(\alpha \times \beta)^{-N}$ as $(\alpha \times \beta)$.  Notice that as $\lambda $ is $T\times S$-invariant, we have
$$
\lambda _{P \otimes Q} (\alpha \times \beta ) = \int {\sum\limits_{  {\left| v \right| = \left| u \right| = N}\atop {P^{\raise1pt \hbox{$\scriptscriptstyle N+$}}(T_v x)=\alpha, Q^{\raise1pt \hbox{$\scriptscriptstyle N+$}}(S_{u} y)=\beta} } p_{x,y}(v,u)\, d\lambda  (x,y)}
$$
By Proposition 3.2.2, we have a measure $m_{x,y}$ on ${\cal A}_N$ such that 
$$
p_{x,y} (v,u) = \int {A(v,u)\,dm_{x,y} (A)} 
$$
\noindent
for all pairs of nodes $v$ and $u$ of length $N$. 
\noindent We then have 
$$
\eqalign{
 \lambda _{P \otimes Q} (\alpha \times \beta ) &= \int \sum\limits_{  {\left| v \right| = \left| u \right| = N}\atop {P^{\raise1pt \hbox{$\scriptscriptstyle N+$}}(T_v x)=\alpha, Q^{\raise1pt \hbox{$\scriptscriptstyle N+$}}(S_{u} y)=\beta} } {p_{x,y} (v,u)\,d \lambda  (x,y)} \cr 
 &= \int \sum\limits_{  {\left| v \right| = \left| u \right| = N}\atop {P^{\raise1pt \hbox{$\scriptscriptstyle N+$}}(T_v x)=\alpha, Q^{\raise1pt \hbox{$\scriptscriptstyle N+$}}(S_{u} y)=\beta} } {\int 
A(v,u)\,dm_{x,y} (A)}\,d\lambda (x,y) \cr 
 &= \int {\int {\sum\limits_{  {\left| v \right| = \left| u \right| = N}\atop {P^{\raise1pt \hbox{$\scriptscriptstyle N+$}}(T_v x)=\alpha, Q^{\raise1pt \hbox{$\scriptscriptstyle N+$}}(S_{u} y)=\beta} } 
{A(v,u)} }\, dm_{x,y} (A)d\lambda (x,y)} \cr 
&= \int{\int {\lambda _{(\tau _x^P ,\tau _y^Q ,A)} (\alpha \times \beta )\,dm_{x,y} 
(A)d\lambda  (x,y)}}. \cr }
$$
{\closeproof}
\medbreak
We will now prove the copying lemma for one-sided joinings. Note that 
while it appears to be more general than the copying lemma in [H,R] 
(Proposition 5.4) in that we also copy distribution of tree names, the proof 
is essentially the same.  Following the convention in Chapter 2, we will work with functions taking values in [0,1).
\proclaim Proposition 3.2.4 (Copying Lemma). Suppose \X and \Y are ergodic $\bf p$-endomorphisms.  Let $\lambda \in$  \joinXY. If $g:X \to [\hbox{0},\hbox{1})$ and $h:Y\to [\hbox{0},\hbox{1})$, then for all $\varepsilon $ and $N$, there exists a function $\bar g$ on $Y$  such that
$$
\left| dist(g_N^{N\nabla})-dist({\bar g_N}^{N\nabla})\right|<\varepsilon
$$
\noindent and
$$ \left| dist(g_N^{N+}\otimes h_N^{N+})
-dist({\bar g_N}^{N+}\vee h_N^{N+}) \right|<\varepsilon .$$
\par 

\noindent {\twelvebx Proof: } Let $\tilde {g} = g_N^{N \nabla } \vee g_N$ and $\tilde {h} = h_N^{N{\nabla }} \vee h_N$ and let $\bar D$ denote $(D_N )^{N \nabla} \times D_N$. (Recall that $D_N$ is the set of midpoints of the dyadic intervals $[t/\hbox{2}^N,(t+\hbox{1})/\hbox{2}^N).$)  Choose $M>N$
such that $\hbox{2}N/M<\varepsilon / \hbox{2}$. Construct a $\varepsilon / \hbox{2}$-tree Rokhlin 
tower ${\bf M}$ of height $M+\hbox{1}$ in ${\bf Y}$. By Proposition 3.2.3, we have measures $m_{x,y}$ on ${\cal A}_M$ such that for each $N\le j \le M$ and for each pair of elements $(\alpha ,\beta)$ in $\bar D ^{N}$,
$$
\lambda _{\tilde {g} \otimes \tilde {h}} ((\alpha \times \beta)^{-j} ) = \int \int 
{\lambda _{(\tau _x^{\tilde g} ,\tau _y^{\tilde h} ,A)} ((\alpha \times \beta)^{-j} 
)\,dm_{x,y} (A)d\lambda  (x,y)}. \eqno(\hbox{1})
$$

Consider the partitions $\tilde g ^{M \tau} :X \to \bar D ^{M \tau}$ and $\tilde h ^{M \tau} :Y \to \bar D ^{M\tau}$, we can assign a measure $\hat \lambda$ to each atom of the partition $\tilde g ^{M \tau} \times \tilde h ^{M \tau} \times {\cal A}_M$ of $X \times Y \times {\cal A}_M$ by setting
$$\hat {\lambda }(p \times q \times A) = \int _{p\times q} {m_{x,y}(A) \,d\lambda 
(x,y)}. \eqno(\hbox{2})
$$
By the Strong 
Tree Rokhlin Lemma, we may assume that the base $C$ of ${\bf M}$ is chosen such that $dist(\tilde {h}^{M\tau }|C)=dist(\tilde {h}^{M\tau })$. For each atom $\gamma  \in \tilde 
{h}^{M{\tau }}$, we define a partition 
$$P_{\gamma } : \gamma  \cap 
C \to \bar D ^{M \tau} \times \bar D ^{M \tau} \times {\cal A}_M$$
\noindent such that 
$$
dist(P_{\gamma })= {\hat \lambda}({{\tilde {g}}^{M \tau }}\times {\tilde 
{h}}^{M{\tau }}\times {\cal A}_{M}|\gamma ).
$$
\noindent The partitions $P_{\gamma }$ over all $\gamma $ collectively define a 
partition $P$ of $C$ such that
$$
dist(P)=dist(\tilde {g}^{M{\tau }}\times \tilde {h}^{M{\tau }}\times 
{{\cal A}_{M}}).
$$
\noindent This gives a bijective correspondence 
$\rho $ of the atoms of $P$ and those of $\tilde {g}^{M{\tau }}\times \tilde 
{h}^{M{\tau }}\times {\cal A}_{M}$ such that if $\alpha  \in P$ and 
$\rho (\alpha )=(\beta ,\gamma ,A)$, then $\alpha  \subseteq 
\gamma $, and $\nu (\alpha | C)=\hat {\lambda }(\rho (\alpha ))$. 

We will now construct the required function $\bar g$ on the tower. Fix an atom $\alpha \in P$ such that $\rho (\alpha )=(\beta ,\gamma ,A)$. Choose any point $x \in \beta $. We define the function $\bar g$ on $ 
\cup _{i = 1}^M \,S^{ - i}\alpha $ by setting ${\bar g}(z) = g(T_{A^{ - 1}v} x)$ 
for $z \in S_{v}\alpha $. By repeating this procedure for each atom 
$\alpha \in P$, we can extend $\bar g$ to a function on ${\cup \bf M} \backslash B$. Extend the 
function $\bar g$ to the rest of $Y$ in any way we like. 

Set $\tilde {\bar g} = {\bar g}_N^{N{\nabla }} \vee {\bar g}_N$. For an integer $j$ such that $N \le j
\le M-N$, we claim that
$$
dist(\tilde {g}^{N + }\otimes \tilde {h}^{N + }) = 
dist(\tilde {\bar g}^{N + } \vee \tilde {h}^{N + }|S^{ - 
j}C_\phi ) \eqno(\hbox{3})
$$
If $\alpha  \in P$ and $\rho (\alpha )=(\beta ,\gamma ,A)$, then by 
our construction, notice that
$$
({\bar g}_{N}(S_{v}y),h_{N}(S_{v}y))=(g_N (T_{A^{ - 1}v} x),h_N (S_v 
y))
$$
\noindent for any $y \in \alpha $ and $x \in \beta $ and $v \in {\cal T}'_{M}$. 
We thus have
$$
(\tilde{\bar g}(S_v y),\tilde {h}(S_v y)) = 
(\tilde{g}(T_{A^{ - 1}v} x),\tilde {h}(S_v y)) \eqno(\hbox{4})
$$
\noindent for $v \in {\cal T}'_{M - N}$. 

Let $\zeta _{1} \vee \zeta _{2}=\{ y \in Y \mid \tilde {\bar g}^{N + }(y)=\zeta_1 , \; \tilde{h}^{N +}(y)= 
\zeta _2 \}$. Then (4) implies that for any $x \in \beta $ and $y \in 
\gamma $, 
$$
\nu (\zeta _{1} \vee \zeta _{2 }|S^{ - j}\alpha ) = \lambda _{(\tau _x^{\tilde 
g},\tau _y^{\tilde {h}} ,A)} ((\zeta _1 \times \zeta _2 )^{ - j}). \eqno(\hbox{5})
$$
Note that (1) and (2) implies that for $x(\beta) \in \beta$ and $y(\gamma) \in \gamma$,
$$
\lambda _{\tilde {g} \otimes \tilde {h}} ((\zeta _1 \times \zeta _2 )^{ - 
j}) = \sum\limits_{(\beta ,\gamma ,A) \in {\tilde g}^{M \tau }\times {\tilde 
h}^{M \tau }\times {\cal A}_M } {\lambda _{(\tau _{x(\beta)}^{\tilde g},\tau 
_{y(\gamma)}^{\tilde h} ,A)} ((\zeta _{1}\times \zeta _{2})^{ - j}) {\hat \lambda 
}(\beta ,\gamma ,A)}. \eqno(\hbox{6})
$$ 
\noindent Hence, since $\nu (S^{ - j}\alpha | S^{ - j}C) = \nu (\alpha 
| C ) = \hat {\lambda }(\beta ,\gamma ,A)$, (5) and (6) imply that 
$$
\eqalign{
\lambda _{{\tilde g} \otimes {\tilde h}} ((\zeta _1 \times \zeta _2 )^{ - 
j}) &= \sum\limits_{\alpha \in P} {\nu (\zeta _1 \vee \zeta _2 |S^{ - 
j}\alpha ) \nu (S^{ - j}\alpha | S^{ - j}C)} \cr 
&= \nu ( \zeta _1 \vee \zeta _2 |S^{ - j}C). \cr} 
$$
It follows that for $N \le j \le M-N$, 
$$
dist({\tilde {g}^{N + }}\otimes {\tilde {h}^{N + }}) = 
dist({\tilde {\bar g}^{N + }}\vee {\tilde {h}^{N + }}|S^{ - j}C),
$$
\noindent which is (3).

As $\nu ( \cup _{i = N}^{M - N} \,S^{ - j}C) > \hbox{1}-\varepsilon 
$, we then have
$$
\left| dist(\tilde {g}^{N + }\otimes {\tilde {h}^{N + 
}}) - dist(\tilde {\bar g}^{N + }\vee {\tilde {h}^{N + }}) 
\right| < \varepsilon, 
$$
\noindent from which the conclusion easily follows, since $\tilde{g}$ refines the partitions $g_N$ and $g_N^{N \nabla}$ (and similarly for $\tilde {\bar g}$ and $\tilde h$). {\closeproof} 
\medjump
\noindent {\section 3.3. The {\= t Distance} }
\medjump
Suppose \X and \Y are ${\bf p}$-endomorphisms. Following [H,R], for a fixed $m \in {\bf N}$ and functions $g:X \to \hbox{[0,1)}$ and $h:Y \to \hbox{[0,1)}$ respectively, define the $\bar t _m$ distance  between 
the pair of processes $({\bf X},g)$ and $({\bf Y},h)$ by
$$
\bar {t}_m (({\bf X},g),({\bf Y},h))  = \mathop {\hbox{inf}}\limits_{\lambda \in C^+ ({\bf X},{\bf Y})}  
\int {\bar {t}_m (\tau _x^g ,\tau _y^h )\,d\lambda }.
$$
\noindent Let
$$
\bar {t}(({\bf X},g),({\bf Y},h))  = \hbox{liminf}\;\bar {t}_m (({\bf X},g),({\bf Y},h)) .
$$
Note that this definition of $\bar {t}_m (({\bf X},g),({\bf Y},h))$ differs from the 
definition in \S 2.4 in which we consider only the product joining $\lambda 
$. Nonetheless, note that Proposition 2.4.6 (which will be used in \S 3.4) 
will also hold with the present definition. In fact, for a tvwB ${\bf X}$ and $\varepsilon > \hbox{0}$, if $dist(g_N^{N\nabla })$ is sufficiently 
close to $dist(h_N^{N\nabla })$ for some $N$, the proof of Proposition 
2.4.6 shows that $\int {\bar {t}_m (\tau _x^g ,\tau _y^h )\,d\lambda } 
<\varepsilon $ for all large $m$ for {\twelvebx any} coupling $\lambda $ (not just one-sided).

Our goal in this section is to prove the following proposition which will be combined with the copying lemma in the previous section to give a joinings proof of Theorem 2.4.1.  
\proclaim Proposition 3.3.1. Suppose $\bar {t}(({\bf X},g),({\bf Y},h)) < \varepsilon$, then there exists $\lambda  
\in$ \joinXY $ $ such that $\int {\left| {g(x) - h(y)} \right|\,d\lambda < 
\varepsilon }$.\par
The proof of Proposition 3.3.1 will depend on following proposition whose proof will be postponed until the end of this section.  
\proclaim Proposition 3.3.2. For each $m \in {\bf N}$, there exists a $\lambda \in C^+({\bf X}, {\bf Y})$ such that
$$ {\hbox{1} \over m} \int {\sum_{i=0}^{m-1} \left| g(T^i x) - h(S^i y) \right| \; d \lambda (x,y)} \le \bar {t}_m 
(({\bf X},g),({\bf Y},h)) + \hbox{1} / \hbox{2}^{m - 2}.$$ \par

\noindent {\twelvebx Proof (of Proposition 3.3.1):} By Proposition 3.3.2, for each $m \in {\bf N}$, we have some $\lambda_m \in C^+ ({\bf X},{\bf Y})$ such that
$$
{ \hbox{1} \over m} \int {\sum_{i=0}^{m-1} \left| g(T^i x) - h(S^i y) \right| \; d \lambda_m (x,y)} \le  \bar {t}_m 
(({\bf X},g),({\bf Y},h)) + \hbox{1} / \hbox{2}^{m - 2}.
$$
Choose an increasing sequence of integers $n_{1},n_{2},\ldots$ such that 
$$
\bar {t}_{n_m } (({\bf X},g),({\bf Y},h)) \to \bar {t}(({\bf X},g),({\bf Y},h)),
$$
\noindent
so that by passing to a further subsequence if necessary, we have some $L<\varepsilon$ such that
$$
{\hbox{1} \over {n_m}} \int {\sum_{i=0}^{n_m -1} \left| g(T^i x) - h(S^i y) \right| \; d \lambda_{n_m} (x,y)} \to L. \eqno(\hbox{1})
$$
For $i \ge \hbox{0}$ and any measure $\lambda$ on $X \times Y$, let $(T \times 
S)^{i} \lambda $ denote the measure defined by $ (T \times 
S)^{i} \lambda (D) = \lambda ((T \times S)^{-i} D)$ for $D \subseteq X \times Y$. Define 
$$\bar \lambda _m = {\hbox{1} \over {n_m} }\sum\limits_{i = 0}^{n_m - 1} {(T\times 
S)^{i} \lambda_{n_m } }.$$ 
Since the set of one-sided couplings is convex, 
$\bar \lambda _{m} \in C^{ + }({\bf X},{\bf Y})$. For any measure $\lambda$ on $X \times Y$ and $i \ge \hbox{0}$, let
$$
\lambda (g^{i} \Delta h^{i}) = \int {\left| {g(T^ix) - h(S^iy)} \right|\,d\lambda }.
$$
\noindent Then,
$$
\eqalign{
\bar \lambda _m (g^{0}\Delta h^{0}) & = {\hbox{1} \over {n_m }}\sum\limits_{i = 0}^{n_m - 1} 
{(T\times S)^{i} \lambda_{n_m } (g^{0}\Delta h^{0})} \cr 
 &= {\hbox{1} \over {n_m }}\sum\limits_{i = 0}^{n_m - 1} {\lambda_{n_m } 
(g^{i}\Delta h^{i})}. 
\cr }\eqno(\hbox{2}) 
$$
Let $\lambda ^{\ast}$ be any $w ^{\ast}$-limit point of the $\bar \lambda _{m}$'s. Since $C^{ 
+ }({\bf X},{\bf Y})$ is $w ^{\ast}$-closed in $C({\bf X},{\bf Y})$, $\lambda ^{\ast} \in 
C^{ + }({\bf X},{\bf Y})$. It is clear that $\lambda ^{\ast}$ is stationary so 
$\lambda ^{\ast} \in J^{ + }({\bf X},{\bf Y})$. By (1) and (2), given $\delta >\hbox{0}$, we 
have for all sufficiently large $m$,
$$
\lambda ^{\ast} (g^{0}\Delta h^{0})\mathop \sim \limits^\delta \bar \lambda _m 
(g^{0}\Delta h^{0})\mathop \sim \limits^\delta L.
$$
\noindent Hence, $\lambda ^{\ast} (g^{0}\Delta h^{0})=L<\varepsilon $, since 
$\delta $ is arbitrary.  {\closeproof}
\medbreak 
We now turn our attention to the proof of Proposition 3.3.2.  The proof will require the construction of 
a one-sided coupling of \X and \Y, which we now describe.  
Given functions $g:X \to [\hbox{0},\hbox{1})$ and $h:Y\to [\hbox{0},\hbox{1})$, and $m \in {\bf N}$, we have the joined 
partition $g_m^{m\tau } \otimes h_m^{m\tau }$ of $X\times Y$ defined by 
$g_m^{m\tau } \otimes h_m^{m\tau }(x,y)=(g_m^{m\tau }(x),h_m^{m\tau }(y))$. 
For an atom $\beta  \in g_m^{m\tau } \otimes h_m^{m\tau }$, let 
$A_{\beta }$ denote a tree automorphism in ${\cal A}_{m}$ which realizes 
$\bar {t}_m (\tau _x^{g_m } ,\tau _y^{h_m } )$ for any $(x,y) \in \beta $. 
Set $A_{x,y}=A_{\beta }$ for any $(x,y) \in \beta $. 

For each $(x,y) \in X\times Y$, we define a measure supported on ${(T\times S)^{-m}}(x,y)$, which we will denote as ${A}_{x,y}^m$, in the following way. For any measurable set $D \subseteq X\times Y$, let ${A}_{x,y}^m (D)$ be the 
total weights of the nodes $v$ of length $m$ such that $(T_v x,S_{A_{x,y} v} 
y) \in D$. Given $\lambda \in C^+({\bf X},{\bf Y})$, define the measure $\bar {\lambda }$ on $(X\times Y,{\cal 
B}\times {\cal C})$ by setting
$$
\bar {\lambda }(D) = \int {{A}_{x,y}^m (D)d\lambda (x,y)} 
$$
\noindent for $D \subseteq X\times Y$. It is straightforward to check that $\bar 
{\lambda }$ defines a probability measure on $X\times Y$. We now show:

\proclaim  Lemma 3.3.3. $\bar {\lambda } \in C^ + ({\bf X},{\bf Y}) $. \par
\noindent {\twelvebx Proof:} First, we need to see that $\bar {\lambda }$ is a coupling. Let 
$\{ \mu _{x}\}$ be a disintegration of $\mu $ induced by the factor 
map $T^{m}:{\bf X} \to {\bf X}$. If $B \subseteq X$, then ${A}_{x,y}^m (B) = \mu _x (B)$. Thus,
$$
\bar {\lambda }(B) = \int {{A}_{x,y}^m (B)\,d\lambda } = \int {\mu 
_x (B)\,d\mu } = \mu (B)
$$
\noindent
and so $\bar {\lambda }$ projects to $\mu$ on $X$. By symmetry, $\bar {\lambda }$ 
projects to $\nu$ on $Y$. Hence $\bar {\lambda } \in C({\bf X},{\bf Y})$.

We now need to check that $\bar {\lambda }$ is one-sided. We first check 
condition i) of Definition 3.1.1.  For convenience, let us represent ${\bf X}$ and ${\bf Y}$ as one-sided shift spaces on $\hbox{[0,1]}^{{\bf N}^\ast 
}$ (i.e. $X=Y=\hbox{[0,1]}^{{\bf N}^\ast }$, and $T=S$ 
is the shift map). Given integers $\hbox{0} \le i \le j$, and a point $z \in \hbox{[0,1]}^{{\bf N}^\ast }$, we let 
$$z[i,j]=\{y \in \hbox{[0,1]}^{{\bf N}^\ast } \mid y_t=z_t, i \le t 
\le j\}$$ 
\noindent and 
$$z[i,\infty )=\{y \in \hbox{[0,1]}^{{\bf N}^\ast 
} \mid y_t=z_t, t \ge i\}.$$ 

By symmetry, it clearly suffices to prove that for each integer $t>\hbox{0}$ and 
for $\bar {\lambda }$-$a.a$ $(x',y')$ in $X\times Y$,
$$
E_{\bar \lambda }(x'[\hbox{0},t - \hbox{1}]| x'[t,\infty )\times 
y'[t,\infty )) = E_{\mu} (x'[\hbox{0},t - \hbox{1}]| x' [t,\infty )).\eqno(\hbox{1})
$$
We prove (1) by considering three cases: $t=m$, $t<m$ and $t>m$:
 
\noindent \underbar {\bf Case I: $t=m$} 

Using the definition of $\bar {\lambda }$, if $x=T^{m}x'$ and 
$y=S^{m}y'$, we have
$$ \eqalign{
 E_{\bar \lambda }(x'[\hbox{0},m - \hbox{1}]| x' [m,\infty )\times 
y'[m,\infty )) &=E_{\bar \lambda }(x'[\hbox{0},m - \hbox{1}]| T^{-m}x \times 
S^{-m}y) \cr &= {A}_{x,y}^m (x'[\hbox{0},m - \hbox{1}]) \cr 
 &= E_{\mu} (x'[\hbox{0},m - \hbox{1}]| x' [m,\infty )). \cr }
$$

\noindent \underbar {\bf Case II: $t<m$}

Observe that for $\bar {\lambda}$-a.a. $(x',y')$, if $x=T^m x'$ and $y=S^m y'$, then 
${A}_{x,y}^m (x',y')>\hbox{0}$.  For each such $(x',y')$, we have
$$
{A}_{x,y}^m (x'[\hbox{0},m - \hbox{1}]\times y'[\hbox{0},m - \hbox{1}]) = E_{\mu} (x'[\hbox{0},m - \hbox{1}]| x'[m,\infty )).
$$
\noindent We may then use Case I) to show that 
$$
\eqalign{
E_{ \bar \lambda} (x'[\hbox{0},t - \hbox{1}]| x' [t,\infty )\times 
y'[t,\infty )) &= {{E_{\bar \lambda }(x'[\hbox{0},m - \hbox{1}]\times y'[t,m - 
\hbox{1}]| x' [m,\infty )\times y'[m,\infty ))}\over {E_{\bar \lambda 
}(x'[t,m - \hbox{1}]\times y'[t,m - \hbox{1}]| x' [m,\infty )\times 
y'[m,\infty ))}} \cr 
&= {{{A}_{x,y}^m (x'[\hbox{0},m - \hbox{1}]\times y'[t,m - \hbox{1}])}\over {{A}_{x,y}^m  (x'[t,m - \hbox{1}]\times y'[t,m - \hbox{1}])}} \cr 
&= {{E_{\mu} (x'[\hbox{0},m - \hbox{1}]| x'[m,\infty))}\over {E_{\mu} 
(x'[t,m - \hbox{1}]| x'[m,\infty ))}} \cr 
&= E_{\mu} (x'[\hbox{0},t - \hbox{1}]| x'[t,\infty )). \cr} \
$$
\noindent \underbar {\bf Case III: $t>m$}

This follows from a direct calculation using Case I) and the fact that $\lambda $ is one-sided.
\medbreak
We now prove condition ii) of Definition 3.1.1. Recall the 
Bernoulli shift $B^ + (\bar {\bf p})=(I(\bar {\bf p})=\{ \hbox{1},...,r\}^{{\bf N}^\ast 
},{\bf m},\sigma )$ as defined in \S 3.1. Consider the cylinder set 
$$
C=\{z \in I(\bar {\bf p})\mid z_0=r_{0},\ldots ,z_{m-1}=r_{m-1}\},
$$

\noindent where $r_{j} \in \{ \hbox{1},\ldots ,r\}$. By the definition of the maps $\psi 
_{\bf X}$ and $\psi _{\bf Y}$, we have a set $V_{C}$ of nodes of length $m$ whose 
total weights is ${\bf m}(C)$ such that 

$$
\psi _{\bf X}^{ - 1} (C) = \mathop \cup \limits_{v \in V_C } T_v X \quad \hbox{and} \quad \psi _{\bf Y}^{ - 1} (C) = \mathop \cup \limits_{v \in V_C } S_v Y.
$$
\noindent Hence,
$$
{A}_{x,y}^m (\psi _{\bf X}^{ - 1} (C)\times \psi _{\bf Y}^{ - 1} (C)) = 
\sum\limits_{v,u \in V_C } {A_{x,y} (v,u)} = {\bf m}(C)
$$
\noindent Condition ii) now follows immediately from the definition of $\bar {\lambda }$ and the one-sidedness of $\lambda $. {\closeproof}
\medbreak
We are now ready to prove Proposition 3.3.2.

\noindent {\twelvebx Proof (of Proposition 3.3.2): } Consider the partition $g_m^{m\tau } \otimes h_m^{m\tau }$ on 
$X\times Y$. Let $\lambda \in C^+ ({\bf X},{\bf Y})$, construct $\bar \lambda$ as in Lemma 3.3.3.  For each atom $\beta  \in g_m^{m\tau } \otimes h_m^{m\tau 
}$ and a pair of nodes $(v,u)$ of length $m$, note that
$$
 \bar {\lambda }((T_v \times S_u )\beta ) = 
A_{\beta} (v,u) \lambda (\beta ). \eqno(\hbox{1})
$$
Choose a point $(x_{\beta },y_{\beta }) \in \beta$ for each atom $\beta$.  Note that by (1) and our choice of $A_\beta$,  
$$
\eqalign{{}&\qquad {\hbox {1} \over m}\sum\limits_{i = 0}^{m - 1}\int _{(T \times S)^{-m}\beta} \left| g_m(T^i x) - h_m(S^i y) \right| d \bar \lambda  \cr
&= \sum\limits_{\left| v \right|=\left| u \right|=m}\int _{T_v \times S_u(\beta)} {\hbox {1} \over m}\sum\limits_{i = 0}^{m - 1}\left| g_m(T^i x) - h_m(S^i y) \right| d \bar \lambda \cr
&= \sum\limits_{\left| v \right|=\left| u \right|=m} {\hbox {1} \over m}\sum\limits_{i = 0}^{m - 1}\left| g_m(T^i T_v x_\beta) - h_m(S^i S_u y_\beta) \right| A_{\beta}(v,u)\lambda (\beta) \cr
&= \bar t_m (\tau_{x_\beta}^{g_m},\tau_{y_\beta}^{h_m}) \lambda (\beta). \cr}
$$
Thus,
$$
\eqalign{
 {\hbox{1} \over m}\sum\limits_{i = 0}^{m - 1} \int {\left| {g_m (T^{i} x) - h_m 
(S^{i} y)} \right|\,d\bar \lambda }  &= {\hbox {1} \over m}\sum\limits_{i = 0}^{m - 1}\sum \limits_{\beta \in g_m^{m\tau } \otimes h_m^{m\tau }} \int _{(T \times S)^{-m}\beta} \left| g(T^i x) - h(S^i y) \right| d \bar \lambda \cr
&= \sum \limits_{\beta \in g_m^{m\tau } \otimes h_m^{m\tau }} \bar t_m (\tau_{x_\beta}^{g_m},\tau_{y_\beta}^{h_m}) \lambda (\beta) \cr 
&= \int {\bar {t}_m (\tau _x^{g_m }  ,\tau _y^{h_m } )}\,d\lambda . \cr} 
$$
Since $\left| g-g_{m}\right| \le \hbox{1}/ \hbox{2}^{m+1}$, we have
$$
{1 \over m}\sum\limits_{i = 0}^{m - 1} {\int {\left| {g (T^{i} x) - h 
(S^{i} y)} \right|\,d\bar {\lambda }} } \le
\int {\bar {t}_m (\tau _x^{g }  ,\tau _y^{h } )}\,d\lambda + \hbox{1}/ \hbox{2}^{m-1}.$$ 
As this argument holds true for all $\lambda \in C^ + ({\bf X},{\bf Y}),$ the result follows. {\closeproof}

\medjump
\noindent{\section 3.4. A Joinings Proof of Theorem 2.4.1}
\medjump

In this section, we prove Theorem 2.4.1 using the machinery of one-sided 
joinings developed in \S 3.1 to \S 3.3. In fact, we will prove a stronger result.  We say that $\lambda \in$ \joinXY $ $ is an {\twelvebx isomorphic joining} if  ${\cal B}\mathop = 
\limits^\lambda {\cal C}$.  Note that each isomorphic joining is the graphical joining of an isomorphism from ${\bf X }$ to ${\bf Y}$.  We now prove
\proclaim Theorem 3.4.1. Let $\bf  p$  be a probability vector. If \X $\in 
End({\bf p})$ is tvwB, the set of isomorphic joinings is a dense $w^{\ast}$-$G_\delta$ of $J^+({\bf X},B^+({\bf p}))$. \par
\noindent {\twelvebx Proof:} Let ${\bf Y}=B^+ ({\bf p})$.  Choose a pair of generating functions $f:X\to [\hbox{0},\hbox{1})$ and $g:Y\to [\hbox{0},\hbox{1})$. Fix $N \in {\bf N}$ and $\delta > \hbox{0}$. 

Let $\lambda  \in$ \joinXY.  For each $N' \in {\bf N}$ and $\delta' >\hbox{0}$,  by the copying lemma, we have a function $\bar f$ 
on $Y$ such that 
$$
\left| dist(f_{N'}^{N' \nabla}) - 
dist({\bar f}_{N'}^{N' \nabla}) \right| < \delta'  \eqno(\hbox{1})
$$
\noindent and 
$$
\left| {dist_{\lambda}(f_{N'}^{N' + }\otimes{g}_{N'}^{N' + }) - 
dist({\bar f}_{N'}^{N' + }\vee {g}_{N'}^{N' + })} \right| < \delta' . 
 \eqno(\hbox{2})
$$
For $\varepsilon >\hbox{0}$, if $\hbox{1}/{N'}$ and $\delta '$ are sufficiently small then by (1) and the proof of Proposition 
2.4.6, 
$$
\bar {t}(({\bf X},f),({\bf Y},{\bar f})) < \varepsilon .
$$
By Proposition 3.3.1, we have some $\lambda ' \in$ \joinXY $ $ such that 
$$
\int {\left| {f(x) - {\bar f}(y)} \right|\,d\lambda ' (x,y) < \varepsilon } . \eqno(\hbox{3})
$$
Then provided that $\varepsilon$, $\hbox{1}/N'$, $\delta'$ 
are small enough, we may use (2) to conclude that
$$
\left| {dist_{\lambda}(f_N^{N + }\otimes {g}_N^{N + }) - 
dist_{\lambda'}(f_N^{N + } \otimes {g}_N^{N + })} \right| < \delta . \eqno(\hbox{4})
$$
As in [H,R], for $\lambda \in$ \joinXY, we say that $f\mathop \subset 
\limits_\lambda ^\varepsilon Y$ if there exists some function $h$ 
on $Y$ such that $\int {\left| {f(x) - h(y)} \right|\,d\lambda < \varepsilon } 
$. Let $O_{\varepsilon }=\{ \lambda \in J^+ ({\bf X},{\bf Y}) \mid {f} \mathop \subset 
\limits_\lambda ^{\varepsilon} Y\}$. Then $O_{\varepsilon }$ is $w^{\ast}$-open in 
\joinXY, and (3) and (4) imply that $O_{\varepsilon }$ is $w^{\ast}$-dense in \joinXY. Define $f\mathop \subset \limits_\lambda ^0 Y$ if $f\mathop 
\subset \limits_\lambda ^\varepsilon Y$ for all $\varepsilon >\hbox{0}$.  Note that $f\mathop \subset \limits_\lambda ^\varepsilon Y$ implies that ${\cal B} \mathop \subset \limits^\lambda {\cal C}$.  

Now, by Baire's Theorem, $\mathop  \cap \limits_{n = 1}^\infty  O_{1/n} $
 is a dense $w^{\ast}$-$G_{\delta }$ in $J^{+} ({\bf X},{\bf Y})$. By symmetry, if  $O'_{\varepsilon }=\{\lambda \in J^+ ({\bf X},{\bf Y}), g\mathop \subset 
\limits_\lambda ^\varepsilon X\}$, we see that the set $\mathop  \cap \limits_{n = 1}^\infty  O'_{1/n} $
 is also dense in $J^+ ({\bf X},{\bf Y})$. Hence, the intersection $\mathop  \cap \limits_{n = 1}^\infty  O_{1/n}  \cap O'_{1/n} $
 is a dense 
$w^{\ast}$-$G_{\delta }$ in $J^+ ({\bf X},{\bf Y})$ so that it is a residual set. This concludes the proof as  
any $\lambda $ in the intersection satisfies $f\mathop \subset 
\limits_\lambda ^0 Y$ and $g\mathop \subset 
\limits_\lambda ^0 X$ so that $\lambda$ is an isomorphic joining.   {\closeproof}
\medbreak

This concludes the proof of the isomorphism theorem, using the machinery of one-sided joinings.  Note that by Proposition 3.1.3, the set of isomorphic joinings is non-empty.  Thus, $\bf X$ and $B^+({\bf p})$ are isomorphic.  

One corollary of the joinings proof is that unless the components of ${\bf p}$ 
are pairwise distinct, there are uncountably many automorphisms of the 
one-sided Bernoulli shift $B^{ + }({\bf p})$ (and hence in general between 
any two tvwB ${\bf p}$-endomorphisms). To prove this, we first need the 
lemma below, which is just a special case of Proposition 3.1.3. In the 
following, we suppose $B^{ + }({\bf p})$ is represented as the one-sided 
shift space $(\{ \hbox{1},\ldots ,s\}^{{\bf N}^\ast },\mu ,T)$ for the Bernoulli 
measure $\mu =\{p_{1},\ldots p_{s}\}^{{\bf N}^\ast }$ and the shift map $T$ 
(so that $j$ has weight $p_{j}$ for $\hbox{1} \le j \le s $).  

\proclaim Lemma 3.4.2. Suppose $h: \{ \hbox{1},\ldots ,s\}^{{\bf N}^\ast } \to \{ \hbox{1},\ldots 
,s\}$ is tree-adapted and for $a.a.x$ in $\{ \hbox{1},\ldots ,s \}^{{\bf N}^\ast }$, $x_0$ and $h(x)$ have the same weight. Then the $h$-name map $\theta :B^{ + }({\bf p}) \to B^{ + }({\bf p})$ defined by $\theta (x)=(h(x),h(Tx),\ldots )$ is a tree-adapted factor map. Moreover, the graphical self-joining $\lambda _{\theta }$  of $B^{ + 
}({\bf p})$ derived from $\theta $ is one-sided.  {\closeproof}\par

\proclaim Proposition 3.4.3. For a probability vector $\bf p$ with at least two identical components, there are uncountably many automorphisms of $B^{ + }({\bf p})$. \par 
\noindent {\twelvebx Proof:} Without loss of generality, suppose $p_{1}=p_{2}$. It 
suffices to show that $J^+ (B^+ ({\bf p}),B^+ ({\bf p}))$ has no isolated points. Indeed, from the proof of the isomorphism theorem, the isomorphic joinings 
is a dense $w^{\ast}$-$G_{\delta }$ in $J^+ (B^ + ({\bf p}),B^ + ({\bf p}))$ . However, by the 
Baire Category Theorem, a dense $w^{\ast}$-$G_{\delta }$ in a complete metric space 
with no isolated points is necessarily uncountable. 

To this end, notice that there are at least two elements in $J^+ (B^+ ({\bf p}),B^+ ({\bf p}))$.  Indeed, the graphical joining derived from the identity automorphism certainly is one.  Another one is the graphical joining derived from the $h$-name map for the function $h : { \{ \hbox{1},\ldots ,s\}^{{\bf N}^\ast}} \to \{ \hbox{1},\ldots ,s \}$ defined by
$$h(x)=\cases {\hbox{2}, &if $x_0=\hbox{1}$; \cr
 \hbox{1}, &if $x_0=\hbox{2}$; \cr
x_0, &otherwise. \cr}$$
\noindent (This is one-sided by Lemma 3.4.2).

  Let $\lambda \in J^+ (B^ + ({\bf p}),B^ + ({\bf p}))$ and  
$\varepsilon > \hbox{0}$. It suffices to construct ${\lambda }' \in J^+ (B^+ ({\bf p}), B^+ ({\bf p}))$ such that $\lambda '  \ne \lambda $ and $\lambda 
' \mathop \sim \limits^\varepsilon \lambda $ in $w^{\ast}$. To see this, since there are at least two elements in $J^+ (B^+ ({\bf p}),B^+ ({\bf p}))$, just choose some one-sided joining $\lambda _1$ distinct from $\lambda$.  By Proposition 3.1.2, convex combinations of $\lambda$ and $\lambda _1$ remain one-sided.  Let $\lambda_{\delta}=(\hbox{1}-\delta) \lambda + \delta \lambda_1.$  If ${\delta}$ is small enough, then $\lambda_{\delta}$ will be sufficiently close to $\lambda$ in $w^{\ast}$, and we are done.  {\closeproof}
\vfill \eject

$ $
\vskip 5.90cm
\noindent {\chapter Chapter 4:} 
\medjump
\noindent {\chapter Examples of TvwB p-Endomorphisms}
\bigjump
\twelverm

In this chapter, we will present two classes of examples of tvwB 
${\bf p}$-endomorphisms. It follows from the isomorphism theorem, Theorem 
2.4.1, that all of these are isomorphic to $B^{ + }({\bf p})$.
\medjump
\noindent{\section 4.1. One-Sided Markov Shifts }
\medjump
Besides the 
one-sided Bernoulli shift $B^{ + }({\bf p})$, the simplest examples of ${\bf 
p}$-endomorphisms can be found among the {\twelvebx one-sided Markov shifts}, which we will define below. 

It will be convenient for us to define a one-sided Markov shifts over left-infinite shift spaces of the form $C^{ - {\bf N}}$ for a finite set $C$. We say that 
a square matrix $A$ is {\twelvebx stochastic} if each entry is nonnegative and the sum of the 
entries in each row equals one. Let $\left| A \right|$ denote the number of rows (or 
columns) of $A$. We may obviously index the rows and columns of $A$ by the 
integers from 1 to $\left| A \right|$ and denote the entries in $A$ by $A_{ij}$, where 1$ \le 
i,j \le \left| A \right|$. A stochastic matrix is {\twelvebx irreducible} if for each pair $(i,j)$, there exists 
$k \in {\bf N}$ such that $(A^{k})_{ij}$ is nonzero. It is well known that 
for an irreducible stochastic matrix, there exists a unique row probability 
vector ${\bf q}$ with all components positive such that ${\bf q}A={\bf q}$. We 
say that ${\bf q}$ is a {\twelvebx left fixed probability vector of $A$}. Let $q_{j}$ denote the $j$-th component of ${\bf q}$. 
Using ${\bf q}$ and $A$, we may define a measure $\mu $ on the cylinder sets of 
$\{ \hbox{1},\ldots ,\left| A \right| \} ^{- {\bf N}}$ by
\vfill \eject
$$
\mu (x \in {\{ \hbox{1},\ldots , \left| A \right| \}}^{ - {\bf N}} \mid x_{-n}=a_{ - n},\ldots 
,x_{-1}=a_{ - 1})=q_{a_ {- n}} A_{a_{ - n} a_{ - (n - 1)}}\cdots A_{a_{ - 2} a_{ - 1}}. 
$$
\noindent It is easily seen that $\mu $ extends to a measure on the Borel 
sigma-algebra ${\cal B}$ of $\{ \hbox{1},\ldots ,\left| A \right| \}^{- {\bf N}}$. If $T$ is the 
shift map on $\{ \hbox{1},\ldots ,\left| A \right| \}^{-{\bf N}}$, then note that $\mu =\mu 
T^{ - 1}$ (i.e. $\mu $ is shift invariant). We define the {\twelvebx one-sided Markov shift over $A$}, denoted $X_A^ 
- $, to be the m.p.s. $(\{ \hbox{1},\ldots ,\left| A \right| \}^{- {\bf N}},{\cal B},\mu,T)$. 
We will refer to the integers $\{ \hbox{1},\ldots ,\left| A \right| \}$ as the {\twelvebx states} of $X_A^ - $. We say 
that $X_A^ - $ is an {\twelvebx irreducible Markov shift over $A$} if $A$ is irreducible. Note that 
${X_A^ -}  \in End({\bf p})$ if and only if the entries in each row of $A$ are the 
components of ${\bf p}$, after deleting all zero entries. 

For each stochastic matrix $A$, consider the function $f_{A}: \{ \hbox{1},\ldots 
,\left| A \right| \}^{-{\bf N}} \to  \{ \hbox{1},\ldots 
,\left| A \right| \}$ defined by $f_{A}(x)=x_{ - 
1}$. By placing the discrete metric $d$ on $\{ \hbox{1},\ldots 
,\left| A \right| \}$, it is clear that $f_A$ is generating . Note that for the one-sided Markov shift $X^-_A$, the $(-{\hbox{1}})^{st}$ coordinates 
and the conditional probabilities of the inverse images of any point $x$ are completely determined by $x_{-1}$. If ${X_A^ -}  \in End({\bf p})$, 
we can thus choose a tree partition $K_{X_A^ - } $ of $X_A^ - $ with the 
following property: whenever $x_{-1}=y_{-1}$, then $x$ and $y$ generate the 
same $f_{A}$-tree name if the tree names are defined with respect to that 
partition. We will henceforth assume that for each ${X_A^ -}  \in End({\bf 
p})$, we choose $K_{X_A^ - } $ with this property. If $I$ is a state of $X_A^ - 
$, let $\tau _I^{f_A } $ be the common $f_{A}$-tree name generated by all 
points $y$ with $f_{A}(y)=I$.

We now wish to give examples of one-sided Markov shifts that are tvwB. 
Our first example is motivated by the well-known fact (Ornstein and 
Friedman) that a strongly mixing two-sided Markov shift is two-sided 
Bernoulli. Unfortunately, the one-sided analogue of this fact requires 
considerably more restrictions. Indeed, consider the one-sided Markov shift 
induced by the stochastic matrix:
$$
\left( \matrix{
 \hbox{2/3} & \hbox{1/3}  \cr 
\cr
\hbox{1/3} & \hbox{2/3} \cr}\right)$$
It is not difficult to see that ${X_A^ - } \in End(\textstyle{1 \over 
3},\textstyle{2 \over 3})$. However, as $X_A^ - $ is not tvwB (since 
the standard generator yields two tree names whose $\bar {t}_m $ distance is 
1 for all $m \in {\bf N}$), $X_A^ - $ is not isomorphic to $B^{ + }(\textstyle{1 \over 
3},\textstyle{2 \over 3})$ and hence cannot be one-sided Bernoulli. (Note 
that, however, the {\twelvebx two-sided} Markov shift over $A$ is isomorphic to the 
{\twelvebx two-sided} Bernoulli shift $B(\textstyle{1 \over 3},\textstyle{2 \over 
3})$.) Nonetheless, as the following shows, for one-sided Markov shifts 
which are uniformly $p$-to-1 endomorphisms, strong mixing does imply one-sided 
Bernoulli. 
\proclaim Proposition 4.1.1. Let $p \in {\bf N}$ and ${\bf p}=(\textstyle{1 \over 
p},...,\textstyle{1 \over p})$. Suppose A is a $N\times N$ stochastic matrix such that the one-sided Markov shift ${X_A^ -}$ is in $End({\bf p})$ and is strongly mixing. Then $X_A^ - 
$ is tvwB. \par  
Some observations will be helpful before we prove Proposition 4.1.1. Note 
that for every stochastic matrix $A$, we may associate a directed weighted 
graph to it. Specifically, we define $G(A)$ to be the graph with $\left| A\right|$ vertices 
identified by the integers $\hbox{1},\ldots ,\left| A\right|$ with a directed edge from $I$ to $J$ 
labeled with weight $w$ if $A_{IJ}=w$ if $w>\hbox{0}$. Note that if ${X_A^ -}  \in End({\bf 
p})$, the set of weights of the edges extending out from any vertex in $G(A)$ 
is the same. The critical observation which will be of use to us in the 
proofs of the ensuing propositions is the following: for any states $I$ and $J$, 
there exists a node $v$ of length $n$ such that $\tau _I^{f_A } (v)=J$ if and only if we have a path in $G(A)$ of length $n$ from vertices $I$ to $J$. 
\medbreak
\noindent {\twelvebx Proof (Proposition 4.1.1):} For brevity, we let $f$ denote $f_{A}$ 
throughout the proof. Since $X_A^ - $ is strongly mixing, there exists some integer 
$n$ such that $(A^{n})_{ij}>\hbox{0}$ for all $\hbox{1} \le i,j \le N$. We thus have some path in 
$G(A)$ of length $n$ between any two vertices in $G(A)$. By the observation made 
just prior to this proof, this implies in particular that there exists some state $J$ such 
that for all state $I$, $\tau _I^{f} (v) = J$ for some node $v$ of length $n$.  Hence, for all $x \in X_A^ - $, $\tau _x^{f} (v) = J$ for some node $v$ of 
length $n$. 

For any pair of points $(x',y')$, let $B_{x',y'}$ 
be any tree automorphism in ${\cal A}_{n}$ such that $\tau _{x'}^{f} (v)=\tau 
_{y'}^{f} (B_{x',y'} (v))=J$ for some node $v$ of length $n$. Given a fixed pair of 
points $(x,y)$ in $X_A^ - $, we now build a tree automorphism $B$ $n$-levels at a 
time which makes $\bar {t}_m (\tau _x^f ,\tau _y^f )$ small for all large $m$. 
For $\hbox{0}< \left| v \right| \le n$, let $B(v)=B_{x,y}(v)$. Inductively, assume that $B(v)$ is 
defined for all $\left| v \right| \le sn$. For $sn<\left| v \right| \le (s+\hbox{1})n$, let $v=uv'$, 
where $\left| v' \right|=sn$ and $\hbox{0}<\left| u \right| \le n$. We then define $B(v)=uB(v')$ 
if $f(T_{v'} x)=f(T_{B(v')} y)$ and define $B(v)=B_{(T_{v'} x,T_{B(v')}
y)} (u)B(v')$ otherwise. 

Now, we have some node $u$ of length $n$ such that $f(T_u x)=f(T_{B(u)} y)$.  Since the $f$-tree name of $x$ 
depends only on $x(-\hbox{1})$ and we are extending by the identity automorphism in the subtree rooted at $u$, it follows that if $n\le t<\hbox{2}n$, the weights of the nodes $\left| v \right| =t$ such that $f(T_v x)=f(T_{B(v)} y)$ total to at least $\hbox{1}/p^{n}$.  For the nodes $\left| u \right| =n$ such that $f(T_u x) \ne f(T_{B(u)} y)$, $B$ is defined in such a way that we have some node $\left| v \right| =n$ such that $f(T_{vu} x)=f(T_{B(vu)} y)$.  Thus, if $\hbox{2}n \le t<\hbox{3}n$, the weights of the nodes $\left| v \right| =t$ with $f(T_v x)=f(T_{B(v)} y)$ sum to at least 
\vfill \eject
$${\hbox{1} \over p^n}(\hbox{1} - {\hbox{1} \over p^n})+{\hbox{1} \over p^n}.$$	 

Inductively, we see that in general, if $sn\le t<(s+\hbox{1})n$, the weights of the nodes $\left| v \right| =t$ with $f(T_v x)=f(T_{B(v)} y)$ sum to at least
$$\sum\limits_{i = 0}^{s - 1} 
{{\hbox{1} \over p^n}(\hbox{1} - {\hbox{1} \over p^n})^i}.$$
\noindent This sum approaches 1 as $s \to 
\infty $, independent of $x$ and $y$. Hence it follows that $X_A^ - $ is tree 
v.w.B. {\closeproof}
\medbreak
The next proposition gives an additional class of tvwB Markov shift. 

\proclaim Proposition 4.1.2. Let ${\bf  p}$ be a probability vector $($not necessarily uniform$)$. Suppose that $A$ is an irreducible 
$N\times N$ matrix such that ${X_A^ -}  \in End({\bf 
p})$ and for every pair of rows in $A$, we can find two identical nonzero entries in the same column, then $X_A^ -$  is tvwB. \par 
\noindent {\twelvebx Proof:} Once again, consider the graph $G(A)$ associated to $A$. Note that 
it suffices to find $n \in {\bf N}$ along with paths of length $n$ from each 
vertex ending at a common vertex such that the corresponding edges in the 
paths have equal weights. Indeed, this will allow us to conclude that for 
any $x$ and $y$ in $X_A^ - $, there exists a node $v$ of length $n$ such that ${\tau 
_x^{f_A }} (v)={\tau _y^{f_A }} (B(v))$ for some tree automorphism $B$, and we 
may argue as before to reach the conclusion. 

To construct the required paths, note that by assumption, there are edges 
extending from vertices 1 and 2 with the same weight $w_{1}$ ending at a 
common vertex, say $J_{1}$, in $G(A)$. We then choose any edge of weight 
$w_{1}$ extending from vertex 3. If this edge ends at vertex $J_{2}$, then 
again by assumption we have edges of equal weights from $J_{1}$ and $J_{2}$ 
which end at a common vertex. This allows us to extend the paths from each 
of vertices 1, 2 and 3 such that they all end at the same vertex and the 
corresponding edges in the paths have the same weight. A simple inductive 
argument enables us to construct the desired paths.  {\closeproof}
\medbreak
We end this section with a proposition which shows that we can decide 
whether a one-sided Markov shift is tvwB simply by checking tree 
names of a finite height (depending on the dimension of the stochastic 
matrix $A$). For a general ${\bf p}$-endomorphism, deciding whether it is tree 
v.w.B is obviously a more difficult problem.

Let ${X_A^ -}  \in 
End({\bf p})$ and consider its associated graph $G(A)$.  Note that every path in $G(A)$ ``sees" a sequence of weights by reading the weights attached to the edges from the start to the end of the path.  We now prove: 
\proclaim Proposition 4.1.3. Suppose that $A$ is an irreducible $N\times N$ matrix such that the Markov shift ${X_A^ -}  \in End({\bf p})$. Then $X_A^ - $ is tvwB if and only if there exist paths of common length $\le N^{3N}$ from each vertex in G(A) which see a common sequence of weights and end in the same vertex. \par
\noindent {\twelvebx Proof:} We note that every edge $e$ in the directed graph $G(A)$ can be represented by 
the triple $(s(e),t(e),w(e))$, where $s(e)$, $t(e)$ and $w(e)$ are the starting 
vertex, terminal vertex and weight of $e$. Since the number of different 
nonzero weights cannot exceed the number of vertices, we have at most 
$N^{3}$ different types of edges under this representation. 

Now, if $({X_A^ -},f_{A})$ is tvwB then we have paths $v_{J}$ from 
each vertex $J$ in $G(A)$ which end at a common vertex and see the same sequence of weights.  Indeed, by the tvwB condition, for any two vertices $I$ and 
$I'$, we must have paths $v_{I}$ and $v_{I'} $ with the desired property. 
If $I''$ is another vertex, then construct any path $v_{I''} $ from 
$I''$ which see the same sequence of weights as $v_{I}$. The three paths just 
constructed end in at most two distinct vertices so we may extend these 
three paths such that they end at a common vertex. An inductive argument 
gives us the required paths. 

Assume that the paths chosen have a common length $>N^{3N}$. For a path $u$, 
let $u(j)$ be the $j{\hbox{-}}th$ edge of $u$. Then there are at most $N^{3N}$ possible 
ordered $N$-tuple of edges $(v_{1}(j),\ldots ,v_{N}(j))$ for each $j$. Now, if 
the paths $v_{1},\ldots ,v_{N}$ have more than $N^{3N}$ edges, then there 
exist integers $j<j'$ such that 
$$
(v_{1}(j),\ldots ,v_{N}(j))=(v_{1}(j'),\ldots 
,v_{N}(j')).
$$

\noindent Hence, we may shorten the path $v_{I}$ by deleting the edges $v_{I}(k)$ for 
$j \le k<j'$, for each vertex $I$. Clearly, the new paths still have the 
desired property but they have a shorter length. Hence, we may continue to 
shorten the paths to have a common length $ \le N^{3N}$.

Conversely, if the asserted property holds, then for every pair of states $I$ 
and $J$ of $X_A^ - $, there is a node $u$ of length $ \le N^{3N}$ such that 
there exists a tree automorphism $B$ with ${\tau _I^{f_A }} (u)={\tau _J^{f_A }} 
(Bu)$. The same argument in the proof of Proposition 4.1.1 shows that $X_A^ 
- $ is tvwB. {\closeproof}
\medbreak
Since the tree names $\tau _I^{f_A } $ restricted to ${\cal T}_{N^{3N}}' $ 
over all states $I$ of $X_A^ - $ determine all paths of length at most 
$N^{3N}$ in $G(A)$, Proposition 4.1.4 shows that it suffices to look at tree 
names of that height to determine whether $X_A^ - $ is tvwB. It is 
worth mentioning that Ashley, Marcus and Tuncel [A,M,T] developed a general, 
though necessarily more complicated, algorithm for deciding whether any two 
one-sided Markov shifts are isomorphic. 
\medjump
\noindent {\section 4.2. A Generalization of [T,Id]}
\medjump
In this section, we shall consider a case of the well-known $[T,Id] $
transformations in the context of a general probability vector ${\bf 
p}=(p_{1},\ldots ,p_{n})$ and characterize those that are one-sided 
Bernoulli.

The $[T,Id]$ transformation can be described as follows. We consider the 
2-shift $B^{ + }({\bf p})=({\{ \hbox{0,1}\} }^{\whole},\mu,\sigma )$ with 
${\bf p}=(\textstyle{1 \over 2},\textstyle{1 \over 2})$, and a Lebesgue 
space $(Y,{\cal C},\nu )$. Suppose $T$ is an automorphism of $Y$, define the 
map $[T,Id]$ on the product space ${\{ \hbox{0,1} \}}^{\whole}\times Y$ with 
product measure $\mu \times \nu $ by
$$
[T,Id](x,g) = (\sigma x,T^{x(0)}g).
$$
\noindent It is not difficult to show that $[T,Id]$ is 
measure preserving and the resulting m.p.s. $( {\{ \hbox{0},\hbox{1}\}}^{\whole} \times 
Y,\mu \times \nu ,[T,Id])$ is a ${\bf p}$-endomorphism. 

A special case of $[T,Id]$ occurs when $Y$ is just the circle represented as 
[0,1) and $T=R_{\alpha }$ is a rotation on $Y$ by an irrational $\alpha $. 
Hoffman and Rudolph [H,R] showed that this particular $[T,Id]$ system (along with other isometric extensions of the uniformly $p$-to-1 endomorphisms with certain 
properties) are all tvwB and hence one-sided Bernoulli.

Let us now extend the $[T,Id]$ system to the situation 
when ${\bf p}$ is an arbitrary finite probability vector. Consider a compact abelian metrizable group $G$ with a translation invariant metric  $d'$ and Haar measure $\nu$ (i.e. $d'(hg,h'g)=d'(h,h')$ for all $h$, $h'$ and $g \in G$).  Let $d$ be the 
discrete metric on $\{ \hbox{1},\ldots,n\}$. Define a metric $D$ on the set 
$R=\{ \hbox{1},\ldots ,n \}\times G$ by
$$
D((x_{1},g_{1}),(x_{2},g_{2}))={\hbox{1} \over \hbox{2}} d(x_{1},x_{2})+ {\hbox{1} \over \hbox{2}} d' (g_{1},g_{2}).
$$
Let ${\bf p}=(p_{1},\ldots ,p_{n})$ and let $B=\{ \hbox{1},\ldots ,n\}^{\whole}$. As usual, let $B^{ + }({\bf p})=(B,\mu ,\sigma )$ be the one-sided 
Bernoulli shift such that state $j$ has weight $p_{j}$, for 1 $\le j \le n$. 
For each state $j$, we associate to it some element $g _{j}$ in $G$.  Consider the transformation $S$ on $B\times G$ defined by $S(x,g)=(\sigma 
x,g g_{j})$ if $x_0=j$. Let $\lambda $ be the product measure $\mu \times \nu$, then 
note that $(B\times G,\lambda ,S) \in End({\bf p})$. 

Consider the function $f:B\times G \to R$ defined by $f(x,g)=(x_0,g)$. 
Note that $f$ is generating. Let $f_1$ and $f_2$ denote the component 
functions of $f$, i.e. $f_{1}(x,g)=x_0$ and $f_{2}(x,g)=g$. Clearly 
$f_{1}:B\times G \to \{ \hbox{1},\ldots ,n\}$ defines a tree partition 
of $(B\times G,\lambda ,S)$. We may thus use $f_{1}$ to define a set 
of partial inverses $S_{v}$ for each node $v$ of the ${\bf p}$-tree ${\cal T}$ as 
described in Chapter 1. Clearly, we have 
$f_{1}(S_{v}(x,g))=f_{1}(S_{v}(x',g'))$ for all nodes $v$ of length $ \ge $1 and for all $(x,g)$ and $(x',g')$ in $B\times G$.
\proclaim Proposition 4.2.1. $(B\times G, \lambda,S)$ is tvwB if there exist $i \ne j$ with $p_{i}=p_{j}$ such that $g 
_{j} g _{i}^{-1}$ has dense orbit in $G$. \par
\noindent {\twelvebx Proof:} We first assume the stated condition and prove that $((B \times G,\lambda ,S),f)$ is tvwB. With no loss of generality, we may 
assume that $p_{1}=p_{2 }$ and $g _{2}g _{1}^{-1}$ has dense orbit. Given $\varepsilon >\hbox{0}$, we note that there exists some $s \in 
{\bf N}$ such that for any $h$ and $h'$ in $G$, we can find 
an integer $\hbox{0}<r<s$ such that $d'(h g_{1}^{-r}, h' g _{2}^{-r})<\varepsilon$. To see this,  
partition $G$ into sets of diameter $< \varepsilon / \hbox{3}$. On 
each set, pick an arbitrary element in it. Say that $y_{1},\ldots ,y_{k}$ 
are the elements picked. Since $g _{2} g_{1}^{-1}$ has dense orbits and $G$ is abelian, for any pair $(y_{i}, y_{j})$ in $G$, there exists an integer 
$r> \hbox{0}$ such that $d'(y_{i}g_{1}^{-r}, y_{j}g _{2}^{-r})<\varepsilon / \hbox{3}$. We may then choose $s$ to be 
larger than all these $r$'s. Thus, for any two elements $h$ and $h'$ in $G$, we have some integer $\hbox{0}<r < s$ such that $d'(h g_{1}^{-r}, h'g _{2}^{-r})<\varepsilon. $

For any pair of points $x$ and $y$ in $B\times G$, by the above paragraph, there exists some tree automorphism $\bar A$ and some node $v$ of length $r<s$ such 
that
$$
d'(f_2 (S_v x),f_2 (S_{{\bar A}v} y)) < \varepsilon .
$$
\noindent (In fact, the inequality can be met by considering the node $v=(\hbox{1},\ldots,\hbox{1})$ of length $r$ and letting $\bar A$ be any tree automorphism such that ${\bar A}v=(\hbox{2},\ldots,\hbox{2})$.)  Moreover, for any node $v'$, we have 
$$
d'(f_2 (S_{v'v} x),f_2 (S_{v'{\bar A}v} y)) < \varepsilon .
$$
\noindent If $v'$ is not the root node, then as we have observed,
$$
d(f_1 (S_{v'v} x),f_1 (S_{v'{\bar A}v} y)) = \hbox{0}.
$$
\noindent  Hence, we have some node $u$ of length $s$ and a tree automorphism $A$ such that 
$$D(f(S_u x),f(S_{Au} y)) < \varepsilon .$$
 If $p_{n}$ is the smallest 
component in the probability vector ${\bf p}$, we can imitate the proof of 
Proposition 4.1.1 to construct the required tree automorphism $A$ $s$-levels at 
a time such that whenever $sN \le t<s(N+\hbox{1})$, the nodes $v$ of length $t$ for 
which $D(f(S_v x),f(S_{Av} y)) < \varepsilon $ have total weights at least 
$$\sum\limits_{i = 0}^{N - 1} {(p_n )^s(\hbox{1} - (p_n )^s)^i}. $$
 From this, we see 
that $(B\times G,\lambda ,S)$ is tvwB.  {\closeproof}

In the case when $G$ is the circle [0,1) with Lebesgue measure, we can improve Proposition 4.2.1 as follows:

\proclaim Proposition 4.2.2. $(B\times [ \hbox{0},\hbox{1}), \lambda,S)$ is tvwB if and only if there exist $i \ne j$ with $p_{i}=p_{j}$ such that $g 
_{j} g _{i}^{-1}$ is irrational. \par
\noindent {\twelvebx Proof:} Since the irrationals have dense orbits, the fact that $(B\times \hbox{[0,1)}, \lambda,S)$ is tvwB given the stated condition is just a special case of Proposition 4.2.1.
Conversely, suppose the stated condition is false but $(B\times \hbox{[0,1)},f)$ is 
tvwB. We then have some positive integer $N \ge \hbox{2}$ such that $g _{j} g _{i}^{-1}$ is some integral multiple of $\textstyle{1 \over 
N}$ whenever $p_{j}=p_{i}$. Suppose $x$ and $y$ are points in $B\times \hbox{[0,1)}$ 
with the property that $f_{2}(x)=t$ and $f_{2}(y)=t+\textstyle{1 \over 
{2N}}$ for some $t \in \hbox{[0,1)}$. Then by our choice of the rotation factors 
$g_{j}$, for every tree automorphism $A$ and node $v$, if 
$f_{2}(S_{v}x)=t'$, then $f_{2}(S_{Av}y)=t' + \textstyle{1 \over {2N}}+\textstyle{k \over N}$ for some integer $k$. Hence, $d'(f_{2}(S_{v}x),f_{2}(S_{Av}y))$ is bounded away from zero and so 
there exists some $\beta >\hbox{0}$ such that for all $n$, 
$$
\bar {t}_n (\tau _x^f ,\tau _y^f ) \ge \beta 
$$
for all pair of points $(x,y)$ such that $f_{2}(x)=t$ and 
$f_{2}(y)=t+\textstyle{1 \over {2N}}$ for some $t \in$ [0,1). Now, as any 
large set must contain a pair of such points, this contradicts the tvwB condition and completes the proof. {\closeproof}
\vfill \eject

$ $
\vskip 5.90cm
\noindent {\chapter Chapter 5:} 
\medjump
\noindent {\chapter Finite Group Extensions of One-Sided}
\medskip
\noindent {\chapter Bernoulli Shifts}
\bigjump
\twelverm
In this chapter, we consider Lebesgue spaces with certain semigroup actions defined on them.  The main result of this chapter is that for any two such spaces with semigroup actions $\{ T_g \}_{g \in G}$ and $\{ S_g \}_{g \in G}$, there is an isomorphism $\phi$ such that $\phi \circ T_g = S_g \circ \phi$ for all $g \in G$.  
\medjump
\noindent {\section 5.1. An Isomorphism Theorem on TvwB G-extensions}
\medjump
Let $G'$ be a metrizable semigroup. For an arbitrary 
nonatomic Lebesgue space $(X,{\cal B},\mu )$, we define a {\twelvebx $G'$-action} on $(X,{\cal B},\mu )$ as a collection of measure preserving 
endomorphisms $\{ T_{g} \}$ with the property that $T_{g}T_{g'}=T_{gg'}$ for all $g, g' \in G'$, and the map $\pi 
:G'\times X \to X$ defined by $\pi (g,x)=T_{g}x$ is 
measurable. If $(Y,{\cal C},\nu )$ is another nonatomic Lebesgue space with 
a $G'$-action $\{ S_{g} \}$ defined on it, then we say that 
$((X,{\cal B},\mu )$,$\{ T_{g}\} )$ and $((Y,{\cal C},\nu 
),\{ S_{g}\} )$ are {\twelvebx $G'$-isomorphic} if there exists a measure preserving 
bijection $\phi :X \to Y$ such that $\phi T_{g}=S_{g}\phi $ for all 
$g \in G'$. 

Consider a finite group $G$.  Note that ${\whole}\times G$ is a 
semigroup with operation defined by $(n,g) \cdot (n',g')=(n+n',g'g)$. For a probability vector ${\bf p}=(p_1,\ldots,p_s)$, we say that 
a ${\whole}\times G$-action $\{ T_{(n,g)} \}$ on a nonatomic Lebesgue 
space $(X,{\cal B},\mu )$ is a {\twelvebx $({\whole}\times G,{\bf p})$-action} if the m.p.s. $(X,{\cal B},\mu ,T_{(1,e)})$ is a ${\bf p}$-endomorphism and $T_{(0,g)}$ is an automorphism of $(X,{\cal B},\mu )$ with no fixed points unless $g=e$. The object of this chapter is to prove the following theorem.

\proclaim Theorem 5.1.1. Let $(X,{\cal B},\mu )$ and $(Y,{\cal C},\nu )$ be two nonatomic Lebesgue spaces with $(\whole \times 
G,{\bf p})$-actions  $\{T_{(n,g)}\}$ and $\{S_{(n,g)}\}$  defined on them. If the systems ${\bf X}=(X,{\cal B},\mu 
,T_{(1,e)})$ and ${\bf Y}=(Y,{\cal C},\nu ,S_{(1,e)})$ are both tvwB, then $((X,{\cal B},\mu 
),\{ T_{(n,g)} \} )$ and $((Y,{\cal C},\nu ), \{ S_{(n,g)}\} )$ are $\whole 
\times G$-isomorphic. \par

A natural example of a space with a $(\whole \times G,{\bf p})$-action is given by a {\twelvebx $G$-extension} of $B^+({\bf p})$, which we will now define.  Let $\nu$ be the uniform probability measure on the finite group $G$ equipped with the discrete topology. Let $I$ denote the set $\{ \hbox{1},\ldots ,s\}$. As we have seen, by assigning weight $p_j$ to $j$, we may represent the one-sided Bernoulli shift $B^{ + }({\bf p})$ as $(I^{\whole},{\cal B},\mu ,T)$, where 
$T$ is the shift map and $\mu $ is the product measure $\{ p_{1},\ldots 
,p_{s}\}^{\whole}$ on the Borel sigma-algebra ${\cal B}$. Let 
$P:I^{\whole} \to I$ be the time zero partition defined by $P(x)=x_0$.  

Consider the product space $I^{\whole}\times G$ 
with the product measure $\mu \times \nu $ on its product sigma-algebra 
$\cal P$. For a measurable partition $\varphi :I^{\whole} \to G$, define 
the map $T^{\varphi }:I^{\whole}\times G \to I^{\whole}\times G$ 
by $T^{\varphi }(x,g)=(Tx,\varphi (x)g)$.  We will refer to such partitions $\varphi$ as a {\twelvebx cocycle}. It is easy to verify that 
$T^{\varphi }$ is measure preserving. We define the $(G,\varphi)$-{\twelvebx extension} of $B^{ + }({\bf p})$, denoted as $B^{ + }({\bf p})^{\varphi }$, to be the m.p.s. $(I^{\whole 
}\times G$,$\cal P$,$\mu \times \nu $,$T^{\varphi })$. More generally, 
we say that a m.p.s. is a $G$-{\twelvebx extension} of $B^{ + }({\bf  p})$ if it is the $(G,\varphi 
)$-extension of $B^{ + }({\bf p})$ for some cocycle $\varphi :I^{\whole} \to G$. 
It is clear that any $G$-extension of $B^{ + }({\bf p})$ is a ${\bf 
p}$-endomorphism. Define the partition $P':I^{\whole} \times G \to 
I\times G$ by $P'(x,g)=(P(x),g)$. Let $d$ be the discrete metric on 
$I\times G$.  Clearly, $(I\times G,d)$ is a compact metric space and $P'$ is 
generating.

Given a $G$-extension $B^{ + }({\bf p})^{\varphi }$, note that each element 
$(n,g) \in {\whole}\times G$ corresponds to a measure preserving map obtained by composing the maps $\varphi _{(1,e)}$ and $\varphi_{(0,g')}$ given by
   
$$\varphi _{(1,e)} (x,g)=(Tx,\varphi (x) g)$$
and
$$\varphi_{(0,g')}(x,g)=(x,gg').$$

It is easy to verify that the maps $\varphi _{(n,g)} $ define a semigroup 
action, namely a ${\whole}\times G$-action, on $I^{\whole} \times G$.  (This is the reason why we defined $(n,g) \cdot (n',g')$ as $(n+n',g'g)$ instead of $(n+n',gg')$.) We say that the $G$-extensions $B^{ + 
}({\bf p})^{\varphi }$ and $B^{ + }({\bf p})^{\phi }$ are ${\whole}\times G$-{\twelvebx isomorphic} if there is a measure preserving bijection $\rho :I^{\whole} \times G \to I^{\whole} \times G$ such that $\rho \varphi _{(n,g)} =\phi _{(n,g)} 
\rho $ a.e. for all $(n,g) \in {\whole}\times G$. We say that $\rho $ is 
a ${\whole}\times G$-{\twelvebx isomorphism} from $B^{ + }({\bf p})^{\varphi }$ to $B^{ + }({\bf 
p})^{\phi }$. As $\rho \circ T^{\varphi }=T^{\phi } \circ \rho $, the ${\whole 
}\times G$-isomorphism $\rho $ is also an isomorphism from $B^{ + }({\bf 
p})^{\varphi }$ to $B^{ + }({\bf p})^{\phi }$ in the sense of our 
original definition in Chapter 1.

We begin the proof of Theorem 5.1.1 with some reductions.  First, note that if $(X,{\cal B},\mu)$ has a $(\whole \times 
G,{\bf p})$-action on it such that ${\bf X}=(X,{\cal B},\mu ,T_{(1,e)})$ is tvwB then we can represent it as a $G$-extension of $B^+({\bf p})$.  In fact, consider the space $X/G$ of $G$-orbits (so each point in $X/G$ is a set of the form $\{ T_{(0,g)}x \} _{g \in G}$).   $T_{(1,e)}$ projects to a measure preserving transformation $\bar T$ on $X/G$ such that $(X/G, \bar T)$ is a $\bf p$-endomorphism. Since $T_{(0,g)}$ has no fixed points unless $g=e$, the canonical projection is $\Pi_G : X \to X/G$ is a tree-adapted factor map.  Thus $(X/G, \bar T)$ is isomorphic to $B^+({\bf p})$ by Theorem 2.4.1 and Proposition 2.4.2.  Representing the space $X$ as $X/G \times G$, it is immediate that $(X,\{T_{(n,g)}\})$ is $\whole \times G$-isomorphic to $(I^{\whole} \times G, \{\phi_{(n,g)}\})$ for some cocycle $\phi : I^{\whole} \to G$.  Thus we may rephrase Theorem 5.1.1 as follows.

\proclaim Theorem 5.1.1$'$. If $B^{ + }({\bf p})^{\varphi }$ and $B^{ + }({\bf p})^{\phi 
}$ are tvwB, then they are $\whole \times G$-isomorphic. \par

Note that by our definition, all $G$-extensions of $B^{ + }({\bf p})$ are 
defined on the space $I^{\whole}\times G$, although different $G$-extensions 
have different dynamics. To emphasize the difference in dynamics, for any cocycle $\psi : I^{\whole} \to G $, we will 
henceforth use $I_{\psi}$ to denote the product space $I^{\whole 
} \times G$ of the m.p.s. $B^{ + }({\bf p})^{\psi }$ and we will denote 
the generating partition $P'$ of $I_{\psi}$ by $P^{\psi }$. 

Let $\Pi :I^{\whole}\times G \to I^{\whole}$ be the canonical 
projection. For $x \in I^{\whole}$ and a subset $A \subseteq I^{\whole}$, let the {\twelvebx  fiber over $x$} denote the set $\Pi ^{ - 1}(x)$ and the {\twelvebx fiber over $A$} denote the set $\Pi 
^{ - 1}(A)$. Note that for any ${\whole}\times G$-isomorphism $\rho 
:I_{\varphi } \to I_{\phi }$, since $\rho \varphi _{(0,g)} =\phi 
_{(0,g)} \rho $ for all $g \in G$, $\rho $ maps fibres to fibres, i.e. for 
$x \in I^{\whole}$, $\rho (\Pi ^{-1}(x))=\Pi ^{-1}(y)$ for 
some $y \in I^{\whole}$. Moreover, it is easily seen that the fibres are 
mapped to one another by a group rotation in the sense that for each $x \in 
I^{\whole}$, there is $g \in G$ and $x' \in I^{\whole}$ such 
that $\rho (x,g')=(x',gg')$ for all $g' \in 
G$.

Conversely, given measurable partitions $\varphi :I^{\whole} \to G$ and 
$\phi :I^{\whole} \to G$, an isomorphism $\rho _{1}:B^{ + }({\bf 
p}) \to B^{ + }({\bf p})$, and a measurable partition $\theta 
:I^{\whole} \to G$ such that
$$
\theta(Tx)\varphi (x) = \phi (\rho _1 (x))\theta (x),
$$
\noindent then it is direct that the map $\rho :I_{\varphi } \to I_{\phi }$ 
defined by $\rho (x,g)=(\rho _{1}(x),\theta(x)g)$ defines a ${\whole}\times G$-isomorphism from $B^{ + }({\bf p})^{\varphi }$ to $B^{ + 
}({\bf p})^{\phi }$. Thus, proving that $B^{ + }({\bf p})^{\varphi }$ 
and $B^{ + }({\bf p})^{\phi }$ are ${\whole}\times G$-isomorphic amounts 
to constructing the maps $\rho _{1}$ and $\theta$.

The basic approach to the proof of Theorem 5.1.1$'$ is the same as that of 
Theorem 2.4.1, with some necessary additions. We will focus primarily on 
these additions and refer the reader to the various propositions in Chapter 
2 when the proofs are exactly the same.

Before we embark on the proof, some remarks on notations are in order. 
Rather than considering functions taking values in the compact metric space 
[0,1] as in Chapter 2, {\twelvebx we will assume in this chapter that all the 
functions take values in the finite set $I\times G$ (i.e. are partitions), 
with the discrete metric $d$}, unless otherwise specified. To emphasize this 
fact, we will denote the functions by capital letters such as $Q$ and $R$ 
instead of small letters $g$ and $h$. Among the functions which range in 
$I\times G$, we shall primarily be interested in those functions $Q:I^{\whole}\times G \to I \times G$ for which there exists a pair of measurable 
partitions $Q_{1}:I^{\whole} \to I$ and $\xi :I^{\whole} \to G$ such 
that 
$$
Q(x,g)=(Q_{1}(x),\xi (x)g)
$$

\noindent
for all $(x,g) \in I^{\whole}\times G$ and $Q_{1}$ is tree-adapted. We 
say that $Q$ is a $G$-{\twelvebx map} if it has the above property. We can thus think of each $G$-map 
$Q$ as being defined by an ordered pair of functions $(Q_{1},\xi )$ of the 
above form. We will write $Q=(Q_{1},\xi )$ if $Q$ is defined by 
$(Q_{1},\xi )$. Note that the generating partition $P'$ is a $G$-map. If 
$Q$ is a $G$-map, then for each $x \in I^{\whole}$, $Q$ assigns the same first 
coordinate to any two distinct points in the fiber over $x$ but distinct group 
coordinates. Thus, a $G$-map takes every fiber onto a set of the form 
$\{ j\} \times G$ for some $j \in I$. 

For the purpose of constructing tree names, it will be convenient to choose 
a ``canonical'' set of partial inverses for $G$-extensions as follows. Recall 
that each node $v$ of the ${\bf p}$-ary tree ${\cal T}$ is a finite sequence of 
integers in $I=\{ \hbox{1},\ldots s \}$. Consider any $G$-extension $B^{ + }({\bf 
p})^{\varphi }$. For each $z \in I^{\whole}$, let $vz$ be the point in 
$I^{\whole}$ obtained by concatenating $v$ to the left of $z$. If $\left| v \right| =t$, then 
for each $(x,g)$ in $I_{\varphi }$, we define the partial inverse $T_v^\varphi 
$ by setting $T_v^\varphi (x,g)$ to be the unique element in $(T^{\varphi 
})^{-t}(x,g)$ whose first coordinate is $vx$. 

Given a function $Q$ on $B^{ + 
}({\bf p})^{\varphi }$, we may then define the partitions $Q^{N\tau 
}:I_{\varphi } \to (I\times G)^{N\tau }$, $Q^{N\nabla }:I_{\varphi 
} \to (I\times G)^{N\nabla }$ and $Q^{N + }:I_{\varphi } \to 
(I\times G)^{N}$ exactly as before. Moreover, we may define the $\bar {t}$ 
distance between any two processes $(B^{ + }({\bf p})^{\varphi },Q)$ and $(B^{ 
+ }({\bf p})^{\phi },R)$ for functions $Q$ and $R$ on the respective spaces as 
described in Chapter 2. Proposition 2.4.6 translates into the following.
\proclaim Proposition 5.1.2. Let $B^{ + }({\bf p})^{\varphi }$ be tvwB. Suppose $Q$ is a function on $I_{\varphi }$. For all $\varepsilon 
>\hbox{0}$, there exist $\delta $ and $N$ with the following property: for any G-extension $B^{ + }({\bf p})^{\phi }$ and function $R$ on $I_{\phi }$, if $R^{N\nabla }\mathop 
\sim \limits^\delta Q^{N\nabla }$, then $\bar {t}((B^{ + }({\bf p})^{\varphi 
},Q),(B^{ + }({\bf p})^{\phi },R))<\varepsilon$. {\closeproof} \par 
\proclaim Proposition 5.1.3. Consider $G$-extensions $B^{ + }({\bf p})^{\varphi }$ and $B^{ + }({\bf 
p})^{\phi }$, along with $G$-maps $Q$ on $I_{\varphi }$ and $R$ on $I_{\phi }$.  Then for all $g \in G$ and all $v \in {\cal T}'$, 
$$
d(Q(T_v^\varphi (x,g)),R(T_{v}^\phi ({x}',g))) = d(Q(T_v^\varphi 
(x,e)),R(T_{v}^\phi ({x}',e))).
$$ 
\par
\noindent {\twelvebx Proof:} This follows immediately from the fact that as $Q$ is a $G$-map, if 
$Q(T_v^\varphi (x,e))=(j,g')$, then $Q(T_v^\varphi 
(x,g))=(j,g'g)$, and likewise for the $G$-map $R$. {\closeproof}
\medbreak
We are now ready to prove the analogue of the perturbation lemma, 
Proposition 2.4.7. For $(j,g) \in I\times G$ and $\bar {g} 
\in G$, set $(j,g) \cdot \bar {g} = (j,g\bar {g})$. For any $G$-extension $B^{ 
+ }({\bf p})^{\varphi }$, a measurable set $B$ in $I_{\varphi }$ and $g \in 
G$, let $B_{g}$ denote the set $B \cap (I^{\whole}\times \{ g \} )$ (so 
$B_{g}$ is the subset of $B$ with group coordinate $g$).  In the following proposition and its proof, we will use the notation $O_{i}(\delta)$, $i=\hbox{1},\hbox{2},\ldots$  to denote real-valued functions of $\delta$ such that $\hbox{lim}_{\delta \to 0}O_{i}(\delta)=\hbox{0}$. 

\proclaim Proposition 5.1.4. Consider an ergodic $G$-extension $B^{ + }({\bf p})^{\varphi }$. Then there exists some $O_{1}(\delta)$ with the following property: For any ergodic $G$-extension $B^{ + }({\bf 
p})^{\phi }$ and any $G$-map $Q$ on $I_{\phi }$, if $\bar {t}((B^{ + }({\bf p})^{\varphi 
},P^{\varphi }),(B^{ + }({\bf p})^{\phi },Q))<\delta $, then for all $\varepsilon 
$ and $N$, there exists a $G$-map $Q'$ on $I_{\phi }$ such that $\left| Q-Q' \right|<O_{1}(\delta )$ and $(P^{\varphi })^{N\nabla}\mathop \sim \limits^\varepsilon Q'^{N\nabla }$. \par
\noindent {\twelvebx Proof:} We proceed along the lines of Proposition 2.4.7 except that we 
cannot use only a single tree name to define a function on the chosen tree 
Rokhlin tower as we did in Chapter 2. For $\eta >\hbox{0}$ to be specified later, 
choose $M \in {\bf N}$ such that
\itemitem{a)} $N/M<\eta / \hbox{2}$ 
\itemitem{b)} there exists a set $S \subseteq I^{\whole}$ with $\mu (S)>\hbox{1}-\eta $ 
such that whenever $x \in S$, $(x,g)$ is $\eta$,$(M-N)$-generic for 
$(P^{\varphi })^{N\nabla }$ for all $g \in G$
\itemitem{c)} $\int {\bar {t}_M } (\tau _{x,g}^{P^\varphi } ,\tau _{y,h}^Q )\,d(\mu \times 
\nu )(x,g)d(\mu \times \nu )(y,h) < \delta$.

From b) and c), if $\eta $ is sufficiently small, we have some $x \in S$ 
such that there exists a set $S' \subseteq I^{\whole}$ with $\mu 
(S')> \hbox{1}-O_{2}(\delta )$ and $\bar {t}_M (\tau 
_{x,e}^{P^\varphi } ,\tau _{z,e}^Q )<O_{2}(\delta )$ whenever $z \in 
S'$. Using the Strong Tree Rokhlin Lemma, build a $\eta / \hbox{2}$-tree Rokhlin 
tower ${\bf M}$ of height $M+\hbox{1}$ in $B^{ + }({\bf p})$ such that at least $\hbox{1} - 
O_{2}(\delta )$ fraction of the base $B$ is in $S'$. Note that ${\bf 
M}'=\Pi ^{-1}({\bf M})$ is also a tree Rokhlin tower in 
$B^{ + }({\bf p})^{\phi }$ with base $B' = \Pi ^{ - 1}(B)$. 
For each $(y,e) \in B'_e $, we define $Q'$ by laying the tree name 
$\tau _{x,e}^{P^\varphi } $ on $\{ (T^{\phi })^{ - i} (y,e) \mid {\hbox{1} \le i 
\le M} \}$ via a suitable tree automorphism $A$ which optimizes $\bar {t}_M (\tau 
_{x,e}^{P^\varphi } ,\tau _{y,e}^Q )$ as we did in Proposition 2.4.7. This 
defines $Q'$ on the column over $B'_e $. We then extend $Q'$ 
to points in the column over ${B}'_g $ for each $g \in G$ by defining 
$Q'(T_v^\phi (y,g))=Q'(T_v^\phi (y,e)) \cdot g$. Note that restricted 
to ${\cup \bf M}' \backslash B' $, $Q'$ is a $G$-map. By our 
construction, ${\cup \bf M}'  \backslash B'$ is a union of complete 
fibers, i.e. ${\cup \bf M}' \backslash B' =\Pi ^{-1}(E)$ for some 
set $E$ in $I^{\whole}$. We may thus extend $Q'$ to the rest of the 
space such that it remains a $G$-map by defining $Q'(x,g)=(P(x),g)$ on 
${\cup \bf M}' \backslash B'$. 

By Proposition 5.1.3 and our definition of $Q'$, whenever $x' \in S' \cap B$, we have for each $g \in G$,
$$
\eqalign{
  {\hbox{1} \over M}\sum\limits_{0 < \left| v \right| \le M} {w_v  d(Q(T_v^\varphi  (x',g)),Q'(T_v^\phi  (x',g)))} 
   &= {\hbox{1} \over M}\sum\limits_{0 < \left| v \right| \le M} {w_v  d(Q(T_v^\varphi  (x',e)),Q'(T_v^\phi  (x',e)))}  \cr 
   &= \bar t_M (\tau _{x,e}^{P^\varphi  } ,\tau _{x',e}^Q ) \cr
&< O_2 (\delta ) \cr} 
$$
\noindent Thus, since $\mu_{B} (S')> \hbox{1}-O_{2}(\delta )$, we have
$$
\int\limits_{{\cup \bf M}' \backslash B'} {d(Q(z,g),Q'(z,g))\,d\mu (z)d\nu (g)} < O_3 
(\delta ).
$$
Since $\mu \times \nu ({\cup \bf M}' \backslash B')>\hbox{1}-\delta $ for sufficiently small $\eta $, we have $\left|Q-Q'\right| 
<O_{1}(\delta )$. 

To prove that $(P^{\varphi })^{N\nabla }\mathop \sim \limits^\varepsilon 
Q'^{N\nabla }$, note that by our construction of $Q'$, we 
have $Q'^{M\nabla }(x',g)=(P^{\varphi })^{M\nabla 
}(x,g)$ for all $(x',g) \in B'$. By a) and b), since $x \in 
S$, the same reasoning in Proposition 2.4.7 shows that the required 
tree distributions differ by less than $\varepsilon $. {\closeproof} 
\medbreak
The following is an immediate consequence of the proof of Proposition 5.1.4.

\proclaim Proposition 5.1.5. Consider ergodic $G$-extensions $B^{ + }({\bf p})^{\varphi }$ and $B^{ + }({\bf 
p})^{\phi }$. For every $\varepsilon $ and $N$, there exists a $G$-map $Q$ on $I_{\phi }$ such that $(P^{\varphi })^{N\nabla }\mathop \sim 
\limits^\varepsilon Q^{N\nabla }.$ {\closeproof} \par

The next proposition is the Strong Sinai's Theorem which is formally identical 
to Proposition 2.4.9. 

\proclaim Proposition 5.1.6 (Strong Sinai's Theorem). Suppose $B^{ + }({\bf p})^{\varphi 
}$ is tvwB. Given $\varepsilon >\hbox{0}$, there exist $\delta $ and $N$ with the following property: If $B^{ + }({\bf p})^{\phi }$ is an ergodic $G$-extension, $Q$ is a $G$-map on $B^{ + }({\bf p})^{\phi }$ with $(P^{\varphi })^{N\nabla }\mathop \sim \limits^\delta Q^{N\nabla 
}$, then there exists a $G$-map $Q'$ on $B^{ + }({\bf p})^{\phi }$ such that $dist(B^{ + }({\bf p})^{\phi },Q')=dist(B^{ + }({\bf p})^{\varphi },P^{\varphi })$ and $\left| Q-Q' \right| <\varepsilon .$ \par
\noindent {\twelvebx Proof:} Use the approach in Proposition 2.4.9 to choose $\delta $ and 
$\hbox{1}/N$ small enough such that we have a Cauchy sequence of $G$-maps $Q^{j}$ on 
$I_{\phi }$ converging to some map $Q'$ such that $\left|Q-Q'\right|<\varepsilon $ with $dist(B^{ + }({\bf p})^{\phi },Q')=dist(B^{ + }({\bf p})^{\varphi },P^{\varphi })$. It is easily 
checked that $Q'$ is also a $G$-map. {\closeproof}

Given the $G$-map $Q=(Q_1,\zeta)$ such that $dist(B^{ + }({\bf p})^{\phi 
},Q)=dist(B^{ + }({\bf p})^{\varphi },P^{\varphi })$, we can proceed as in Chapter 2 to construct a tree-adapted factor map $\rho :B^{ + }({\bf p})^{\phi } \to B^{ + }({\bf 
p})^{\varphi }$ defined by $\rho (x,g)=(Q_1^{\whole}(x),\zeta(Q_1^{\whole}(x))g)$.   Moreover, for any $G$-extension $B^{ + }({\bf p})^{\phi }$, the canonical projection 
$\Pi :I_{\phi } \to I^{\whole}$ is clearly a factor map from $B^{ + 
}({\bf p})^{\phi }$ to $B^{ + }({\bf p})$. Combining Propositions 5.1.5 and 
5.1.6, along with the fact that the partition $Q'$ constructed in 
Proposition 5.1.6 is a $G$-map, the following corollary is immediate.
\proclaim Corollary 5.1.7 (Sinai's Theorem). Suppose $B^{ + }({\bf p})^{\varphi }$ is tvwB and $B^{ + 
}({\bf p})^{\phi }$ is ergodic. Then there exists a $G$-map $Q$ on $B^{ + }({\bf p})^{\phi }$ such that $dist(B^{ + }({\bf 
p})^{\phi },Q)=dist(B^{ + }({\bf p})^{\varphi },P^{\varphi }).$ Moreover, the factor map $\rho 
:B^{ + }({\bf p})^{\phi } \to B^{ + }({\bf p})^{\varphi }$ as constructed above projects to a factor map $\rho 
_{\Pi }:B^{ + }({\bf p}) \to B^{ + }({\bf p})$  such that the following diagram commutes:
$$
\matrix{{B^ + ({\bf p})^\phi } & {\mathop \to \limits^{\rho } } & {B^ + 
({\bf p})^\varphi }\cr
\cr
 {\Pi \downarrow }  &   & { \downarrow \Pi } \cr
\cr
 {B^ + ({\bf p})}  & {\mathop \to \limits_{\rho _\Pi } }  & {B^ + ({\bf p})}. \cr }
$$
{\closeproof}

Our next goal is the copying lemma. Recall that for $B \subseteq I^{\whole}\times G$ and $g \in G$, we defined the set $B_{g}$ to be the subset of $B$ 
with group coordinate $g$. For 
$g \in G$ and a set $S \subseteq I^{\whole}\times G$, let
$$
S \cdot g=\{ (x,g'g) \in I^{\whole}\times G \mid (x,g') 
\in S \}.
$$
\proclaim Proposition 5.1.8 (Copying Lemma). For ergodic $G$-extensions $B^{ + }({\bf p})^{\phi }$ and $B^{ + }({\bf 
p})^{\varphi }$, let $Q$ be a $G$-map on $I_{\phi }$ such that $dist(B^{ + }({\bf p})^{\phi },Q)=dist(B^{ + 
}({\bf p})^{\varphi },P^{\varphi })$. 
Then for all $\varepsilon $ and $N$, we have a $G$-map $Q'$ on $I_{\varphi }$ such that 
$$
dist(Q' \vee P^\varphi )^{N\nabla }\mathop \sim \limits^\varepsilon 
dist(P^\phi \vee Q)^{N\nabla }. 
$$
\par
\noindent {\twelvebx Proof:} Let $M$ be chosen such that $N/M<\varepsilon / \hbox{2}$. Define a 
partition $\bar {\varphi }:I^{\whole } \to I\times G$ by $\bar {\varphi 
}(x)=(P(x),\varphi (x))$. Construct a $\varepsilon /\hbox{2}$-tree Rokhlin tower 
${\bf M}$ of height $M+\hbox{1}$ in $B^{ + }({\bf p})$ with base $B$ such that 
$dist(\bar {\varphi }^{M\nabla })=dist(\bar {\varphi }^{M\nabla }|B)$. 
Using the canonical projection $\Pi :I_{\varphi } \to I^{\whole }$, we 
may lift the partition $\bar {\varphi }^{M\nabla }$ and the tree Rokhlin 
tower ${\bf M}$ so that we may regard $\bar {\varphi }^{M\nabla }$ as 
being a partition of $I_{\varphi }$ and the tree Rokhlin tower ${\bf M}$ 
with base $B$ as being in $B^{ + }({\bf p})^{\varphi }$. Clearly, we still 
have $dist(\bar {\varphi }^{M\nabla })=dist(\bar {\varphi }^{M\nabla 
}|B)$. For any atom $\alpha  \in \bar {\varphi }^{M\nabla }$ and $g 
\in G$, note that all points in $\alpha _{g}$ have the same $P^{\varphi 
}$-$M$-tree name. Thus, $dist((P^{\varphi })^{M\nabla })=dist((P^{\varphi 
})^{M\nabla }|B)$. 

Consider the tree-adapted factor map $\rho :B^{ + }({\bf p})^{\phi } \to 
B^{ + }({\bf p})^{\varphi }$ as constructed in Sinai's Theorem. For any atom $\alpha  \in \bar 
{\varphi }^{M\nabla }$, note that if $(x,g) \in \alpha $, then 
$(P^{\varphi })^{M\nabla }(x,g)=Q^{M\nabla }(y,g')$ for any 
$(y,g') \in \rho ^{ - 1} (x,g)$ by Proposition 2.1.4. Thus 
$(P^{\varphi })^{M\nabla }(\alpha _{g})=Q^{M\nabla }(\rho ^{ - 1} 
(\alpha _g ))$.

For each atom $\alpha  \in \bar {\varphi }^{M\nabla }$, construct a 
partition $U_{\alpha }^e$ of $\alpha _{e} \cap B$ such that 
$$
dist(U_{\alpha }^e)=dist((P^{\phi } \vee 
Q)^{M\nabla }|\rho ^{ - 1} (\alpha _e )) \eqno(\hbox{1})
$$

\noindent We then define the required partition $Q'$ on the column of the tower over $B_{e}$ as follows. Fix 
an atom $\alpha  \in \bar {\varphi }^{M\nabla }$. By (1), we can 
define a bijective correspondence $\Lambda $ of the atoms of $U_{\alpha }^e$ 
and those of $(P^{\phi } \vee Q)^{M\nabla }$ with the same conditional 
measures. Consider an atom $\beta  \in U_{\alpha }^e$. As we have already noted, $(P^\varphi)^{M \nabla}(x,e)=Q^{M \nabla}(x',g')$ for any $(x,e) \in \beta$ and $(x',g') \in \Lambda (\beta)$. We may then define $Q'$ on the tower over $\beta$ by laying the $P ^ \phi $-$M$-tree name of some point $(x',g')$ in $\Lambda (\beta)$ via suitable tree automorphisms such that $(Q' \vee P^{\varphi })^{M\nabla }(x,e)=(P^{\phi } \vee 
Q)^{M\nabla }(x',g')$ for any $(x,e) \in \beta $. 
Repeating for each $\beta  \in U_{\alpha }^e$ and then for each $\alpha  
\in \bar {\varphi }^{M\nabla }$, we can then define $Q'$ on the 
tower over $B_{e}$. 

We now extend $Q'$ on the column of the tower over $B_{g}$ for all $g 
\in G$. To do this, note that if $(x,g')$ is in the column over $B_{e}$, 
then $(x,g'g)$ is in the column over $B_{g}$. Define $Q' 
(x,g'g)=Q'(x,g')\cdot g$. This defines 
$Q'$ on ${\cup \bf M}\backslash B$. It is clear that $Q'$ is a 
$G$-map restricted to ${\cup \bf M}\backslash B$. Since ${\cup \bf M}\backslash B$ 
is ${\cal B}$-measurable, we may then extend $Q'$ so that it is a $G$-map 
on the full space by setting $Q'(x,g)=(P(x),g)$.

We now prove that $dist(Q' \vee P^\varphi )^{N\nabla }\mathop \sim 
\limits^\varepsilon dist(P^{\phi} \vee Q)^{N\nabla }$. If $\kappa  \in (I \times 
G)^{M\tau }$ and $g \in G$, let $\kappa  \cdot g \in (I\times 
G)^{M\tau }$ be defined by $(\kappa  \cdot g)(v)=\kappa (v) 
\cdot g$. If $\bar {\kappa } \in (I\times G)^{M\nabla }$ and $\kappa $ 
is a representative in $\bar {\kappa }$, define $\bar {\kappa } \cdot g 
\in (I\times G)^{M\nabla }$ to be the equivalence class containing 
$\kappa  \cdot g$. For any atom $\alpha  \in \bar {\varphi 
}^{M\nabla }$ and $\beta_{\alpha}  \in U_{\alpha }^e$, note that if $(x,e) \in 
\beta_{\alpha} $ and $(x',g') \in \Lambda (\beta_{\alpha}  )$, then for 
each $g \in G$,
$$
\eqalign{
 (Q' \vee P^\varphi )^{M\nabla }(x,g) &= ((Q' \vee P^\varphi )^{M\nabla 
}(x,e))\cdot g \cr 
&= ((P^\phi \vee Q)^{M\nabla }(x',g')) \cdot g \cr 
&= (P^\phi \vee Q)^{M\nabla }(x',g'g). \cr }
$$
\noindent This implies that $\beta_{\alpha}   \cdot g$ and $\Lambda (\beta_{\alpha}  ) \cdot g$ have 
the same image under $(Q' \vee P^{\varphi })^{M\nabla }$ and 
$(P^{\phi } \vee Q)^{M\nabla }$ respectively. Moreover, since $Q$ is a 
$G$-map, $\rho ^{ - 1} (\alpha _g )=\rho ^{ - 1} (\alpha _e ) \cdot g$ and so
$$
\eqalign{
 \mu \times \nu (\Lambda (\beta_{\alpha}  ) \cdot g \cap \rho ^{ - 1} (\alpha _g )) 
&= \mu \times \nu (\Lambda (\beta_{\alpha}  ) \cdot g|\rho ^{ - 1} (\alpha _g )) 
\mu \times \nu (\rho ^{ - 1} (\alpha _g )) \cr 
&= \mu \times \nu (\Lambda (\beta_{\alpha}  )| \rho ^{ - 1} (\alpha _e ))
\mu \times \nu (\rho ^{ - 1} (\alpha _e )) \cr  
&= \mu \times \nu (\Lambda (\beta_{\alpha}  )|\rho ^{ 
- 1} (\alpha _e )) \mu \times \nu (\alpha_e ) \cr 
&= \mu \times \nu (\beta_{\alpha}  | \alpha _e \cap B) \mu \times \nu (\alpha 
_e | B) \cr 
& = \mu \times \nu (\beta_{\alpha}  | B) = \mu \times \nu (\beta_{\alpha}  \cdot g | B) \cr }
$$
\noindent As this holds for all $g \in G$, $\alpha \in {{\bar \varphi} ^{M\nabla}}$, and $\beta_\alpha \in U_{\alpha}^e$, we have 
$$
dist(P^{\phi } \vee Q)^{M\nabla }=dist((Q' \vee P^{\varphi 
})^{M\nabla }|B)
$$
\noindent Hence, by the same reasoning as in the proof of Proposition 2.4.12, we see 
that the distribution of $(Q' \vee P^{\varphi })^{N\nabla }$ on 
each level of the tower except the top $N$ levels equals $dist(P^{\phi } 
\vee Q)^{N\nabla }$. Since $N/M<\varepsilon /\hbox{2}$, the result follows. {\closeproof}
\medbreak

The rest of the proof now follows along the same lines as Theorem 2.4.1.  To begin, we have the following proposition which is formally identical to Proposition 2.1.13.  

\proclaim Proposition 5.1.9. Consider tvwB $G$-extensions $B^{ + }({\bf p})^{\phi }$  and $B^{ + }({\bf p})^{\varphi }$. Suppose $Q$ is a $G$-map on $I_{\phi }$  such that $dist(B^{ + }({\bf p})^{\phi },Q)=dist(B^{ + 
}({\bf p})^{\varphi }$, $P^{\varphi })$. Then for every $\eta > \hbox{0}$, $\varepsilon >\hbox{0}$, there exists a $G$-map $Q'$ on $I_{\phi }$ such that 
\itemitem {i)} $P^\phi \mathop \subset\limits^{\varepsilon} \sum (Q') $
\itemitem {ii)} $dist(B^{ + }({\bf p})^{\phi },Q')=dist(B^{ + }({\bf 
p})^{\varphi },P^{\varphi })$
\itemitem {iii)} $\left|Q-Q' \right|<\eta $. \par
\noindent {\closeproof} 
\medbreak
\noindent {\twelvebx Proof (Theorem 5.1.1$'$):} Imitating the proof of Theorem 2.4.1, we choose, 
via proposition 5.1.9, a Cauchy sequence of $G$-maps $\{U^{j}\}$ on 
$I_{\phi }$ converging to a $G$-map $U$ on $I_{\phi }$ such that $dist(B^{ + 
}({\bf p})^{\phi },U)=dist(B^{ + }({\bf p})^{\varphi },P^{\varphi 
})$ and $\sum (U)=\sum (P^{\phi })$. 

Let $U=(U_{1},\xi )$. Since $\sum(U)=\sum (P^{\phi })$, the factor 
map $\rho :B^{ + }({\bf p})^{\phi } \to B^{ + }({\bf 
p})^{\varphi }$ defined by $\rho (x,g)=(\bar {U}_1 (x),\xi (x)g)$ is an isomorphism, where $\bar {U}_1 (x)$ is the $U_{1}$-name 
of $x$. Since 
$\rho T^{\phi }=T^{\varphi }\rho $, it follows that 
$$
\varphi (\bar {U}_1 (x))\xi (x) = \xi (Tx)\phi (x)
$$
\noindent Obviously, $\bar {U}_1 $ defines an automorphism of $B^{ + }({\bf p})$.  Hence, $\rho $ defines a ${\whole}\times G$-isomorphism from $B^{ + 
}({\bf p})^{\varphi }$ to $B^{ + }({\bf p})^{\phi }$. {\closeproof} 
\medjump
\noindent {\section 5.2. Some Applications of Theorem 5.1.1}
\medjump
\noindent{\twelvebx Example 5.2.1.} As an application of Theorem 5.1.1, consider a finite group $G$ of order $\ge$ 3 with two distinct generators $h$ and $h'$ (thus $G \cong {\bf 
Z}/n{\bf Z}$, for $n \ge \hbox{3}$). Fix a probability vector ${\bf 
p}=(p_{1},\ldots ,p_{s})$ with $p_1 = p_2$.
Define a map $\varphi :I^{\whole} \to G$ by
$$
\varphi (z) = \cases{h & if $z_0 = \hbox{2}$;\cr \cr e & otherwise. \cr}
$$

Note that the $P^{\varphi }$-tree name of $(x,g)$ is completely determined by 
$P^{\varphi }(x,g)=(x_0,g)$. Let $\tau _{(j,g)}$ be the common $P^{\varphi }$-tree 
name of all points $(x,g)$ such that $x_0=j$. Using the ideas in \S 4.1, to 
prove that $B^{ + }({\bf p})^{\varphi }$ is tvwB, it suffices to 
show that for any two elements $(j,g)$ and $(j',g')$ in $I\times 
G$, we have some node $v$ and some tree automorphism $A$ such that $\tau 
_{(j,g)}(v) =\tau _{(j',g')}(Av)$. In fact, it is enough 
to show this in the case when $j=j'=\hbox{1}$.

Choose $k \in {\bf N}$ such that $g'=h^{k}g$. Then for the node 
$v=(\hbox{1},\hbox{1},\hbox{1},\ldots ,\hbox{1})$ of length $k+\hbox{1}$, $\tau _{(j,g)}(v)=(\hbox{1},g)$; for the 
node $u=(\hbox{1},\hbox{2},\hbox{2},\ldots ,\hbox{2})$ of length $k+\hbox{1}$, $\tau _{(j,g')}(u) = 
(\hbox{1},g)$. Since these two nodes can be matched by a tree automorphism, we have $\tau _{(j,g)}(v) = 
\tau _{(j',g')}(Av)$ for some tree automorphism $A$. Thus, 
$B^{ + }({\bf p})^{\varphi }$ is tvwB.

In the same way, for the map $\psi :I^{\whole} \to G$ defined by
$$
\psi (z) = \cases{h' & if $z_0 = \hbox{2}$;\cr \cr e & otherwise, \cr}
$$
the same argument shows that $B^+({\bf p}) ^{\psi}$ is tvwB.  Thus Theorem 5.1.1 shows that $B^{ + }({\bf p})^{\varphi }$ and $B^+({\bf p}) ^{\psi}$ are $\whole \times G$-isomorphic.  

Note that the two cocycles $\phi$ and $\psi$ are not cohomologous provided that $hh'^{-1}$ is a generator (if they were cohomologous, then a $\whole \times G$-isomorphism between the corresponding extensions exists trivially).   Indeed, if they were cohomologous, the $G$-extension obtained from the map $\phi - \psi$ defined by $(\phi - \psi) (x) = \phi (x) \psi (x)^{-1}$ would not even be ergodic, contradicting the above argument which shows that it must be tvwB and hence ergodic.  {\closeproof}
 \medbreak

The following proposition, as observed by del Junco, gives another interesting 
application of Theorem 5.1.1. We say that a tree-adapted factor map $\phi :X \to Y$ is 
{\twelvebx uniformly $p$-to-1} if the fiber over $a.a.y$ in $Y$ contains $p$ points each with 
conditional probability $\hbox{1}/p$.  One example of such a factor map is if ${\bf X}={\bf Y}$ is the one-sided (1/2,1/2)-Bernoulli shift and $\phi$ is the addition map defined by $\phi (x)_i=x_i+x_{i+1}$.  The following proposition essentially says that in the case when ${\bf p}=(\hbox{1}/\hbox{2},\hbox{1}/\hbox{2})$, then this is (up to automorphism) the only example.

\proclaim Proposition 5.2.2. Suppose \X $\in End({\bf 
p})$ is tvwB. Let \Y be a tree adapted factor of ${\bf X}$. For any pair of uniformly 2-to-1 tree-adapted factor maps $\Pi _{1}:{\bf X} \to {\bf 
Y}$ and $\Pi _{2}:{\bf X} \to {\bf Y}$, there exist automorphisms $\rho :{\bf X} \to {\bf 
X}$ and $\varphi :{\bf Y} \to {\bf Y}$ such that $\Pi _{2} \circ \rho =\varphi \circ \Pi_{1}. $ \par

\noindent {\twelvebx Proof:} We define two $(\whole \times {\bf Z}_{2},{\bf p})$-actions on 
$(X,{\cal B},\mu )$ as follows. Let $ {\bf Z}_{2}=\{ \hbox{0},\hbox{1} \}$. Define 
$T_{(1,0)}=T$, and $T_{(0,0)} $ to be the identity map. Define $T_{(0,1)} 
:X \to X$ by sending $x$ to the unique point $x'$ distinct from $x$ such 
that $\Pi _{1}(x')=\Pi _{1}(x)$. The fact that $\Pi _{1}$ is 
uniformly 2-to-1 implies that $T_{(0,1)}$ is an automorphism of $X$. Using 
$T_{(1,0)}$, $T_{(0,0)}$ and $T_{(0,1)} $, we can define $T_{(n,g)}$ by 
group composition. By tree-adaptedness, it is easy to check that 
$\{ T_{(n,g)} \}$ defines a $({\whole}\times {\bf Z}_2,{\bf p})$-action on $X$. In a 
similar way, using $\Pi _{2}$, we can define a second $({\whole}\times 
{\bf Z}_2,{\bf p})$-action $\{ S_{(n,g)} \}$ on $X$. 

By Theorem 5.2.2, we have a ${\whole}\times {\bf Z}_{2}$-isomorphism $\rho :(X,\{ T_{(n,g)} \}) \to (X,\{ S_{(n,g)} \})$. In particular, $\rho 
$ defines an automorphism of ${\bf X}$. Note that as $S_{(0,1)} \circ \rho 
=\rho  \circ T_{(0,1)}$, $\rho $ maps fibres of $\Pi _{1}$ to those 
of $\Pi _{2}$ and so the map $\varphi :{\bf Y} \to {\bf Y}$ defined by 
$\varphi =(\Pi _{2}) \circ \rho  \circ (\Pi _{1})^{ - 1}$ is 
well-defined and is clearly an isomorphism with the required property. {\closeproof}
\medbreak
\noindent {\twelvebx Remark 5.2.3.} It is reasonable to believe that Theorem 5.1.1 should also 
hold in the case when $G$ is a compact metrizable group instead of just a 
finite group. We were unable to prove this, although it seems that the same 
ideas should apply and is likely a technical extension. One possible 
approach would be to define a partition $Q_{N}$ on $G$ by dividing $G$ into 
subsets of diameter $ \le  \hbox{1}/ \hbox{2}^{N}$ for each $N$ and define tree 
distributions by a suitable discretization.  

Second, it would be interesting to give a joinings proof of this result, 
though the definition of one-sided joinings might need to be modified to achieve the required isomorphism. 
\vfill \eject

$ $
\vskip 5.90cm
\noindent {\chapter Chapter 6:} 
\medjump
\noindent {\chapter Some Open Problems}
\bigjump
\twelverm

This thesis has extended the investigation of the isomorphism problem of 
non-invertible endomorphisms initiated by Hoffman and Rudolph. Obviously, 
much is left unaddressed. In this chapter, we state some problems in this area which 
warrant further investigation.
\medjump
\noindent {\twelvebx Problem 1.} Consider a Lebesgue probability space $(X,{\cal 
B},\mu )$ with two commuting endomorphisms $T$ and $S$ on it, i.e. $TS=ST$ a.e. The basic model is the two-dimensional lattice $\{\hbox{0},\hbox{1} \}^{({\bf N}^{\ast})^2}$ with product measure $(\hbox{1/2},\hbox{1/2})^{({\bf N}^{\ast})^2}$, and with $T$ and $S$ being the shift maps in each direction.  Given another Lebesgue space $(Y,{\cal C},\nu)$ with two commuting endomorphisms $T'$ and $S'$ on it with the same property as $T$ and $S$, is there a reasonable way to decide when $(X,{\cal B},\mu)$ and $(Y,{\cal C},\nu)$ are ${\whole}\times {\whole}$-isomorphic, i.e. when can we find a measure preserving bijection 
$\phi :X \to Y$
such that $\phi T=T' \phi $ and $\phi S=S' \phi $ a.e.?
\medjump
\noindent {\twelvebx Problem 2.} So far, all endomorphisms \X 
considered are {\twelvebx homogeneous} in the sense that for a.a. $x$, $x$ has the same number of 
inverse images with the same set of conditional probabilities. Obviously, 
most endomorphisms do not share this property. Can one give a reasonable 
classification for a subset of these? If we have a fixed 
one-sided Markov shift which is not homogeneous, is there some reasonable 
criterion to ensure that a given endomorphism is isomorphic to it?
\medjump
\noindent {\twelvebx Problem 3.} Proposition 5.2.2 in particular implies that there is essentially only one tree-adapted uniformly 2-to-1 endomorphisms of $B^+ ({\bf p})$.  It would be interesting to classify tree-adapted uniformly 3-to-1 endomorphisms (and in general $p$-to-1 endomorphisms) of $B^+ ({\bf p})$.  The proof of Proposition 5.2.2 cannot be used since even for tree-adapted uniformly 3-to-1 maps of $B^+ ({\bf p})$, there is no canonical way of constructing a $\whole \times {\bf Z}_3$-action on the Bernoulli shift space.   

\vfill \eject
$ $
\vskip 1.50 cm
{\chapter References:}
\bigjump
\item{[A,M,T]} J. Ashley, B. Marcus, S. Tuncel.  The Classification of One-Sided Markov Chains.  {\twelvesl Ergodic Theory $\&$ Dynamical Systems} {\twelvebx 17} (1997), 269-295.
\medbreak
\item{[Bil]} P. Billingsley.  Ergodic Theory and Information.  {\twelvesl John Wiley $\&$ Sons, Inc.}, New York, 1965.
\medbreak
\item{[Bri]} J.Y. Briend.  La propri\'et\'e de Bernoulli pour les endomorphismes de ${\bf P}^k(\hbox{\field C})$, {\twelvesl Ergodic Theory $\&$ Dynamical Systems} {\twelvebx 22}, (2002), 323$-$327.
\medbreak
\item{[Fur]}	H. Furstenberg.  Recurrence in Ergodic Theory and Combinatorial Number Theory.  {\twelvesl Princeton University Press}, Princeton, New Jersey, 1981.
\medbreak
\item{[H,H]}	D. Heicklen, C. Hoffman.  Rational Maps are d-adic Bernoulli.  {\twelvesl Annals of Math.} {\twelvebx 156} (2002), 103$-$114.
\medbreak
\item{[H,R]}	C. Hoffman, D. Rudolph.  Uniform Endomorphisms which are Isomorphic to a Bernoulli Shift.  {\twelvesl Annals of Math.} {\twelvebx 156} (2002), 79$-$101.
\medbreak
\item{[J,J]} P. Jong, A. del Junco.  Endomorphisms Isomorphic to the Bernoulli p-Shift.  Submitted to {\twelvesl Studia Mathematica}.
\medbreak
\item{[Orn]}	D. Ornstein.  Bernoulli Shifts with the Same Entropy are Isomorphic.  {\twelvesl Advances in Math.} {\twelvebx 4} $($1970$)$, 337$-$352.
\medbreak
\item{[Pa1]}	W. Parry.  Entropy and Generators in Ergodic Theory.  {\twelvesl W. A. Benjamin, Inc.}, New York, 1969.
\medbreak
\item{[Pa2]}	W. Parry.  Topics in Ergodic Theory. {\twelvesl Cambridge University Press}, Cambridge, 1981.   
\medbreak
\item{[Pet]} 	K. Peterson.  Ergodic Theory.  {\twelvesl Cambridge University Press}, Cambridge, 1983.
\medbreak
\item{[Roy]} H. L. Royden.  Real Analysis, 3rd edition.  {\twelvesl Macmillan Publishing Company}, New York, 1967.
\medbreak
\item{[Rud]}	D. Rudolph.  Fundamentals of Measurable Dynamics: Ergodic Theory of Lebesgue Spaces.  {\twelvesl Oxford University Press}, New York, 1990. 
\medbreak
\item{[Ryd]} Ryder, Herbert John.  Combinatorial Mathematics.  The Carus Mathematical Monographs, $\#$14. {\twelvesl John Wiley \& Sons Inc}, New York, 1963.
\medbreak
\item{[Shi]} P. Shields.  The Theory of Bernoulli Shifts.  {\twelvesl University of Chicago Press}, Chicago and London, 1973.
\medbreak
\item{[Var]} 	S. Varadhan.  Ergodic Theory: A Seminar. {\twelvesl Courant Institute of Mathematical Sciences},  New York University, 1975, Chapters 1$-$6.  

\end